\documentclass[12pt]{amsart}
\usepackage{texdraw,multicol,rotating}

\usepackage{amsmath}
\usepackage{amssymb}

\input{txdtools.tex}
\theoremstyle{plain}
\newtheorem{thm}{Theorem}[section]
\newtheorem{prop}[thm]{Proposition}

\newtheorem{cor}[thm]{Corollary}
\theoremstyle{definition}
\newtheorem{defi}[thm]{Definition}

\newtheorem{example}[thm]{Example}
\theoremstyle{remark}

\numberwithin{equation}{section}

\newbox{\tmpa}
\newbox{\tmpb}

\DeclareMathOperator{\wt}{wt}
\newcommand{\nc}{\newcommand}
\nc{\Uq}{U_q}
\nc{\Z}{\mathbf{Z}}
\nc{\C}{\mathbf{C}}
\nc{\Q}{\mathbf{Q}}
\nc{\op}{\oplus}
\nc{\ot}{\otimes}
\nc{\pv}{P^{\vee}}
\nc{\ali}{\alpha_i}
\nc{\B}{\mathbf{B}}
\nc{\F}{\mathbf{F}}
\nc{\bP}{\mathbf{P}}
\nc{\V}{\mathbf{V}}
\nc{\La}{\Lambda}
\nc{\la}{\lambda}
\nc{\nbinom}[2]{\genfrac{}{}{0pt}{1}{{#1}}{{#2}}}
\nc{\qbinom}[2]{\left[\genfrac{}{}{0pt}{1}{{#1}}{{#2}}\right]}
\nc{\tmppath}{\mathcal{P}}
\nc{\fit}{\tilde{f}_i}
\nc{\eit}{\tilde{e}_i}
\nc{\Y}{\mathbf{Y}}
\nc{\A}{\mathbf{A}}
\nc{\ra}{\rightarrow}
\nc{\vep}{\varepsilon}
\nc{\vphi}{\varphi}
\nc{\g}{\mathfrak{g}}
\nc{\h}{\mathfrak{h}}
\nc{\oP}{\overline{P}}
\nc{\tmppathp}{\mathbf{p}}

\nc{\tris}{
\bsegment
\move(0 0)\lvec(10 0)\lvec(10 10)\lvec(0 0)\ifill f:0.7
\esegment
}
\nc{\recs}{
\bsegment
\move(0 0)\lvec(10 0)\lvec(10 5)\lvec(0 5)\lvec(0 0)\ifill f:0.7
\esegment
}
\nc{\hcvec}[5]{%
\getpos(#1 #3)\spx\spy
\getpos(#2 #3)\epx\epy
\getpos(#4 #5)\xoff\yoff
\realadd \spx \xoff \twox
\realadd \epx {-\xoff} \thrx
\realadd \spy \yoff \posy
\move({\spx} {\spy})
\clvec ({\twox} {\posy})({\thrx} {\posy})({\epx} {\epy})
\rmove(0 0)
}

\nc{\ahead}[2]{%
\cossin (0 0)({#1} {#2})\cosa\sina
\bsegment
  \drawdim in \setunitscale 0.065
  \realmult {-0.5} \cosa \hcosa
  \realmult {-0.5} \sina \hsina
  \move({\hcosa} {\hsina}) \ravec({\cosa} {\sina})
\esegment
}
\nc{\boxi}{%
{%
\savebox{\tmppic}{\begin{texdraw}
\small
\drawdim em
\textref h:C v:C
\setunitscale 0.55
\htext(0 0){$i$}
\move(-1 -1)\lvec(-1 1)\lvec(1 1)\lvec(1 -1)\lvec(-1 -1)
\end{texdraw}}%
\raisebox{-0.19\height}{\usebox{\tmppic}}%
}%
}
\nc{\boxj}{%
{%
\savebox{\tmppic}{\begin{texdraw}
\small
\drawdim em
\textref h:C v:C
\setunitscale 0.55
\htext(0 0.1){$j$}
\move(-1 -1)\lvec(-1 1)\lvec(1 1)\lvec(1 -1)\lvec(-1 -1)
\end{texdraw}}%
\raisebox{-0.19\height}{\usebox{\tmppic}}%
}%
}
\nc{\boxipo}{%
{%
\savebox{\tmppic}{\begin{texdraw}
\small
\drawdim em
\textref h:C v:C
\setunitscale 0.55
\htext(0.15 0){$i\!\!+\!\!1$}
\move(-1.4 -1)\lvec(-1.4 1)\lvec(1.4 1)\lvec(1.4 -1)\lvec(-1.4 -1)
\end{texdraw}}%
\raisebox{-0.19\height}{\usebox{\tmppic}}%
}%
}
\everytexdraw{
\drawdim in \arrowheadsize l:0.065 w:0.03 \arrowheadtype t:F
\textref h:C v:C
}
\newsavebox{\tmppic}
\newsavebox{\tmpfig}
\newsavebox{\tmpdraw}
\newsavebox{\tmpfiga}
\newsavebox{\tmpfigb}
\newsavebox{\tmpfigc}
\newsavebox{\tmpfigd}
\newsavebox{\tmpfige}
\newsavebox{\tmpfigf}
\newsavebox{\tmpfigg}
\newsavebox{\tmpfigh}
\newsavebox{\tmpfigi}
\newsavebox{\tmpfigj}
\newsavebox{\tmpfigk}
\newsavebox{\tmpfigl}
\newsavebox{\tmpfigm}
\newsavebox{\tmpfign}
\newsavebox{\tmpfigo}
\newsavebox{\tmpfigp}
\newsavebox{\tmpfigq}
\newsavebox{\tmpfigr}
\newsavebox{\tmpfigs}
\newsavebox{\tmpfigt}
\newsavebox{\tmpfigu}
\newsavebox{\tmpfigv}
\newsavebox{\tmpfigw}
\newsavebox{\tmpfigx}
\newsavebox{\tmpfigy}
\newsavebox{\tmpfigz}

\newsavebox{\tmpfigaa}
\newsavebox{\tmpfigab}
\newsavebox{\tmpfigac}
\newsavebox{\tmpfigad}
\newsavebox{\tmpfigae}
\newsavebox{\tmpfigaf}
\newsavebox{\tmpfigag}
\newsavebox{\tmpfigah}
\newsavebox{\tmpfigai}
\newsavebox{\tmpfigaj}
\newsavebox{\tmpfigak}
\newsavebox{\tmpfigal}
\newsavebox{\tmpfigam}
\newsavebox{\tmpfigan}
\newsavebox{\tmpfigao}
\newsavebox{\tmpfigap}
\newsavebox{\tmpfigaq}
\newsavebox{\tmpfigar}
\newsavebox{\tmpfigas}
\newsavebox{\tmpfigat}
\newsavebox{\tmpfigau}
\newsavebox{\tmpfigav}
\newsavebox{\tmpfigaw}
\newsavebox{\tmpfigax}
\newsavebox{\tmpfigay}
\newsavebox{\tmpfigaz}

\newsavebox{\tmpfigba}
\newsavebox{\tmpfigbb}
\newsavebox{\tmpfigbc}
\newsavebox{\tmpfigbd}
\newsavebox{\tmpfigbe}
\newsavebox{\tmpfigbf}
\newsavebox{\tmpfigbg}
\newsavebox{\tmpfigbh}

\nc{\node}{\lcir r:1 }
\nc{\sline}{\bsegment\savepos(10 0)(*ex *ey)
            \move(1 0)\rlvec(8 0)
            \esegment\move(*ex *ey)}
\nc{\dline}{\bsegment\savepos(10 0)(*ex *ey)
            \move(0.93 0.4)\rlvec(8.14 0)\rmove(0 -0.8)\rlvec(-8.14 0)
            \esegment\move(*ex *ey)}
\nc{\uline}{\bsegment\savepos(0 10)(*ex *ey)
            \move(0 1)\rlvec(0 8)
            \esegment\move(*ex *ey)}
\nc{\lpoint}{\savecurrpos(*ex *ey)
             \rmove(2.5 1.5)\rlvec(-1.5 -1.5)\rlvec(1.5 -1.5)
             \move(*ex *ey)}
\nc{\rpoint}{\savecurrpos(*ex *ey)
             \rmove(-2.5 -1.5)\rlvec(1.5 1.5)\rlvec(-1.5 1.5)
             \move(*ex *ey)}
\nc{\bline}{\bsegment\savepos(10 0)(*ex *ey)
            \linewd 0.6 \move(1.1 0)\rlvec(7.8 0)
            \esegment\move(*ex *ey)}

\nc{\araise}[1]{\raisebox{4.5pt}{#1}}
\nc{\braise}[1]{\raisebox{12.1pt}{#1}}
\nc{\craise}[1]{\raisebox{8pt}{#1}}
\nc{\draise}[1]{\raisebox{12pt}{#1}}

\begin{document}

\title[Crystal bases and combinatorics of Young walls]
      {Crystal bases for quantum affine algebras \\
and combinatorics of Young walls}
\author[S.-J. Kang]{Seok-Jin Kang$^{*}$}
\address{Department of Mathematics\\
         Seoul National University\\
         Seoul 151-742, Korea}
\thanks{$^{*}$This research was supported by KOSEF Grant
\# 98-0701-01-5-L and the Young Scientist Award, 
Korean Academy of Science and Technology}
\email{sjkang@math.snu.ac.kr}

\begin{abstract}
In this paper, we give a realization of crystal bases 
for quantum affine algebras using some new combinatorial
objects which we call the Young walls.
The Young walls consist of colored blocks 
with various shapes that are built on the 
given ground-state wall and can be viewed as generalizations of 
Young diagrams.
The rules for building Young walls and the action of
Kashiwara operators are given explicitly in terms of
combinatorics of Young walls. 
The crystal graphs for basic representations are 
characterized as the set of all reduced proper Young walls.
The characters of basic representations can be computed easily
by counting the number of colored blocks that have been 
added to the ground-state wall. 
\end{abstract}

\maketitle

\vskip 1cm

\baselineskip=14pt

\section*{Introduction}

The {\it crystal basis theory} for integrable modules over 
quantum groups was introduced by Kashiwara 
using a combinatorial method (\cite{Kas90, Kas91}).
In \cite{Lus90}, a more geometric approach was developed by Lusztig
and is called the {\it canonical basis theory}.
The crystal bases can be viewed as bases at $q=0$ and they  
are given a structure of colored oriented graph, 
called the {\it crystal graphs},
with arrows defined by the {\it Kashiwara operators}.
The crystal graphs have many nice combinatorial features reflecting the
internal structure of integrable modules over quantum groups.
For instance, one of the major goals in representation theory
is to find an explicit expression for the characters of representations
and this goal can be achieved by finding an explicit combinatorial 
description of crystal bases.
Moreover, the crystal bases have extremely nice behavior with 
respect to taking the tensor product.
The action of the Kashiwara operators is given by the 
{\it tensor product rule} and the irreducible decomposition of the
tensor product of integrable modules is equivalent to decomposing
the tensor product of crystal graphs into a disjoint union of connected
components.
Thus the crystal basis theory provides us with a very
powerful combinatorial method of studying 
the structure of integrable modules over  quantum groups.

\vskip 3mm

In \cite{MM}, Misra and Miwa constructed the crystal bases for
basic representations of quantum affine algebras $U_q(A_n^{(1)})$ using
the {\it Fock space representation} on the space of 
Young diagrams with colored boxes.
The crystal graphs were characterized as the set of
{\it $n$-reduced Young diagrams} (with colored boxes).
Their idea was extended to construct crystal bases for irreducible
highest weight $U_q(A_n^{(1)})$-modules with arbitrary
higher level \cite{JMMO}.
The crystal graphs constructed in \cite{JMMO} and \cite{MM} 
can be parametrized by certain {\it paths}
which arise naturally in the theory of solvable lattice models.

\vskip 3mm

Motivated by this observation, Kang, Kashiwara, Misra, Miwa, Nakashima
and Nakayashiki developed the theory of {\it perfect crystals} for
general quantum affine algebras and gave a realization of crystal
graphs for irreducible highest weight modules over
classical quantum affine algebras with arbitrary higher levels
in terms of \emph{paths} (\cite{KMN1, KMN2}).
In this way, the theory of vertex models can be explained in the
language of representation theory of quantum affine algebras
and the 1-point function of the vertex model was
expressed as the quotient of the string function and the character
of the corresponding irreducible highest weight representation.

\vskip 3mm

On the other hand, in \cite{LLT}, 
Lascoux, Leclerc and Thibon discovered a surprising 
connection between the crystal basis theory for quantum affine 
algebras and modular representation theory of Hecke algebras. 
For a complex number $\zeta$, let 
$H_{N}(\zeta)$ denote the Hecke algebra of type $A_{N-1}$.
It is known that the indecomposable modules over 
$H_{N}(\zeta)$, called the {\it Specht modules}, 
are parametrized by Young diagrams with $N$ boxes.
We denote by $S(Y)$ the Specht module corresponding 
to the Young diagram $Y$ with $N$ boxes. 
If $\zeta$ is not a root of unity, then the category 
of finite dimensional $H_{N}(\zeta)$-modules is semisimple,
and hence the set of Specht modules coincides with 
the set of irreducible modules.
However, if $\zeta$ is a primitive $n$-th root of unity, then
the Specht modules are no longer irreducible, and
the irreducible $H_{N}(\zeta)$-modules are characterized 
as the unique iireducible quotients of the Specht modules
corresponding to the $n$-reduced Young diagrams with $N$ boxes 
(see, for example, \cite{DJ1, DJ2}). 

\vskip 3mm 

Recall that an $n$-reduced Young diagram $Y$ with $N$ boxes
can be viewed as a crystal basis element of a basic 
representation of the quantum affine algebra $U_q(A_n^{(1)})$. 
Using the $U_q(A_n^{(1)})$-module action on the space of 
all Young diagrams (see, for example, \cite {JMMO, MM}), one
should be able to write the {\it global basis element} 
(or {\it canonical basis element}) 
$G(Y)$ as a linear combination of Young diagrams:
$$G(Y)= \sum_{Y'} d_{Y',Y}(q) Y' \quad \text{for some} 
\ \ d_{Y', Y}(q) \in \C[q].$$ 
In \cite{LLT}, Lascoux, Leclerc and Thibon gave a recursive
algorithm of computing the polynomials $d_{Y',Y}(q)$ and 
conjectured that 

\vskip 3mm

\qquad\qquad\qquad $d_{Y',Y}(1) = [S(Y'): D(Y)],$

\vskip 3mm 
\noindent
where $D(Y)$ is the irreducible $H_{N}(\zeta)$-module 
corresponding to an $n$-reduced Young diagram $Y$, 
$S(Y')$ is the Specht module 
corresponding to a Young diagram $Y'$, 
and $[S(Y'): D(Y)]$ denotes the multiplicity of $D(Y)$
in a composition series of $S(Y')$. 
This conjecture was proved by Ariki (\cite{Ar}).
In \cite{GW}, Goodman and Wenzl found a faster 
algorithm of computing the polynomials $d_{Y',Y}(q)$
and gave a combinatorial proof 
of the fact that the polynomials $d_{Y',Y}(q)$ coincide with
the affine Kazhdan-Lusztig polynomials. 

\vskip 3mm

The purpose of this paper is to give a realization of crystal bases 
for quantum affine algebras using some new combinatorial
objects which we call the {\it Young walls}.
The Young walls consist of colored blocks 
with various shapes that are built on
the given {\it ground-state wall} 
and can be viewed as
generalizations of Young diagrams.
For the classical quantum affine algebras of type 
$A_n^{(1)}$ $(n\ge 1)$, $A_{2n-1}^{(2)}$ $(n\ge 3)$,
$D_n^{(1)}$ $(n\ge 4)$, $A_{2n}^{(2)}$ $(n\ge 2)$,
$D_{n+1}^{(2)}$ $(n\ge 2)$ and $B_n^{(1)}$ $(n\ge 3)$, 
the rules for building Young walls and the action of Kashiwara 
operators are given explicitly in terms of combinatorics of 
Young walls.
They are quite similar to playing with LEGO blocks and the Tetris game.
The crystal graphs for basic representations are characterized as
the set of all {\it reduced proper Young walls}.
The characters of basic representations can be computed easily
by counting the number of colored blocks in reduced proper Young walls 
that have been added to the ground-state wall. 

\vskip 3mm

We expect that we will be able to define the $U_q(\g)$-module
action on the space of all proper Young walls and find an 
effective algorithm for computing the global basis element 
corresponding to each reduced proper Young wall. 
For some mysterious reason, it is more difficult to deal with
the quantum affine algebras of type $C_n^{(1)}$
than the other classical quantum affine algebras. 
In \cite{HK}, we gave a characterization of crystal graphs 
for basic representations of $U_q(C_2^{(1)})$ as the set of
reduced proper Young walls. 
But it still remains to extend the results of \cite{HK} to
quantum affine algebras $U_q(C_n^{(1)})$ for $n\ge 3$. 
As the representation theory of Hecke algebras can be explained
in terms of Young diagrams, we also expect that there exists 
an interesting algebraic structure related to the combinatorics of
Young walls such that the irreducible modules at some specialization
are parametrized by reduced proper Young walls and that the 
decomposition matrices are determined by the polynomials giving the 
global basis elements.

\vskip 3mm

This paper is organized as follows.
In Section 1, we fix the notations and recall the basic 
definitions on quantum affine algebras. 
In Section 2, we review some of the fundamental properties 
of crystal bases for integrable modules over quantum groups. 
The tensor product rule for Kashiwara operators is explained
using the {\it $i$-signatures}. 
The definition of {\it perfect crystals} is given in Section 3.
We also include the examples of perfect crystals of level 1.
In Section 4, we recall one of the main results of 
\cite{KMN1} and \cite{KMN2}: the {\it path realization
of crystal graphs}. 
The crystal structure on the set of paths is given by the 
tensor product rule. 
We also demonstrate the top part of crystal graphs 
for basic representations of quantum affine algebras
of some small ranks. 

\vskip 3mm
Section 5 and Section 6 are devoted to the  combinatorics of 
Young walls. In Section 5, we define the notion of 
{\it ground-state walls}, {\it Young walls} and {\it proper Yong walls},
and give an explicit description of the rules and patterns 
for building Young walls. 
In Section 6, we define an affine crystal structure on the 
set of all proper Young walls. 
We also show that the set of 
all {\it reduced proper Young walls} is stable under the
action of Kashiwara operators. 

\vskip 3mm

Finally, in Section 7, we prove the main result of this
paper (Theorem \ref{thm:main}).
We prove that the crystal graph $B(\la)$ of the basic representation 
is isomorphic to the affine crystal $\Y(\la)$
consisting of reduced proper Young walls built on the 
ground-state wall $Y_{\la}$.
The path realization of crystal graphs plays a crucial 
role in the proof of our theorem. 
At the end of Section 7, we give some examples of the affine
crystal $\Y(\la)$ consisting of reduced proper Young walls.

\vspace{5mm}
\noindent\textbf{Acknowledgments.} \ 
The author would like to thank 
Jin Hong for many valuable discussions and his generous 
help in making the diagrams in this paper.
Part of this work was completed
while the author was visiting Massachusetts Institute of 
Technology, Yale University and 
Korea Institute for Advanced Study in the year of 1999.
He is very grateful to all the faculty and staff members of these
institutions for their hospitality and support during his visit.

\vskip 1cm

\section{Quantum affine algebras}

Let $I=\{0,1,\cdots, n \}$ be an index set and let 
$A=(a_{ij})_{i,j \in I}$ be a generalized Cartan matrix of affine type. 
Consider a free abelian group of rank $n+2$
\begin{equation}
\pv = \Z h_0 \op \Z h_1 \op \cdots \op \Z h_n \op \Z d
\end{equation}
and let ${\mathfrak h}= \C \otimes_{\Z} \pv$
be its complexification.
The free abelian group $\pv$ is called the {\it dual weight lattice}
and the complex vector space ${\mathfrak h}$ 
is called the {\it Cartan subalgebra}. 

\vskip 3mm

We define the linear functionals $\alpha_i$ and $\Lambda_i$ $(i\in I)$
on ${\mathfrak h}$ by
\begin{equation}
\begin{aligned}\mbox{}
\alpha_i (h_j) & = a_{ji}, \quad \alpha_i (d) = \delta_{0, i}, \\
\Lambda_i(h_j) & = \delta_{ij}, \quad \Lambda_i(d)=0 \qquad (i,j \in I).
\end{aligned}
\end{equation}
The $\alpha_i$ (resp. $h_i$) are called the {\it simple roots} 
(resp. {\it simple coroots}) and the $\Lambda_i$ are called the
{\it fundamental weights}. 
We denote by $\Pi=\{\alpha_i | \, i\in I\}$ 
(resp. $\Pi^{\vee} = \{ h_i | \, i\in I\}$)
the set of simple roots (resp. simple coroots).
The {\it affine weight lattice} is defined to be 
\begin{equation}
P=\{ \lambda \in {\mathfrak h}^* | \ \lambda (\pv) \subset \Z \}.
\end{equation}

\vskip 3mm 

The quintuple $(A, \Pi, \Pi^{\vee}, P, P^{\vee})$ is called an 
{\it affine Cartan datum}.
To each affine Cartan datum, we
can associate an infinite dimensional Lie algebra  $\mathfrak{g}$
called the {\it affine Kac-Moody algebra} (\cite{Kac90}). 
The center of the affine Kac-Moody algebra $\mathfrak{g}$
is 1-dimensional and is generated by the {\it canonical central element}
\begin{equation}
c=c_0 h_0 + c_1 h_1 + \cdots + c_n h_n.
\end{equation} 
Moreover, the imaginary roots of ${\mathfrak g}$
are nonzero integral multiples of the {\it null root} 
\begin{equation}
\delta = d_0 \alpha_0 + d_1 \alpha_1 + \cdots + d_n \alpha_n.
\end{equation}
Here, $c_i$ and $d_i$ $(i\in I)$ 
are the non-negative integers given in \cite{Kac90}. 
Using the fundamental weights and the null root, 
the affine weight lattice can be written as 
\begin{equation}
P=\Z \Lambda_0 \op \Z \Lambda_1 
\op \cdots \op \Z \Lambda_n \op \Z \delta.
\end{equation}
Set 
\begin{equation}
P^+=\{ \lambda \in P \, | \, \lambda(h_i) \in \Z_{\ge 0}
\ \ \ \text{for all} \ \ i\in I \}.
\end{equation} 
The elements of $P$ (resp. $P^{+}$) are called the {\it affine weights}
(resp. {\it affine dominant integral weights}). 
The {\it level} of an affine dominant integral weight $\lambda \in P^+$
is defined to be the nonnegative integer $\lambda(c)$. 

\vskip 5mm

\begin{example}
In this example, we present the affine Dynkin diagrams, the canonical central
elements, the null roots and the level 1 dominant integral weights 
for affine Cartan data of classical type.

\vskip 3mm

(a)  $A_n^{(1)}$ ($n\geq 1$) : 

\vskip 2mm
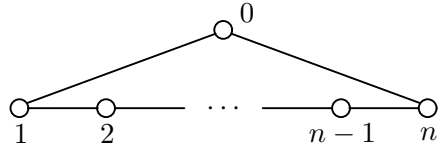
\begin{center}
\begin{texdraw}
\small
\drawdim em
\setunitscale 0.3
\move(0 0)
\lcir r:1
\move(1 0)\lvec(9 0)
\move(10 0)\lcir r:1
\move(11 0)\lvec(19 0)
\htext(23.6 0){$\cdots$}
\move(28 0)\lvec(36 0)
\move(37 0)\lcir r:1
\move(38 0)\lvec(46 0)
\move(47 0)\lcir r:1
\move(23.5 9)\lcir r:1
\move(0.9 0.6)\lvec(22.6 8.6)
\move(24.4 8.6)\lvec(46.1 0.6)
\htext(0 -3){$1$}
\htext(10 -3){$2$}
\htext(37 -3){$n-1$}
\htext(47 -3){$n$}
\htext(26 11){$0$}
\move(-2 12.5)\move(49 -4.6)
\end{texdraw}
\end{center}

\vskip 2mm
\begin{equation*}
\begin{aligned}\mbox{}
& c = h_0 + h_1 + \cdots + h_n, \\
& \delta = \alpha_0 + \alpha_1 + \cdots + \alpha_n, \\
& \lambda = \Lambda_i \ \ \ (i\in I).
\end{aligned}
\end{equation*}

\vskip 5mm
(b) $A_{2n-1}^{(2)}$ ($n\geq 3$) :

\vskip 2mm
\begin{center}
\begin{texdraw}
\small
\drawdim em
\setunitscale 0.3
\move(0 0)
\lcir r:1
\move(1 0)\lvec(9 0)
\move(10 0)\lcir r:1
\move(11 0)\lvec(19 0)
\htext(23.6 0){$\cdots$}
\move(28 0)\lvec(36 0)
\move(37 0)\lcir r:1
\move(38 0.4)\lvec(46 0.4)
\move(38 -0.4)\lvec(46 -0.4)
\move(39.2 1.2)\lvec(38 0)\lvec(39.2 -1.2)
\move(47 0)\lcir r:1
\move(10 1)\lvec(10 9)
\move(10 10)\lcir r:1
\htext(0 -3){$1$}
\htext(10 -3){$2$}
\htext(37 -3){$n-1$}
\htext(47 -3){$n$}
\htext(12.5 11){$0$}
\move(-2 -4.6)\move(49 12.6)
\end{texdraw}
\end{center}

\vskip 2mm
\begin{equation*}
\begin{aligned}\mbox{}
& c = h_0 + h_1 + 2 h_2 +\cdots + 2h_{n-1}+2 h_n, \\
& \delta = \alpha_0 + \alpha_1 
+ 2\alpha_2 + \cdots + 2\alpha_{n-1} + \alpha_n, \\
& \lambda=\Lambda_0,  \ \Lambda_1.
\end{aligned}
\end{equation*}

\vskip 5mm
(c) $D_n^{(1)}$ ($n\geq 4$)

\vskip 2mm
\begin{center}
\begin{texdraw}
\small
\drawdim em
\setunitscale 0.3
\move(0 0)
\lcir r:1
\move(1 0)\lvec(9 0)
\move(10 0)\lcir r:1
\move(11 0)\lvec(19 0)
\htext(23.6 0){$\cdots$}
\move(28 0)\lvec(36 0)
\move(37 0)\lcir r:1
\move(38 0)\lvec(46 0)
\move(47 0)\lcir r:1
\move(10 1)\lvec(10 9)\move(10 10)\lcir r:1
\move(37 1)\lvec(37 9)\move(37 10)\lcir r:1
\htext(0 -3){$1$}
\htext(10 -3){$2$}
\htext(37 -3){$n-2$}
\htext(47 -3){$n$}
\htext(12.5 11){$0$}
\htext(42.5 11){$n-1$}
\move(-2 -4.6)\move(49 12.6)
\end{texdraw}
\end{center}

\vskip 2mm
\begin{equation*}
\begin{aligned}\mbox{}
& c = h_0+h_1+2h_2+\cdots+2h_{n-2}+h_{n-1}+h_n,\\
& \delta = \alpha_0+\alpha_1+2\alpha_2+\cdots+2\alpha_{n-2}+\alpha_{n-1}
+\alpha_n, \\
& \lambda = \Lambda_0, \ \Lambda_1, \ \Lambda_{n-1}, \
\Lambda_n.
\end{aligned}
\end{equation*}

\vskip 5mm
(d)  $A_{2n}^{(2)}$ ($n\geq 2$)

\vskip 2mm
\begin{center}
\begin{texdraw}
\small
\drawdim em
\setunitscale 0.3
\move(0 0)
\lcir r:1
\move(1 0)\lvec(9 0)
\move(10 0)\lcir r:1
\move(11 0)\lvec(19 0)
\htext(23.6 0){$\cdots$}
\move(28 0)\lvec(36 0)
\move(37 0)\lcir r:1
\move(38 0.4)\lvec(46 0.4)
\move(38 -0.4)\lvec(46 -0.4)
\move(39.2 1.2)\lvec(38 0)\lvec(39.2 -1.2)
\move(47 0)\lcir r:1
\move(-1 0.4)\lvec(-9 0.4)
\move(-1 -0.4)\lvec(-9 -0.4)
\move(-10 0)\lcir r:1
\move(-7.8 1.2)\lvec(-9 0)\lvec(-7.8 -1.2)
\htext(-10 -3){$0$}
\htext(0 -3){$1$}
\htext(10 -3){$2$}
\htext(37 -3){$n-1$}
\htext(47 -3){$n$}
\move(-12 -4.6)\move(49 3)
\end{texdraw}
\end{center}

\vskip 2mm
\begin{equation*}
\begin{aligned}\mbox{}
& c = h_0+2h_1+\cdots+2h_{n-1}+2h_n,\\
& \delta = 2\alpha_0 + 2\alpha_1 + \cdots + 2\alpha_{n-1}+\alpha_{n}, \\
& \lambda=\Lambda_0.
\end{aligned}
\end{equation*}

\vskip 5mm
%\newpage
(e) $D_{n+1}^{(2)}$ ($n\geq2$)

\vskip 2mm
\begin{center}
\begin{texdraw}
\small
\drawdim em
\setunitscale 0.3
\move(0 0)
\lcir r:1
\move(1 0)\lvec(9 0)
\move(10 0)\lcir r:1
\move(11 0)\lvec(19 0)
\htext(23.6 0){$\cdots$}
\move(28 0)\lvec(36 0)
\move(37 0)\lcir r:1
\move(38 0.4)\lvec(46 0.4)
\move(38 -0.4)\lvec(46 -0.4)
\move(44.8 1.2)\lvec(46 0)\lvec(44.8 -1.2)
\move(47 0)\lcir r:1
\move(-1 0.4)\lvec(-9 0.4)
\move(-1 -0.4)\lvec(-9 -0.4)
\move(-10 0)\lcir r:1
\move(-7.8 1.2)\lvec(-9 0)\lvec(-7.8 -1.2)
\htext(-10 -3){$0$}
\htext(0 -3){$1$}
\htext(10 -3){$2$}
\htext(37 -3){$n-1$}
\htext(47 -3){$n$}
\move(-12 -4.6)\move(49 3)
\end{texdraw}
\end{center}

\vskip 2mm
\begin{equation*}
\begin{aligned}\mbox{}
& c = h_0+2h_1+\cdots+2h_{n-1}+h_n,\\
& \delta=\alpha_0+\alpha_1+ \cdots+\alpha_{n-1}+\alpha_n, \\
& \lambda=\Lambda_0, \ \Lambda_n. 
\end{aligned}
\end{equation*}

\vskip 5mm
(f) $B_n^{(1)}$ ($n\geq3$)

\vskip 2mm
\begin{center}
\begin{texdraw}
\small
\drawdim em
\setunitscale 0.3
\move(0 0)
\lcir r:1
\move(1 0)\lvec(9 0)
\move(10 0)\lcir r:1
\move(11 0)\lvec(19 0)
\htext(23.6 0){$\cdots$}
\move(28 0)\lvec(36 0)
\move(37 0)\lcir r:1
\move(38 0.4)\lvec(46 0.4)
\move(38 -0.4)\lvec(46 -0.4)
\move(44.8 1.2)\lvec(46 0)\lvec(44.8 -1.2)
\move(47 0)\lcir r:1
\move(10 1)\lvec(10 9)
\move(10 10)\lcir r:1
\htext(0 -3){$1$}
\htext(10 -3){$2$}
\htext(37 -3){$n-1$}
\htext(47 -3){$n$}
\htext(12.5 11){$0$}
\move(-2 -4.6)\move(49 12.6)
\end{texdraw}
\end{center}

\vskip 2mm
\begin{equation*}
\begin{aligned}\mbox{}
& c = h_0 + h_1 + 2h_2 + \cdots + 2h_{n-1} + h_n,\\
& \delta = \alpha_0+\alpha_1+2\alpha_2+\cdots+2\alpha_{n-1}+2\alpha_n, \\
& \lambda= \Lambda_0, \ \Lambda_1, \  \Lambda_n.
\end{aligned}
\end{equation*}

\end{example}

\vskip 5mm

We denote by $q^h$ $(h\in P^{\vee})$ the basis elements
of the group algebra $\C(q)[P^{\vee}]$ with the multiplication
$q^h q^{h'} = q^{h+h'}$ $(h,h'\in \pv)$.
Let $(\ \ | \ \ )$ be a nondegenerate symmetric bilinear form on
${\mathfrak h}^*$ satisfying 
\begin{equation*}
\frac{2(\alpha_i | \alpha_j)}{(\alpha_i | \alpha_i)} = a_{ij}
\qquad \text{for all} \quad i,j \in I.
\end{equation*}
Set \
%\begin{equation*}
$q_i = q^{\frac{(\alpha_i |\alpha_i)}{2}},
\ \  
K_i = q^{\frac{(\alpha_i |\alpha_i)}{2} h_i}$ \ 
%\end{equation*}
and define 
\begin{equation*}
[k]_i = \frac{q_i^k - q_i^{-k}}{q_i - q_i^{-1}},
\quad
[n]_i ! = \prod_{k=1}^{n} [k]_i,
\quad
%\qbinom{m}{n}_{i}
{\begin{bmatrix} m \\ n \end{bmatrix}}_{i}
=\frac{[m]_{i}!}{[n]_{i}! \, [m-n]_{i}!}. 
\end{equation*}
We will also use the following notations:
\begin{equation*}
e_i^{(n)} = e_i^n/[n]_i!, \qquad
f_i^{(n)} = f_i^n/[n]_i!.
\end{equation*}

\vskip 5mm 
\begin{defi}
The {\it quantum affine algebra} $U_q({\mathfrak g})$ 
associated with the affine Cartan datum 
$(A, \Pi, \Pi^{\vee}, P, P^{\vee})$ is the associative algebra 
with 1 over $\C(q)$ 
generated by the symbols $e_i$, $f_i$ $(i\in I)$ and $q^h$
$(h\in \pv)$ subject to the following defining relations:
\end{defi}
\begin{equation}
\begin{aligned}\mbox{}
\ & q^0 = 1, \ \ q^h q^{h'} = q^{h+ h'} \quad (h, h'\in P^{\vee}),\\
\ & q^h e_i q^{-h} = q^{\ali(h)} e_i, \quad 
q^h f_i q^{-h} = q^{-\ali(h)} f_i \quad (h\in P^{\vee}, i\in I), \\
\ & e_i f_j - f_j e_i = \delta_{ij} \frac{K_i - K_i^{-1}}{q_i - q_i^{-1}}
\quad (i,j \in I), \\
\ & \sum_{k=0}^{1-a_{ij}} (-1)^k 
{\begin{bmatrix} 1-a_{ij} \\ k \end{bmatrix}}_{i}
e_i^{1-a_{ij}-k} e_j e_i^{k} = 0 \qquad (i \neq j), \\
\ & \sum_{k=0}^{1-a_{ij}} (-1)^k 
{\begin{bmatrix} 1-a_{ij} \\ k \end{bmatrix}}_{i}
f_i^{1-a_{ij}-k} f_j f_i^{k} = 0 \qquad (i \neq j).
\end{aligned}
\end{equation}

The subalgebra of $U_q({\mathfrak g})$ generated by 
$e_i$, $f_i$, $K_i^{\pm1}$ $(i\in I)$ is denoted by 
$U'_q({\mathfrak g})$, and is 
also called the {\it quantum affine algebra}. 
Let 
\begin{equation*}
\oP^{\vee} = \Z h_0 \op \Z h_1 \op \cdots \op \Z h_n
\quad \text{and} \quad 
\overline {\mathfrak h} = \C \otimes_{\Z} \oP^{\vee}.
\end{equation*}
Consider $\alpha_i$ and $\Lambda_i$ $(i\in I)$ as linear functionals
on $\overline{\mathfrak h}$ and set 
\begin{equation*}
\oP=\Z \Lambda_0 \op \Z \Lambda_1 \op \cdots \op \Z \Lambda_n.
\end{equation*}
For example, $\delta = 0$ as an element in $\oP$. 
The elements of $\oP$ are called the {\it classical weights}. 
The algebra $U'_q({\mathfrak g})$ can be regarded as the quantum 
affine algebra associated with the {\it classical Cartan datum}
$(A, \Pi, \Pi^{\vee}, \oP, \oP^{\vee})$. 

\vskip 3mm

The projection $\text{cl} : P \longrightarrow \oP$
will be denoted by $\lambda \longmapsto \overline{\lambda}$
and we will fix an embedding $\text{aff}: \oP \longrightarrow P$
such that 
$$\text{cl} \circ \text{aff} = \text{id}, 
\qquad \text{aff} \circ \text{cl} (\alpha_i)=\alpha_i
\quad \text{for} \quad i\neq 0.$$
We define 
$$\oP^{+} = \text{cl} (P^{+}) 
=\{\la \in \oP | \ \la(h_i) \ge 0 \ \ \text{for all} \ \, i\in I\}.$$
The elements of $\oP^{+}$ are called the 
{\it classical dominant integral weights}. 
For simplicity, we will omit the notations for the projection
and the embedding. 
We will state explicitly whether a linear functional is an affine 
weight or a classical weight whenever it could cause a confusion. 
A classical dominant integral weight $\lambda \in \oP^{+}$ is said to 
have {\it level $l \in \Z_{\ge 0}$} if $\lambda(c)= l$.
Note that it has the same level as its affine counterpart.

\vskip 1cm 
\section{Crystal bases}

In this section, we briefly review the crystal basis theory for quantum affine 
algebras.
The {\it category ${\mathcal{O}}_{int}$} consists of $U_q(\g)$-modules 
(resp. $U'_q(\g)$-modules $M$) satisfying the properties:

(i)  $M=\bigoplus_{\lambda \in P} M_{\lambda}$ 
(resp. $M=\bigoplus_{\lambda \in \oP} M_{\lambda}$),
where 
\begin{equation*}
M_{\la} = \{ v\in M \mid q^h v = q^{\la(h)} v \text{ for all }h\in \pv
\ (\text {resp. } h\in \oP^{\vee}) \},
\end{equation*}

(ii)  for each $i\in I$, $M$ is a direct sum of 
finite dimensional irreducible $U_i$-modules, where $U_i$ denotes
the subalgebra generated by $e_i$, $f_i$, $K_i^{\pm1}$ which is
isomorphic to $\Uq(\mathfrak{sl}_2)$. 

\vskip 3mm
\noindent
It is known that every $U_q({\mathfrak g})$-module 
(resp. $U'_{q}(\mathfrak{g})$-module) in the category $\mathcal{O}_{int}$
is a direct sum of irreducible highest weight modules $V(\lambda)$
with $\la \in P^{+}$ (resp. $\la \in {\bar P}^{+}$). 

\vskip 3mm 

Fix an index $i\in I$.
By the representation theory of $\Uq(\mathfrak{sl}_2)$, every 
element $v\in M_\la$ can  be written uniquely as
\begin{equation*}
v = \sum_{k\geq0} f_i^{(k)} v_k,
\end{equation*}
where $k \ge -\la(h_i)$ and $v_k \in \ker e_i \cap M_{\la+k\ali}$.
We define the endomorphisms $\eit$ and $\fit$ on $M$,
called the {\it Kashiwara operators},  by
\begin{equation}
\eit v = \sum_{k\geq1} f_i^{(k-1)} v_k, \qquad 
\fit v = \sum_{k\geq0} f_i^{(k+1)} v_k.
\end{equation}
It is easy to see that 
\begin{equation*}
\eit M_{\la} \subset M_{\la + \ali}, 
\quad \fit M_{\la} \subset M_{\la - \ali}
\ \ \text {for all} \ \  i\in I. 
\end{equation*}

\vskip 3mm

Let $M$ be a $U_q(\g)$-module (or $U'_q(\g)$-module) in the category 
${\mathcal O}_{int}$ and let 
\begin{equation*}
\A = \{ f/g \in \C(q)  \, | \, f, g \in \C[q], \, g(0)\neq 0 \}
\end{equation*}
be the subring of $\C(q)$ consisting of the rational functions
in $q$ that are regular at $q=0$.
We now define the notion of {\it crystal bases} 
which was introduced in \cite{Kas90, Kas91}.

\vskip 3mm
\begin{defi}
A free $\A$-submodule $L$ of $M$ is called a {\it crystal lattice} if

(i) $L$ generates $M$ as a $\C(q)$-vector space; i.e., 
$M \cong \C(q) \otimes_{\A} L$,

(ii) $L=\bigoplus_{\lambda \in P} L_{\la}$ (or 
$L=\bigoplus_{\lambda \in \oP} L_{\la}$), where 
$L_{\la} = L \cap M_{\la}$,

(iii) $\eit L \subset L$, $\fit L \subset L$ for all $i\in I$.
\end{defi}

\vskip 3mm
\begin{defi}
A {\it crystal basis} of $M$ is a pair $(L, B)$ such that 

(i) $L$ is a crystal lattice of $M$,

(ii) $B$ is a $\C$-basis of $L/qL$,

(iii) $B=\bigsqcup_{\la \in P} B_{\la}$, where 
$B_{\la}=B \cap \left(L_{\la} / q L_{\la} \right)$,

(iv) $\eit B\subset B\cup \{0\}$, \ $\fit B\subset B\cup \{0\}$,

(v) for $b, b'\in B$, $\fit b = b'$ if and only if $b = \eit b'$.

\end{defi}

\noindent
The set $B$ is given a colored oriented graph structure with the 
arrow defined by 
$$b \stackrel{i} \longrightarrow b' \quad \text{if and only if}
\quad \fit b = b'.$$
The graph $B$ is called the {\it crystal graph} of $M$ and it reflects 
the combinatorial structure of $M$.
For instance, we have 
\begin{equation*}
\dim_{\C(q)} M_{\lambda} =\# B_{\lambda}
\qquad \text{for all} \quad 
\lambda \in P \quad (\text{or} \quad \lambda \in \oP). 
\end{equation*}

\vskip 3mm

We recall the basic properties of crystal graphs. 
Let $B$ be a crystal graph for a $U_q(\g)$-module
(or a $U'_q(\g)$-module) $M$ in the category ${\mathcal{O}}_{int}$.
For each $b\in B$ and $i\in I$, we define
\begin{equation}
\varepsilon_i(b)=\max \{k\ge 0 | \, \eit^k b \in B \}, 
\qquad 
\varphi_i(b) = \max \{ k\ge 0 | \, \fit^k b \in B \}.
\end{equation}

\noindent
Then the crystal graph $B$ satisfies the following properties. 

\vskip 3mm

\begin{prop} {\rm (\cite{Kas91, Kas93, Kas94})}

{\rm (a)} For all $i\in I$ and $b\in B$, 
we have 
\begin{equation*}
\begin{aligned}\mbox{}
& \varphi_i(b) = \varepsilon_i(b) + \langle h_i, \wt(b) \rangle, \\
& \wt(\eit b)=\wt(b) + \alpha_i, \\
& \wt(\fit b)=\wt(b) - \alpha_i. 
\end{aligned}
\end{equation*}

{\rm (b)} 
If $\eit b \in B$, then
$$\varepsilon_i(\eit b) = \varepsilon_i(b) - 1, \quad
\varphi_i(\eit b) = \varphi_i(b) + 1.$$

{\rm (c)} If $\fit b \in B$, then
$$\varepsilon_i(\fit b) = \varepsilon_i(b) + 1, \quad
\varphi_i(\fit b) = \varphi_i(b) - 1.$$

\end{prop}

\vskip 3mm
Moreover, the crystal bases have extremely simple behavior 
with respect to taking the tensor product.

\vskip 3mm

\begin{prop} {\rm (\cite{Kas90, Kas91})} \hfill

\vskip 2mm

Let $M_j$ $(j=1, 2)$ be a  $U_q(\g)$-module 
$(\text{or a $U'_q(\g)$-module})$ 
in the category ${\mathcal O}_{int}$ and $(L_j, B_j)$ be its crystal
basis. 
Set 
\begin{equation*}
L = L_1 \otimes_{\A} L_2, \quad B= B_1 \times B_2. 
\end{equation*}
Then $(L,B)$ is a crystal basis of $M_1 \otimes_{\C(q)} M_2$
with the Kashiwara operators on $B$ given by 
\begin{equation*}
\begin{aligned}\mbox{}
\eit(b_1\ot b_2)
  &=
  \begin{cases}
  \eit b_1 \ot b_2 & \text{if $\vphi_i(b_1) \geq \vep_i(b_2)$,}\\
  b_1\ot\eit b_2 & \text{if $\vphi_i(b_1) < \vep_i(b_2)$,}
  \end{cases}\\
\fit(b_1\ot b_2)  
  &= 
  \begin{cases}
  \fit b_1 \ot b_2 & \text{if $\vphi_i(b_1) > \vep_i(b_2)$,}\\
  b_1\ot\fit b_2 & \text{if $\vphi_i(b_1) \leq \vep_i(b_2)$.}
  \end{cases} 
\end{aligned}
\end{equation*}
\end{prop}

As an immediate corollary, we obtain the following simple
criterion for determining the maximal vectors in the tensor product
of crystal graphs. 

\begin{cor} {\rm (\cite{Kas90, Kas91})} \label{cor:maximal vectors}

{\rm (a)} Let $M_j$ be a $U_q(\g)$-module $(\text{or a $U'_q(\g)$-module})$
in the category ${\mathcal O}_{int}$
and let $(L_j, B_j)$ be a crystal basis of $M_j$ $(j=1,2)$.
Then $b_1 \ot b_2 \in B_1 \ot B_2$ is a maximal vector
$(\text{i.e., $\eit (b_1 \ot b_2)=0$ for all $i\in I$})$ if and only if
$\eit b_1 =0$ and $\langle h_i, \wt b_1 \rangle
\ge \varepsilon_i(b_2)$ for all $i\in I$.

{\rm (b)} Let $M_j$ be a $U_q(\g)$-module $(\text{or a $U'_q(\g)$-module})$
in the category ${\mathcal O}_{int}$ 
and let $(L_j, B_j)$ be a crystal basis of $M_j$ $(j=1, \cdots, N)$.
Then the vector $b_1 \ot \cdots \ot b_N \in B_1 \ot \cdots \ot B_N$
is a maximal vector if and only if $b_1 \ot \cdots \ot b_k$ is
a maximal vector for all $k=1,\cdots, N$.
\end{cor}

\vskip 3mm

The tensor product rule gives a very convenient
combinatorial description of the action of Kashiwara
operators on the multi-fold tensor product of crystal graphs.
Let $M_j$ be a $U_q(\g)$-module $(\text{or a $U'_q(\g)$-module})$
in the category ${\mathcal O}_{int}$ 
with a crystal basis $(L_j, B_j)$ \ $(j=1, \cdots, N)$.
Fix an index $i\in I$ and consider a vector
$b=b_1 \ot \cdots \ot b_N \in B_1 \ot \cdots \ot B_N$.
To each $b_j \in B_j$ \, $(j=1, \cdots , N)$,
we assign a sequence of $-$'s and $+$'s with as many $-$'s
as $\vep_i(b_j)$ followed by as many $+$'s as $\vphi_i(b_j)$:
\begin{equation*}
\begin{aligned}\mbox{}
b \, &  =  \, b_1 \otimes b_2 \otimes \cdots \otimes b_N \\
& \longmapsto 
 (\underbrace{-, \cdots, -}_{\vep_i(b_1)},
\underbrace{+, \cdots, +}_{\vphi_i(b_1)}, \cdots \cdots,
\underbrace{-, \cdots, -}_{\vep_i(b_N)},
\underbrace{+, \cdots, +}_{\vphi_i(b_N)}).
\end{aligned}
\end{equation*}
In this sequence, we cancel out all the $(+,-)$-pairs to obtain a
sequence of $-$'s followed by $+$'s:
\begin{equation}
\text{$i$-sgn}(b) = (-, -, \cdots, -, +, +, \cdots, +).
\end{equation}
The sequence $\text{$i$-sgn}(b)$ is called the {\it $i$-signature} of $b$.

Now the tensor product rule tells that
$\eit$ acts on $b_j$ corresponding to the right-most $-$ in
$\text{$i$-sgn}(b)$ and $\fit$ acts on $b_k$ corresponding to
the left-most $+$ in $\text{$i$-sgn}(b)$:
\begin{equation}\label{eq:signature}
\begin{aligned}\mbox{}
& \eit b = b_1 \otimes \cdots \ot \eit b_j \ot \cdots \ot b_N, \\
& \fit b = b_1 \otimes \cdots \ot \fit b_k \ot \cdots \ot b_N.
\end{aligned}
\end{equation}
We define $\eit b=0$ (resp. $\fit b=0$) if there is no $-$ 
(resp. $+$) in the $i$-signature of $b$.

\vskip 3mm 
We close this section with the existence and uniquenes theorem 
for crystal bases. 

\begin{prop} {\rm (\cite{Kas91})} \label{prop:existence}

{\rm (a)} Every $U_q(\g)$-module {\rm (}resp. $U'_q(\g)$-module{\rm )} 
in the category  ${\mathcal O}_{int}$ has a unique crystal basis. 
More precisely, let $V(\la)$ be the irreducible highest weight 
$U_q(\g)$-module {\rm (}resp. $U'_q(\g)$-module{\rm )} with highest weight
$\la \in P^{+}$ {\rm (}resp. $\la \in {\bar P}^{+}${\rm )} and highest weight 
vector $u_{\la}$. 
Let $L(\la)$ be the free $\A$-submodule of $V(\la)$ spanned by 
the vectors of the form $\tilde f_{i_1} \cdots \tilde f_{i_r} u_{\la}$
$(i_k \in I, r\in \Z_{\ge 0})$ and set 
\begin{equation*}
B(\la) = \{ \tilde f_{i_1} \cdots \tilde f_{i_r} u_{\la} \in 
L(\la) / q L(\la) \} \setminus \{0\}.
\end{equation*}
Then $(L(\la), B(\la))$ is a crystal basis of $V(\la)$ and every
crystal basis of $V(\la)$ is isomorphic to $(L(\la), B(\la))$.

{\rm (b)} Let $M$ be a $U_q(\g)$-module {\rm (}resp. a $U'_q(\g)$-module{\rm )} 
in the category  ${\mathcal O}_{int}$ and $(L,B)$ be its crystal basis.
Define an automorphism of $U_q(\g)$ {\rm (}resp. of $U'_q(\g)${\rm )} by
$\overline {e_i} = e_i$, $\overline {f_i} = f_i$, and 
$\overline {q^h} = q^{-h}$ for $i\in I, h\in P^{\vee}$ {\rm (}resp.
$h\in {\bar P}^{\vee}${\rm )}. 
Then there exists a unique $\C(q)$-basis $G=\{G(b) | \, b\in B \}$
of $M$ such that 
\begin{equation*}
\overline {G(b)} =G(b), \quad 
G(b) \equiv b \ \ (\text{mod} \, qL) \quad \text{for all} \ b\in B.
\end{equation*}
\end{prop}

\vskip 3mm

\begin{defi}
The basis $G$ of $M$ given in Proposition \ref{prop:existence} is 
called the {\it global basis} or the {\it canonical 
basis} of $M$ associated with the crystal graph $B$. 
\end{defi}

\vskip 1cm

\section{Perfect crystals}

By extracting properties of the crystal graphs, we define the 
notion of abstract {\it crystals} as follows (cf. \cite{Kas93, Kas94}). 

\begin{defi} 
An {\it affine crystal} (resp. {\it classical crystal}) 
is a set $B$ together with the maps
$\wt : B \ra P$ (resp. $\wt : B \ra \oP$),
$\vep_i : B \ra \Z\cup\{-\infty\}$,
$\vphi_i : B \ra \Z\cup\{-\infty\}$,
$\eit : B \ra B\cup\{0\}$, and
$\fit : B \ra B\cup\{0\}$
satisfying the following conditions:

\vskip 3mm

\hskip 3mm (i) for all $i\in I$, $b\in B$, we have 
\begin{equation*}
\begin{aligned}\mbox{}
& \varphi_i(b) = \varepsilon_i(b) + \langle h_i, \wt(b) \rangle, \\
& \wt(\eit b)=\wt(b) + \alpha_i, \\
& \wt(\fit b)=\wt(b) - \alpha_i. 
\end{aligned}
\end{equation*}

\hskip 3mm (ii) if $\eit b \in B$, then
$$\varepsilon_i(\eit b) = \varepsilon_i(b) - 1, \quad
\varphi_i(\eit b) = \varphi_i(b) + 1.$$

\hskip 3mm (iii) if $\fit b \in B$, then
$$\varepsilon_i(\fit b) = \varepsilon_i(b) + 1, \quad
\varphi_i(\fit b) = \varphi_i(b) - 1.$$

\hskip 3mm (iv) $\fit b = b'$ if and only if $b = \eit b'$ 
for all $i\in I$, $b, b' \in B$,

\vskip 3mm
\hskip 3mm (v) if $\vep_i(b) = - \infty$, then $\eit b = \fit b = 0$. 

\end{defi}

\vskip 3mm 
\noindent
The crystal graphs for $U_q(\g)$-modules
(resp. $U'_q(\g)$-modules) in the category ${\mathcal O}_{int}$
are affine crystals (resp. classical crystals). 

\vskip 3mm

\begin{defi}
Let $B_1$ and $B_2$ be (affine or classical) crystals. 
A {\it morphism} $\psi:B_1 \ra B_2$ of crystals is a map
$\psi:B_1 \cup \{0\} \ra B_2 \cup \{0\}$ satisfying the conditions:

\vskip 3mm
\hskip 3mm (i) $\psi(0)=0$,

\vskip 3mm 
\hskip 3mm (ii) if $b\in B_1$ and $\psi(b) \in B_2$, then 
$$\wt(\psi(b))= \wt(b), \ \  \vep_i(\psi(b))=\vep_i(b),
\ \ \vphi_i(\psi(b)) = \vphi_i(b),$$

\hskip 3mm (iii) if $b, b'\in B_1$, $\psi(b), \psi(b') \in B_2$ and $\fit b =b'$,
then $\fit \psi(b) = \psi(b')$. 
\end{defi}

\vskip 3mm
\noindent
A morphism of crystals is said to be {\it strict} if it commutes with 
the Kashiwara operators $\eit$ and $\fit$ $(i\in I)$.

\vskip 3mm

\begin{defi}
The {\it tensor product} $B_1 \otimes B_2$ of the crystals 
$B_1$ and $B_2$  is defined to be 
the set $B_1 \times B_2$ whose crystal structure is given by 
\begin{equation}
\begin{aligned}\mbox{}
\wt(b_1\ot b_2) &= \wt(b_1) + \wt(b_2),\\
\vep_i(b_1\ot b_2)
  &= \max(\vep_i(b_1), \vep_i(b_2) - \langle h_i, \wt(b_1) \rangle),\\
\vphi_i(b_1\ot b_2)
  &= \max(\vphi_i(b_2), \vep_i(b_1) + \langle h_i, \wt(b_2) \rangle),\\
\eit(b_1\ot b_2)
  &=
  \begin{cases}
  \eit b_1 \ot b_2 & \text{if $\vphi_i(b_1) \geq \vep_i(b_2)$,}\\
  b_1\ot\eit b_2 & \text{if $\vphi_i(b_1) < \vep_i(b_2)$,}
  \end{cases}\\
\fit(b_1\ot b_2)  
  &= 
  \begin{cases}
  \fit b_1 \ot b_2 & \text{if $\vphi_i(b_1) > \vep_i(b_2)$,}\\
  b_1\ot\fit b_2 & \text{if $\vphi_i(b_1) \leq \vep_i(b_2)$.}
  \end{cases} 
\end{aligned}
\end{equation}
\end{defi}
\noindent
Here, we denote $b_1 \otimes b_2 = (b_1, b_2)$ and use the
convention that $b_1 \ot 0 = 0 \ot b_2 =0$. 

\vskip 3mm

We now define the notion of {\it perfect crystals}. 
Let $B$ be a classical crystal.
For $b\in B$, we define
\begin{equation}
\vep(b) = \sum_i \vep_i(b)\La_i \quad \text {and} \quad 
\vphi(b) = \sum_i \vphi_i(b)\La_i.
\end{equation}
Note that 
\begin{equation*}
\wt(b) = \vphi(b) - \vep(b).
\end{equation*}
For a positive integer $l>0$, set 
\begin{equation}
\oP^+_l = \{\la\in \oP^+\mid \langle c, \la \rangle = l \}.
\end{equation}

\vskip 3mm 
\begin{defi}
For $l\in \Z_{> 0}$, we say that a finite classical crystal $\B$
is a {\it perfect crystal of level $l$} if

\vskip 2mm

(i) there is a finite dimensional $U'_q(\g)$-module with a crystal basis 
whose crystal graph is isomorphic to $\B$, 

\vskip 2mm
     
(ii) $\B \ot \B$ is connected,

\vskip 2mm

(iii)  there exists some $\la_0\in\oP$ such that 
      $$\wt(\B)\subset \la_0 + \sum_{i\neq0} \Z_{\leq0}\alpha_i, 
\qquad \#(\B_{\la_0}) = 1,$$

(iv) for any $b\in \B$, we have $\langle c, \vep(b) \rangle \geq l$,

\vskip 2mm
(v)  for each $\la \in \oP^+_l$, there exist unique
$b^{\la} \in \B$ and $b_{\la} \in \B$ such that 
$$\vep(b^{\la})=\la, \qquad \vphi(b_{\la})=\la.$$
\end{defi}

\vskip 3mm
\noindent
A finite dimensional $U'_q(\g)$-module $\V$ is called a
\emph{perfect representation} of \emph{level $l$} if it has a crystal
basis $(L,B)$ such that $B$ is isomorphic to 
a perfect crystal of level $l$. 

\vskip 3mm 
\noindent
{\it Remark.} For a perfect crystal $\B$, define 
\begin{equation*}
\B^{\text{min}}=\{ b\in \B | \, \langle c, \vep(b) \rangle = l \}.
\end{equation*}
Then the maps $\vep, \vphi : \B^{\text{min}} = \{ b\in B \mid
      \langle c, \vep(b) \rangle = l \} \longrightarrow \oP^+_l$ are bijective.

\vskip 3mm

In the following, we list some examples of perfect crystals 
of level 1 for the classical quantum affine algebras 
of type  $A_n^{(1)}$ $(n\ge 1)$, $A_{2n-1}^{(2)}$ $(n\ge 3)$, 
$D_n^{(1)}$ $(n\ge 4)$, $A_{2n}^{(2)}$ $(n\ge 2)$, 
$D_{n+1}^{(2)}$ $(n\ge 2)$ and $B_n^{(1)}$ $(n\ge 3)$.

\vskip 3mm

\begin{example}  \ $A_n^{(1)}$ $(n\ge 1)$ : 

\vskip 3mm

\begin{center}
\savebox{\tmppic}{\begin{texdraw}
\small
\drawdim em
\setunitscale 0.5
\nc{\numbox}[1]{\bsegment
\move(-1 -1)\lvec(1 -1)\lvec(1 1)\lvec(-1 1)\lvec(-1 -1)
\htext(0 0){#1}
\esegment}
\move(0 0)\numbox{$1$}
\move(1.5 0)\ravec(5 0)
\htext(3.7 1){$1$}
\move(8 0)\numbox{$2$}
\move(9.5 0)\ravec(5 0)
\htext(11.7 1){$2$}
\htext(16.7 0){$\cdots$}
\move(18.5 0)\ravec(5 0)
\htext(20.5 1){$n\!-\!1$}
\move(25 0)\numbox{$n$}
\move(26.5 0)\ravec(5 0)
\htext(28.7 1){$n$}
\move(33 0)\numbox{$0$}
\hcvec{0.5}{32.5}{-1.5}{4.5}{-3.5}
\move(0.5 -1.5)\ahead{-450}{260}
\htext(16.5 -5.2){$0$}
\move(-1 -6.2)\move(34 1.8)
\end{texdraw}}%
%$\B$ : \raisebox{-2.8em}{\usebox{\tmppic}}
\usebox{\tmppic}
\end{center}

\vskip 3mm
\noindent
Here, we have 
$$b^{\La_i} = \boxipo\,, \ \ 
b_{\La_i} = \boxi \ \  (i=0,1,\cdots,n).$$

\end{example}

\vskip 5mm

\begin{example} \ $A_{2n-1}^{(2)}$ $(n\ge 3)$ : 

\vskip 3mm

\begin{center}
\savebox{\tmppic}{\begin{texdraw}
%\fontsize{9}{9}\selectfont
\small
\drawdim em
\setunitscale 0.5
\nc{\numbox}[1]{\bsegment
\move(-1 -1)\lvec(1 -1)\lvec(1 1)\lvec(-1 1)\lvec(-1 -1)
\htext(0 0){#1}
\esegment}
\move(0 0)\numbox{$1$}
\move(1.5 0)\ravec(5 0)
\htext(3.7 1){$1$}
\move(8 0)\numbox{$2$}
\move(9.5 0)\ravec(5 0)
\htext(11.7 1){$2$}
\htext(16.7 0){$\cdots$}
\move(18.5 0)\ravec(5 0)
\htext(20.5 1){$n\!-\!1$}
\move(25 0)\numbox{$n$}
\move(26.5 0)\ravec(5 0)
\htext(28.7 1){$n$}
\move(33 0)\numbox{$\bar{n}$}
\move(34.5 0)\ravec(5 0)
\htext(36.4 1){$n\!-\!1$}
\htext(41.7 0){$\cdots$}
\move(43.5 0)\ravec(5 0)
\htext(45.7 1){$2$}
\move(50 0)\numbox{$\bar{2}$}
\move(51.5 0)\ravec(5 0)
\htext(53.7 1){$1$}
\move(58 0)\numbox{$\bar{1}$}
%\move(8.5 -1.5)
%\clvec (13 -5)(53 -5)(57.5 -1.5)
%\move(8.8 -1.68)\ravec(-0.45 0.26)
\hcvec{8.5}{57.5}{-1.5}{4.5}{-3.5}
\move(8.5 -1.5)\ahead{-460}{275}
\htext(33 -5.1){$0$}
%\move(0.5 1.5)
%\clvec (5 5)(45 5)(49.5 1.5)
%\move(0.8 1.68)\ravec(-0.45 -0.26)
\hcvec{0.5}{49.5}{1.5}{4.5}{3.5}
\move(0.5 1.5)\ahead{-460}{-275}
\htext(25 5.3){$0$}
\move(-1.1 -6.2)\move(59.1 6.2)
\end{texdraw}}%
%$\B$ : \raisebox{0.3em-0.5\height}{\usebox{\tmppic}}
\usebox{\tmppic}
\end{center}

\vskip 3mm
\noindent
Here, we have
$${b^{\La_0} =
\savebox{\tmppic}{\begin{texdraw}
\small
\drawdim em
\textref h:C v:C
\setunitscale 0.55
\htext(0 0){$1$}
\move(-1 -1)\lvec(-1 1)\lvec(1 1)\lvec(1 -1)\lvec(-1 -1)
\end{texdraw}}%
\raisebox{-0.19\height}{\usebox{\tmppic}}%
}\,, \ 
{b_{\La_0} =
\savebox{\tmppic}{\begin{texdraw}
\small
\drawdim em
\textref h:C v:C
\setunitscale 0.55
\htext(0 0){$\bar{1}$}
\move(-1 -1)\lvec(-1 1)\lvec(1 1)\lvec(1 -1)\lvec(-1 -1)
\end{texdraw}}%
\raisebox{-0.19\height}{\usebox{\tmppic}}%
}\,; \ \ 
{b^{\La_1} =
\savebox{\tmppic}{\begin{texdraw}
\small
\drawdim em
\textref h:C v:C
\setunitscale 0.55
\htext(0 0){$\bar{1}$}
\move(-1 -1)\lvec(-1 1)\lvec(1 1)\lvec(1 -1)\lvec(-1 -1)
\end{texdraw}}%
\raisebox{-0.19\height}{\usebox{\tmppic}}%
}\,, \ 
{b_{\La_1} =
\savebox{\tmppic}{\begin{texdraw}
\small
\drawdim em
\textref h:C v:C
\setunitscale 0.55
\htext(0 0){$1$}
\move(-1 -1)\lvec(-1 1)\lvec(1 1)\lvec(1 -1)\lvec(-1 -1)
\end{texdraw}}%
\raisebox{-0.19\height}{\usebox{\tmppic}}%
}\,. $$
\end{example}

\vskip 5mm
%\newpage
\begin{example} \ $D_{n}^{(1)}$ $(n\ge 4)$ : 

\vskip 3mm

\begin{center}
\savebox{\tmppic}{\begin{texdraw}
\fontsize{9}{9}\selectfont
\drawdim em
\setunitscale 0.5
\nc{\numbox}[1]{\bsegment
\move(-1 -1)\lvec(1 -1)\lvec(1 1)\lvec(-1 1)\lvec(-1 -1)
\htext(0 0){#1}
\esegment}
\move(0 0)
\bsegment
\move(0 0)\numbox{$1$}
\move(1.5 0)\ravec(5 0)
\htext(3.7 1){$1$}
\move(8 0)\numbox{$2$}
\move(9.5 0)\ravec(5 0)
\htext(11.7 1){$2$}
\htext(16.7 0){$\cdots$}
\move(18.5 0)\ravec(5 0)
\htext(20.5 1){$n\!-\!2$}
\move(26.2 0)
\bsegment
\move(-2.2 -1)\lvec(2.2 -1)\lvec(2.2 1)\lvec(-2.2 1)\lvec(-2.2 -1)
\htext(0.1 0){$n\!-\!1$}
\esegment
\move(28 1.5)\ravec(4 3.3)
\move(28 -1.5)\ravec(4 -3.3)
\htext(27.6 3.7){$n\!-\!1$}
\htext(29 -3.5){$n$}
\esegment
\move(67 0)
\bsegment
\move(0 0)\numbox{$\bar{1}$}
\move(-6.5 0)\ravec(5 0)
\htext(-4.3 1){$1$}
\move(-8 0)\numbox{$\bar{2}$}
\move(-14.5 0)\ravec(5 0)
\htext(-12.3 1){$2$}
\htext(-16.7 0){$\cdots$}
\move(-23.5 0)\ravec(5 0)
\htext(-21.3 1){$n\!-\!2$}
\move(-26.2 0)
\bsegment
\move(-2.2 -1)\lvec(2.2 -1)\lvec(2.2 1)\lvec(-2.2 1)\lvec(-2.2 -1)
\htext(0 0){$\overline{n\!-\!1}$}
\esegment
\move(-32 4.8)\ravec(4 -3.3)
\move(-32 -4.8)\ravec(4 3.3)
\htext(-29 3.7){$n$}
\htext(-27.8 -3.5){$n\!-\!1$}
\esegment
\move(33.5 5)\numbox{$n$}
\move(33.5 -5)\numbox{$\bar{n}$}
\hcvec{0.5}{58.5}{1.5}{5}{8}
\move(0.5 1.5)\ahead{-150}{-168}
%\move(0.5 1.5)\ravec(-0.15 -0.168)
\htext(29.5 8.8){$0$}
\hcvec{8.5}{66.5}{-1.5}{5}{-8}
\move(8.5 -1.5)\ahead{-150}{168}
%\move(8.5 -1.5)
%\clvec (13.5 -9.5)(61.5 -9.5)(66.5 -1.5)
%\move(8.5 -1.5)\ravec(-0.15 0.168)
\htext(37.5 -8.6){$0$}
\move(-1.1 -9.5)\move(68.2 9.5)
\end{texdraw}}%
%$\B$ : \raisebox{0.3em-0.5\height}{\usebox{\tmppic}}
\usebox{\tmppic}
\end{center}

\vskip 3mm
\noindent
Here, we have 
\begin{equation*}
\begin{aligned}\mbox{}
& {b^{\La_0} =
\savebox{\tmppic}{\begin{texdraw}
\small
\drawdim em
\textref h:C v:C
\setunitscale 0.55
\htext(0 0){$1$}
\move(-1 -1)\lvec(-1 1)\lvec(1 1)\lvec(1 -1)\lvec(-1 -1)
\end{texdraw}}%
\raisebox{-0.21\height}{\usebox{\tmppic}}%
}\,, \ 
{b_{\La_0} =
\savebox{\tmppic}{\begin{texdraw}
\small
\drawdim em
\textref h:C v:C
\setunitscale 0.55
\htext(0 0){$\bar{1}$}
\move(-1 -1)\lvec(-1 1)\lvec(1 1)\lvec(1 -1)\lvec(-1 -1)
\end{texdraw}}%
\raisebox{-0.21\height}{\usebox{\tmppic}}%
}\,; \ \ 
{b^{\La_1} =
\savebox{\tmppic}{\begin{texdraw}
\small
\drawdim em
\textref h:C v:C
\setunitscale 0.55
\htext(0 0){$\bar{1}$}
\move(-1 -1)\lvec(-1 1)\lvec(1 1)\lvec(1 -1)\lvec(-1 -1)
\end{texdraw}}%
\raisebox{-0.21\height}{\usebox{\tmppic}}%
}\,, \ 
{b_{\La_1} =
\savebox{\tmppic}{\begin{texdraw}
\small
\drawdim em
\textref h:C v:C
\setunitscale 0.55
\htext(0 0){$1$}
\move(-1 -1)\lvec(-1 1)\lvec(1 1)\lvec(1 -1)\lvec(-1 -1)
\end{texdraw}}%
\raisebox{-0.21\height}{\usebox{\tmppic}}%
}\,,\\
& {b^{\La_{n-1}} =
\savebox{\tmppic}{\begin{texdraw}
\small
\drawdim em
\textref h:C v:C
\setunitscale 0.55
\htext(0 0){$n$}
\move(-1 -1)\lvec(-1 1)\lvec(1 1)\lvec(1 -1)\lvec(-1 -1)
\end{texdraw}}%
\raisebox{-0.21\height}{\usebox{\tmppic}}%
}\,, \ 
{b_{\La_{n-1}} =
\savebox{\tmppic}{\begin{texdraw}
\small
\drawdim em
\textref h:C v:C
\setunitscale 0.55
\htext(0 0){$\bar{n}$}
\move(-1 -1)\lvec(-1 1)\lvec(1 1)\lvec(1 -1)\lvec(-1 -1)
\end{texdraw}}%
\raisebox{-0.21\height}{\usebox{\tmppic}}%
}\,; \ \ 
{b^{\La_n} =
\savebox{\tmppic}{\begin{texdraw}
\small
\drawdim em
\textref h:C v:C
\setunitscale 0.55
\htext(0 0){$\bar{n}$}
\move(-1 -1)\lvec(-1 1)\lvec(1 1)\lvec(1 -1)\lvec(-1 -1)
\end{texdraw}}%
\raisebox{-0.21\height}{\usebox{\tmppic}}%
}\,, \ 
{b_{\La_n} =
\savebox{\tmppic}{\begin{texdraw}
\small
\drawdim em
\textref h:C v:C
\setunitscale 0.55
\htext(0 0){$n$}
\move(-1 -1)\lvec(-1 1)\lvec(1 1)\lvec(1 -1)\lvec(-1 -1)
\end{texdraw}}%
\raisebox{-0.21\height}{\usebox{\tmppic}}%
}\,. 
\end{aligned}
\end{equation*}

\end{example}

\vskip 5mm

\begin{example} \  $A_{2n}^{(2)}$ $(n\ge 2)$ : 

\vskip 3mm

\begin{center}
\savebox{\tmppic}{\begin{texdraw}
%\fontsize{7}{7}\selectfont
\drawdim em
\setunitscale 0.6
\nc{\numbox}[1]{\bsegment
\setunitscale 0.52
\move(-1 -1)\lvec(1 -1)\lvec(1 1)\lvec(-1 1)\lvec(-1 -1)
\htext(0 0){#1}
\esegment}
\move(0 0)
\bsegment
\move(0 0)\numbox{$1$}
\move(1.5 0)\ravec(5 0)
\htext(3.7 1){$1$}
\move(8 0)\numbox{$2$}
\move(9.5 0)\ravec(5 0)
\htext(11.7 1){$2$}
\htext(16.7 0){$\cdots$}
\move(18.5 0)\ravec(5 0)
\htext(20.5 1){$n\!-\!2$}
\move(26.2 0)
\bsegment
\setunitscale 0.52
\move(-2.2 -1)\lvec(2.2 -1)\lvec(2.2 1)\lvec(-2.2 1)\lvec(-2.2 -1)
\htext(0.1 0){$n\!-\!1$}
\esegment
\move(28.9 0)\ravec(5 0)\htext(30.9 1){$n\!-\!1$}
\move(35.4 0)\numbox{$n$}
\esegment
\move(0 -8)
\bsegment
\move(0 0)\numbox{$\bar{1}$}
\move(6.5 0)\ravec(-5 0)
\htext(4.3 1){$1$}
\move(8 0)\numbox{$\bar{2}$}
\move(14.5 0)\ravec(-5 0)
\htext(12.3 1){$2$}
\htext(16.7 0){$\cdots$}
\move(23.5 0)\ravec(-5 0)
\htext(21.6 1){$n\!-\!2$}
\move(26.2 0)
\bsegment
\setunitscale 0.52
\move(-2.2 -1)\lvec(2.2 -1)\lvec(2.2 1)\lvec(-2.2 1)\lvec(-2.2 -1)
\htext(0 0){$\overline{n\!-\!1}$}
\esegment
\move(33.9 0)\ravec(-5 0)\htext(32.1 1){$n\!-\!1$}
\move(35.4 0)\numbox{$\bar{n}$}
\esegment
\move(35.4 -1.5)\avec(35.4 -6.5)\htext(36.5 -3.7){$n$}
\htext(-9 -4){$\emptyset$}
\move(-8 -3.3)\avec(-1.5 -0.3)\htext(-5.2 -0.9){$0$}
\move(-1.5 -7.7)\avec(-8 -4.7)\htext(-3.8 -5.4){$0$}
\move(-9.6 -9.1)\move(37.1 1.7)
\end{texdraw}}%
%$\B$ : \raisebox{0.3em-0.5\height}{\usebox{\tmppic}}
\usebox{\tmppic}
\end{center}

\vskip 3mm
\noindent
Here, we have 
$$b^{\La_0} = \emptyset, \ \ b_{\La_0} = \emptyset.$$

\end{example}

\vskip 5mm

\begin{example} \ 
$D_{n+1}^{(2)}$ $(n\ge 2)$ : 

\vskip 3mm

\begin{center}
\savebox{\tmppic}{\begin{texdraw}
\small
\drawdim em
\setunitscale 0.5
\nc{\numbox}[1]{\bsegment
\move(-1 -1)\lvec(1 -1)\lvec(1 1)\lvec(-1 1)\lvec(-1 -1)
\htext(0 0){#1}
\esegment}
\move(0 0)
\bsegment
\move(0 0)\numbox{$1$}
\move(1.5 0)\ravec(5 0)
\htext(3.7 1){$1$}
\move(8 0)\numbox{$2$}
\move(9.5 0)\ravec(5 0)
\htext(11.7 1){$2$}
\htext(16.7 0){$\cdots$}
\move(18.5 0)\ravec(5 0)
\htext(20.5 1){$n\!-\!1$}
\move(25 0)\numbox{$n$}
\esegment
\move(0 -8)
\bsegment
\move(0 0)\numbox{$\bar{1}$}
\move(6.5 0)\ravec(-5 0)
\htext(4.3 1){$1$}
\move(8 0)\numbox{$\bar{2}$}
\move(14.5 0)\ravec(-5 0)
\htext(12.3 1){$2$}
\htext(16.7 0){$\cdots$}
\move(23.5 0)\ravec(-5 0)
\htext(21.5 1){$n\!-\!1$}
\move(25 0)\numbox{$\bar{n}$}
\esegment
\move(33 -4)\numbox{$0$}
\htext(-8 -4){$\emptyset$}
\move(26.5 -0.3)\avec(31.5 -3.3)\htext(29 -0.8){$n$}
\move(31.5 -4.7)\avec(26.5 -7.7)\htext(29 -5.1){$n$}
\move(-7 -3.3)\avec(-1.5 -0.3)\htext(-4.6 -0.7){$0$}
\move(-1.5 -7.7)\avec(-7 -4.7)\htext(-3.5 -5.4){$0$}
\end{texdraw}}%
%$\B$ : \raisebox{-2.8em}{\usebox{\tmppic}}
\usebox{\tmppic}
\end{center}

\vskip 3mm
\noindent
Here, we have 
$$b^{\La_0} = \emptyset, \ \ b_{\La_0} = \emptyset\,; \ \
{b^{\La_n} =
\savebox{\tmppic}{\begin{texdraw}
\small
\drawdim em
\textref h:C v:C
\setunitscale 0.55
\htext(0 0){$0$}
\move(-1 -1)\lvec(-1 1)\lvec(1 1)\lvec(1 -1)\lvec(-1 -1)
\end{texdraw}}%
\raisebox{-0.19\height}{\usebox{\tmppic}}%
}\,, \ \
{b_{\La_n} =
\savebox{\tmppic}{\begin{texdraw}
\small
\drawdim em
\textref h:C v:C
\setunitscale 0.55
\htext(0 0){$0$}
\move(-1 -1)\lvec(-1 1)\lvec(1 1)\lvec(1 -1)\lvec(-1 -1)
\end{texdraw}}%
\raisebox{-0.19\height}{\usebox{\tmppic}}%
}\,. $$

\end{example}

\vskip 5mm

\begin{example} \  $B_{n}^{(1)}$ $(n\ge 3)$ : 

\vskip 3mm

\begin{center}
\savebox{\tmppic}{\begin{texdraw}
\small
\drawdim em
\setunitscale 0.5
\nc{\numbox}[1]{\bsegment
\move(-1 -1)\lvec(1 -1)\lvec(1 1)\lvec(-1 1)\lvec(-1 -1)
\htext(0 0){#1}
\esegment}
\move(0 0)
\bsegment
\move(0 0)\numbox{$1$}
\move(1.5 0)\ravec(5 0)
\htext(3.7 1){$1$}
\move(8 0)\numbox{$2$}
\move(9.5 0)\ravec(5 0)
\htext(11.7 1){$2$}
\htext(16.7 0){$\cdots$}
\move(18.5 0)\ravec(5 0)
\htext(20.5 1){$n\!-\!1$}
\move(25 0)\numbox{$n$}
\move(26.5 0)\ravec(5 0)
\htext(28.7 1){$n$}
\esegment
\move(8 0)
\bsegment
\move(25 0)\numbox{$0$}
\move(26.5 0)\ravec(5 0)
\htext(28.7 1){$n$}
\move(33 0)\numbox{$\bar{n}$}
\move(34.5 0)\ravec(5 0)
\htext(36.4 1){$n\!-\!1$}
\htext(41.7 0){$\cdots$}
\move(43.5 0)\ravec(5 0)
\htext(45.7 1){$2$}
\move(50 0)\numbox{$\bar{2}$}
\move(51.5 0)\ravec(5 0)
\htext(53.7 1){$1$}
\move(58 0)\numbox{$\bar{1}$}
\esegment
\hcvec{8.5}{65.5}{-1.5}{4.5}{-3.5}
\move(8.5 -1.5)\ahead{-460}{275}
\htext(37 -5.1){$0$}
\hcvec{0.5}{57.5}{1.5}{4.5}{3.5}
\move(0.5 1.5)\ahead{-460}{-275}
\htext(29 5.3){$0$}
\move(-1.1 -6.2)\move(59.1 6.2)
\end{texdraw}}%
\usebox{\tmppic}
\end{center}

\vskip 3mm 
\noindent
Here, we have 
\begin{equation*}
\begin{aligned}\mbox{}
& {b^{\La_0} =
\savebox{\tmppic}{\begin{texdraw}
\small
\drawdim em
\textref h:C v:C
\setunitscale 0.55
\htext(0 0){$1$}
\move(-1 -1)\lvec(-1 1)\lvec(1 1)\lvec(1 -1)\lvec(-1 -1)
\end{texdraw}}%
\raisebox{-0.19\height}{\usebox{\tmppic}}%
}\,, \  \ 
{b_{\La_0} =
\savebox{\tmppic}{\begin{texdraw}
\small
\drawdim em
\textref h:C v:C
\setunitscale 0.55
\htext(0 0){$\bar{1}$}
\move(-1 -1)\lvec(-1 1)\lvec(1 1)\lvec(1 -1)\lvec(-1 -1)
\end{texdraw}}%
\raisebox{-0.19\height}{\usebox{\tmppic}}%
}\,; \  \ 
{b^{\La_1} =
\savebox{\tmppic}{\begin{texdraw}
\small
\drawdim em
\textref h:C v:C
\setunitscale 0.55
\htext(0 0){$\bar{1}$}
\move(-1 -1)\lvec(-1 1)\lvec(1 1)\lvec(1 -1)\lvec(-1 -1)
\end{texdraw}}%
\raisebox{-0.19\height}{\usebox{\tmppic}}%
}\,, \ \
{b_{\La_1} =
\savebox{\tmppic}{\begin{texdraw}
\small
\drawdim em
\textref h:C v:C
\setunitscale 0.55
\htext(0 0){$1$}
\move(-1 -1)\lvec(-1 1)\lvec(1 1)\lvec(1 -1)\lvec(-1 -1)
\end{texdraw}}%
\raisebox{-0.19\height}{\usebox{\tmppic}}%
}\,; \\
& {b^{\La_n} =
\savebox{\tmppic}{\begin{texdraw}
\small
\drawdim em
\textref h:C v:C
\setunitscale 0.55
\htext(0 0){$0$}
\move(-1 -1)\lvec(-1 1)\lvec(1 1)\lvec(1 -1)\lvec(-1 -1)
\end{texdraw}}%
\raisebox{-0.19\height}{\usebox{\tmppic}}%
}\,, \  \ 
{b_{\La_n} =
\savebox{\tmppic}{\begin{texdraw}
\small
\drawdim em
\textref h:C v:C
\setunitscale 0.55
\htext(0 0){$0$}
\move(-1 -1)\lvec(-1 1)\lvec(1 1)\lvec(1 -1)\lvec(-1 -1)
\end{texdraw}}%
\raisebox{-0.19\height}{\usebox{\tmppic}}%
}\,. 
\end{aligned}
\end{equation*}

\end{example}

\newpage
\section{Path realization of crystal graphs}

Fix a positive integer $l>0$ and let $\B$ be a perfect crystal of 
level $l$. 
By definition, for any classical dominant integral weight 
$\la \in \oP^+_l$, there exists a unique element 
$b_{\la} \in \B$ such that $\vphi(b_{\la})=\la$.
Set
\begin{equation*}
\mu=\lambda - \wt(b_{\la})=\vep(b_{\la}),
\end{equation*}
and denote by $u_{\mu}$ the highest weight vector of the
crystal graph $B(\mu)$. 
Then, using the fact that $\B$ is perfect, one can show that the vector 
$u_{\mu} \otimes b_{\la}$ is the unique maximal vector 
in $B(\mu) \otimes \B$. 
Moreover, we have:

\vskip 3mm

\begin{thm} {\rm (\cite{KMN1})}
Let $\B$ be a perfect crystal of level $l>0$.
Then for any dominant integral weight $\la \in \oP^+_l$,
there exists a crystal isomorphism 
\begin{equation*}
\Psi: B(\la) \stackrel{\sim}\longrightarrow B(\vep(b_{\la})) \otimes \B
\quad \text{given by} \quad 
u_{\la} \longmapsto u_{\vep(b_{\la})} \otimes b_{\la},
\end{equation*}
where $b_{\la}$ is the unique vector in $\B$ such that 
$\vphi(b_{\la})=\la$. 
\end{thm}

Set 
\begin{equation*}
\la_0 = \la, \quad \la_{k+1}=\vep(b_{\la_k}), 
\end{equation*}
and
\begin{equation*}
b_0 = b_{\la}, \quad b_{k+1}=b_{\la_{k+1}}. 
\end{equation*}
By taking the composition of crystal isomorphism given in Theorem 
4.1 repeatedly, we obtain a sequence of crystal isomorphisms
\begin{equation*}
B(\la) \stackrel{\sim} \longrightarrow 
B(\lambda_1) \otimes \B 
\stackrel{\sim} \longrightarrow B(\la_2) \otimes \B \otimes \B 
\stackrel{\sim} \longrightarrow \cdots\cdots
\end{equation*}
given by 
\begin{equation*}
u_{\la} \longmapsto u_{\la_1} \otimes b_0 \longmapsto u_{\la_2} \otimes 
b_1 \otimes b_0 \longmapsto \cdots \cdots,
\end{equation*}
which yields the infinite sequences 
\begin{equation*}
\text{\bf w}_{\la} 
=(\la_{k})_{k=0}^{\infty}
=(\cdots, \la_{k+1}, \la_k, \cdots, \la_1, \la_0) \quad \text{in} 
\ \ (\oP^{+}_l)^{\infty}
\end{equation*}
and 
\begin{equation*}
\text{\bf p}_{\la} 
=(b_{k})_{k=0}^{\infty}
=\cdots \otimes b_{k+1} \otimes b_k \otimes \cdots 
\otimes b_1 \otimes b_0 \quad \text{in} 
\ \ \B^{\otimes {\infty}}.
\end{equation*}
Thus for each $k\ge 1$, we get a crystal isomorphism 
\begin{equation*}
\Psi_k: B(\la) \stackrel{\sim} \longrightarrow 
B(\la_k) \otimes \B^{\otimes k}
\end{equation*}
given by
\begin{equation*}
u_{\la} \longmapsto u_{\la_k} \otimes b_{k-1} \otimes 
\cdots \otimes b_1 \otimes b_0. 
\end{equation*}
Since $\B$ is perfect, we have 
$\vphi(b_j)=\la_j$ and $\vep(b_j)=\la_{j+1}$.
It follows that the sequences 
$$\text{\bf w}_{\la} 
=(\la_{k})_{k=0}^{\infty}
=(\cdots, \la_{k+1}, \la_k, \cdots, \la_1, \la_0)$$
and 
$$\text{\bf p}_{\la} 
=( b_{k})_{k=0}^{\infty}
=\cdots \otimes b_{k+1} \otimes b_k \otimes \cdots 
\otimes b_1 \otimes b_0$$
are periodic with the same period.
That is, there is a positive integer $N>0$ such that
$\la_{j+N}=\la_j$, $b_{j+N}=b_j$ for all $j=0,1, \cdots, N-1$.

\vskip 3mm

\begin{defi} \hfill 

(a) The sequence $\text{\bf p}_{\la} 
=(b_{k})_{k=0}^{\infty}
=\cdots \otimes b_{k+1} \otimes b_k \otimes \cdots 
\otimes b_1 \otimes b_0$ 
is called the {\bf ground-state path} of weight $\la$.

(b) A {\it $\la$-path} in $\B$ is a sequence 
$\text{\bf p}=(\text{\bf p}(k))_{k=0}^{\infty}
= \cdots \otimes \text{\bf p}(k) \otimes \cdots \otimes 
\text{\bf p}(1) \otimes \text{\bf p}(0)$ such that 
$\text{\bf p}(k) = b_{k}$ for all $k\gg 0$.

\end{defi}

\vskip 3mm

Let $\bP(\la)=\bP(\la, \B)$ be the set of all $\la$-paths in $\B$. 
We define the $U'_q(\g)$-crystal structure on $\bP(\la)$
as follows. 
Let $\text{\bf p} = (\text{\bf p}(k))_{k=0}^{\infty}$ be a 
$\la$-path in $\bP(\la)$ and let $N>0$ be a positive integer
such that $\text{\bf p}(k)=b_k$ for all $k\ge N$. 
For each $i\in I$, we define 
\begin{equation}
\begin{aligned}\mbox{}
& \overline{\text{wt}} (\text{\bf p}) = \la_N + \sum_{k=0}^{N-1} 
\overline{\text{wt}} \text{\bf p}(k), \\
& \eit \text{\bf p} = \cdots \otimes \text{\bf p}(N+1) 
\otimes \eit(\text{\bf p}(N) \otimes \cdots \otimes \text{\bf p}(0)), \\
& \fit \text{\bf p} = \cdots \otimes \text{\bf p}(N+1) 
\otimes \fit(\text{\bf p}(N) \otimes \cdots \otimes \text{\bf p}(0)), \\
& \vep_i(\text{\bf p})=\max(\vep_i(\text{\bf p}')-\vphi_i(b_N), 0), \\
& \vphi_i(\text{\bf p})=\vphi_i(\text{\bf p}')+ 
\max(\vphi_i(b_N)-\vep_i(\text{\bf p}'), 0),
\end{aligned}
\end{equation}
where $\text{\bf p}'=\text{\bf p}(N-1) \otimes \cdots \otimes 
\text{\bf p}(1) \otimes \text{\bf p}(0)$. 

\vskip 3mm 
Then we have the {\it path realization} of the 
classical crystal $B(\la)$:

\vskip 3mm

\begin{thm} {\rm (\cite{KMN1})} \hfill 

{\rm (a)} The maps $\overline{\text{wt}} : \bP(\la) \ra \oP$,
$\eit, \fit : \bP(\la) \ra \bP(\la)  \cup \{0\}$,
$\vep_i, \vphi_i: \bP(\la) \ra \Z$ define a 
classical crystal structure on $\bP(\la)$.

\vskip 2mm

{\rm (b)} There exists an isomorphism of classical crystals 
\begin{equation*}
\Psi: B(\la) \stackrel{\sim} \longrightarrow \bP(\la)
\quad \text{given by} \ \  u_{\la} \longmapsto 
\text{\bf p}_{\la}. 
\end{equation*}

\end{thm}

\vskip 5mm

In the following examples, we list the ground-state paths for 
basic representations of quantum affine algebras and 
illustrate the top part of their crystal graphs for some small ranks. 

\vskip 5mm

\begin{example} $A_n^{(1)}$ $(n\ge 1)$

\vskip 3mm

(a) Ground-state paths

\begin{equation*}
\text{\bf p}_{\La_i} = (\text{\bf p}_{\La_i}(k))_{k=0}^{\infty} =
(\,\cdots\,, 1, 2, \cdots\,,n, 0, 1, 2, \cdots\,, i)
\end{equation*}

\newpage

(b) Crystal graph $\B(\La_0)$ for $A_2^{(1)}$

\vskip 5mm
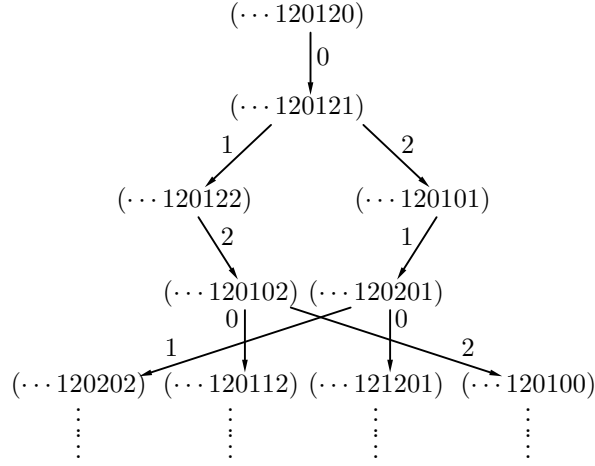
\begin{center}
\begin{texdraw}
\fontsize{10}{10}\selectfont
\drawdim em
\setunitscale 0.5
\htext(-0.5 0){$(\cdots 1 2 0 1 2 0)$}
\htext(-0.5 -7){$(\cdots 1 2 0 1 2 1)$}
\htext(-9 -14){$(\cdots 1 2 0 1 2 2)$}
\htext(9 -14){$(\cdots 1 2 0 1 0 1)$}
\htext(-5.5 -21){$(\cdots 1 2 0 1 0 2)$}
\htext(5.5 -21){$(\cdots 1 2 0 2 0 1)$}
\htext(-17 -28){$(\cdots 1 2 0 2 0 2)$}
\htext(-5.5 -28){$(\cdots 1 2 0 1 1 2)$}
\htext(5.5 -28){$(\cdots 1 2 1 2 0 1)$}
\htext(17 -28){$(\cdots 1 2 0 1 0 0)$}
\htext(-17 -30){$\vdots$}
\htext(-5.5 -30){$\vdots$}
\htext(5.5 -30){$\vdots$}
\htext(17 -30){$\vdots$}
\htext(-17 -32){$\vdots$}
\htext(-5.5 -32){$\vdots$}
\htext(5.5 -32){$\vdots$}
\htext(17 -32){$\vdots$}
\move(0.5 -1.2)\ravec(0 -4.7)\rmove(1 2.8)\htext{$0$}
\move(-2.5 -8.2)\ravec(-5 -4.7)\rmove(1.8 3.2)\htext{$1$}
\move(4.5 -8.2)\ravec(5 -4.7)\rmove(-1.6 3.2)\htext{$2$}
\move(-8 -15.2)\ravec(3 -4.7)\rmove(-0.7 3.2)\htext{$2$}
\move(10 -15.2)\ravec(-3 -4.7)\rmove(0.8 3.2)\htext{$1$}
\move(-4.5 -22.2)\ravec(0 -4.7)\rmove(-0.9 4)\htext{$0$}
\move(6.5 -22.2)\ravec(0 -4.7)\rmove(0.9 4)\htext{$0$}
\move(-1 -22)\ravec(16 -5.2)\rmove(-2.5 2.1)\htext{$2$}
\move(3.5 -22)\ravec(-16 -5.2)\rmove(2.5 2.1)\htext{$1$}
\move(-21.2 -28.8)\move(23 0.8)
\end{texdraw}
\end{center}

\end{example}

\vskip 5mm

\begin{example} $A_{2n-1}^{(2)}$ $(n\ge 3)$

\vskip 3mm

(a) Ground-state paths

\begin{equation*}
\begin{aligned}\mbox{}
& \text{\bf p}_{\La_0} = (\text{\bf p}_{\La_0}(k))_{k=0}^{\infty} =
(\,\cdots\,,1, \bar 1, 1, \bar 1, 1, \bar 1), \\
& \text{\bf p}_{\La_1} = (\text{\bf p}_{\La_1}(k))_{k=0}^{\infty} =
(\,\cdots\,, \bar 1, 1, \bar 1, 1, \bar 1, 1). 
\end{aligned}
\end{equation*}

\vskip 5mm

(b) Crystal graph $\B(\La_0)$ for $A_5^{(2)}$

\vskip 5mm 

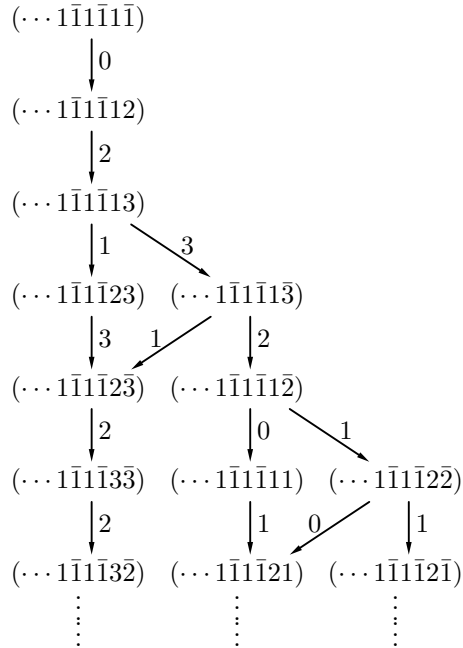
\begin{center}
\begin{texdraw}
\fontsize{10}{10}\selectfont
\drawdim em
\setunitscale 0.5
\htext(0 0){$(\cdots 1 \bar{1} 1 \bar{1} 1 \bar{1})$}
\htext(0 -7){$(\cdots 1 \bar{1} 1 \bar{1} 1 2)$}
\htext(0 -14){$(\cdots 1 \bar{1} 1 \bar{1} 1 3)$}
\htext(0 -21){$(\cdots 1 \bar{1} 1 \bar{1} 2 3)$}
\htext(0 -28){$(\cdots 1 \bar{1} 1 \bar{1} 2 \bar{3})$}
\htext(0 -35){$(\cdots 1 \bar{1} 1 \bar{1} 3 \bar{3})$}
\htext(0 -42){$(\cdots 1 \bar{1} 1 \bar{1} 3 \bar{2})$}
\htext(12 -21){$(\cdots 1 \bar{1} 1 \bar{1} 1 \bar{3})$}
\htext(12 -28){$(\cdots 1 \bar{1} 1 \bar{1} 1 \bar{2})$}
\htext(12 -35){$(\cdots 1 \bar{1} 1 \bar{1} 1 1)$}
\htext(12 -42){$(\cdots 1 \bar{1} 1 \bar{1} 2 1)$}
\htext(24 -35){$(\cdots 1 \bar{1} 1 \bar{1} 2 \bar{2})$}
\htext(24 -42){$(\cdots 1 \bar{1} 1 \bar{1} 2 \bar{1})$}
\htext(0 -44){$\vdots$}
\htext(12 -44){$\vdots$}
\htext(24 -44){$\vdots$}
\htext(0 -46){$\vdots$}
\htext(12 -46){$\vdots$}
\htext(24 -46){$\vdots$}
\move(1 -1.5)\ravec(0 -4)\rmove(1 2.3)\htext{$0$}
\move(1 -8.5)\ravec(0 -4)\rmove(1 2.3)\htext{$2$}
\move(1 -15.5)\ravec(0 -4)\rmove(1 2.3)\htext{$1$}
\move(1 -22.5)\ravec(0 -4)\rmove(1 2.3)\htext{$3$}
\move(1 -29.5)\ravec(0 -4)\rmove(1 2.3)\htext{$2$}
\move(1 -36.5)\ravec(0 -4)\rmove(1 2.3)\htext{$2$}
\move(4 -15.5)\ravec(6 -4)\rmove(-1.8 2.3)\htext{$3$}
\move(10 -22.5)\ravec(-6 -4)\rmove(1.9 2.3)\htext{$1$}
\move(13 -22.5)\ravec(0 -4)\rmove(1 2.3)\htext{$2$}
\move(13 -29.5)\ravec(0 -4)\rmove(1 2.3)\htext{$0$}
\move(13 -36.5)\ravec(0 -4)\rmove(1 2.3)\htext{$1$}
\move(16 -29.5)\ravec(6 -4)\rmove(-1.8 2.3)\htext{$1$}
\move(22 -36.5)\ravec(-6 -4)\rmove(1.9 2.3)\htext{$0$}
\move(25 -36.5)\ravec(0 -4)\rmove(1 2.3)\htext{$1$}
\end{texdraw}
\end{center}

\end{example}

\newpage

\begin{example}
$D_{n}^{(1)}$ $(n\ge 4)$

\vskip 3mm

(a) Ground-state paths

\begin{equation*}
\begin{aligned}\mbox{}
& \text{\bf p}_{\La_0} = (\text{\bf p}_{\La_0}(k))_{k=0}^{\infty} =
(\,\cdots\,,1, \bar 1, 1, \bar 1, 1, \bar 1), \\
& \text{\bf p}_{\La_1} = (\text{\bf p}_{\La_1}(k))_{k=0}^{\infty} =
(\,\cdots\,, \bar 1, 1, \bar 1, 1, \bar 1, 1), \\
& \text{\bf p}_{\La_{n-1}} = (\text{\bf p}_{\La_{n-1}}(k))_{k=0}^{\infty} =
(\,\cdots\,,n, \bar n, n, \bar n, n, \bar n), \\
& \text{\bf p}_{\La_n} = (\text{\bf p}_{\La_n}(k))_{k=0}^{\infty} =
(\,\cdots\,, \bar n, n, \bar n, n, \bar n, n), \\
\end{aligned}
\end{equation*}

\vskip 5mm

(b) Crystal graph $\B(\La_0)$ for $D_4^{(1)}$

\vskip 8mm

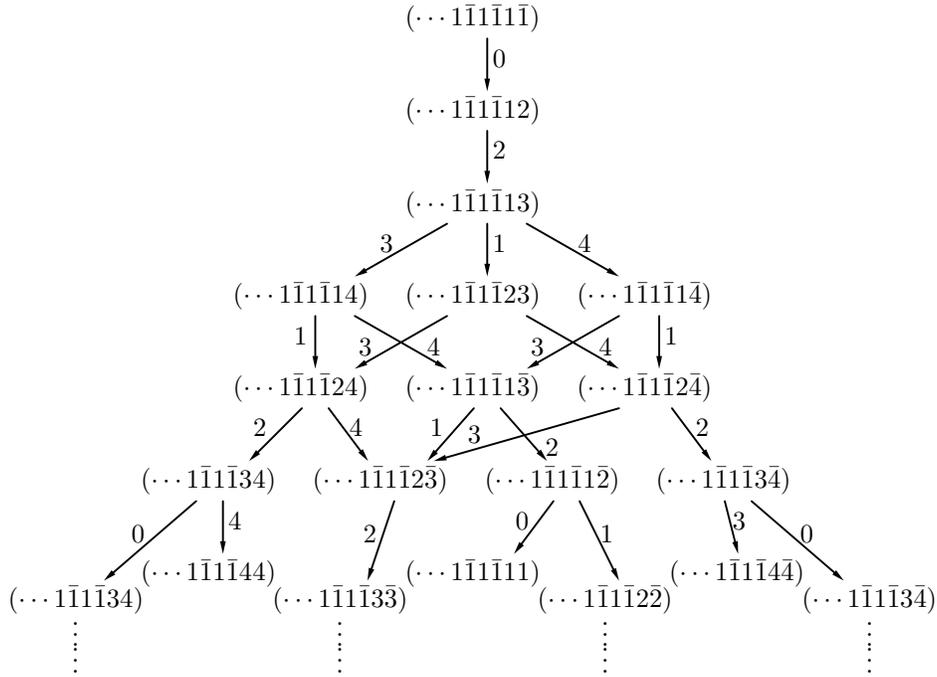
\begin{center}
\begin{texdraw}
\fontsize{10}{10}\selectfont
\drawdim em
\setunitscale 0.5
\htext(0 0){$(\cdots 1 \bar{1} 1 \bar{1} 1 \bar{1})$}
\htext(0 -7){$(\cdots 1 \bar{1} 1 \bar{1} 1 2)$}
\htext(0 -14){$(\cdots 1 \bar{1} 1 \bar{1} 1 3)$}
\htext(-13 -21){$(\cdots 1 \bar{1} 1 \bar{1} 1 4)$}
\htext(0 -21){$(\cdots 1 \bar{1} 1 \bar{1} 2 3)$}
\htext(13 -21){$(\cdots 1 \bar{1} 1 \bar{1} 1 \bar{4})$}
\htext(-13 -28){$(\cdots 1 \bar{1} 1 \bar{1} 2 4)$}
\htext(0 -28){$(\cdots 1 \bar{1} 1 \bar{1} 1 \bar{3})$}
\htext(13 -28){$(\cdots 1 \bar{1} 1 \bar{1} 2 \bar{4})$}
\htext(-20 -35){$(\cdots 1 \bar{1} 1 \bar{1} 3 4)$}
\htext(-7 -35){$(\cdots 1 \bar{1} 1 \bar{1} 2 \bar{3})$}
\htext(6 -35){$(\cdots 1 \bar{1} 1 \bar{1} 1 \bar{2})$}
\htext(19 -35){$(\cdots 1 \bar{1} 1 \bar{1} 3 \bar{4})$}
\htext(-30 -44){$(\cdots 1 \bar{1} 1 \bar{1} 3 4)$}
\htext(-20 -42){$(\cdots 1 \bar{1} 1 \bar{1} 4 4)$}
\htext(-10 -44){$(\cdots 1 \bar{1} 1 \bar{1} 3 \bar{3})$}
\htext(0 -42){$(\cdots 1 \bar{1} 1 \bar{1} 1 1)$}
\htext(10 -44){$(\cdots 1 \bar{1} 1 \bar{1} 2 \bar{2})$}
\htext(20 -42){$(\cdots 1 \bar{1} 1 \bar{1} 4 \bar{4})$}
\htext(30 -44){$(\cdots 1 \bar{1} 1 \bar{1} 3 \bar{4})$}
\htext(30 -46){$\vdots$}
\htext(10 -46){$\vdots$}
\htext(-10 -46){$\vdots$}
\htext(-30 -46){$\vdots$}
\htext(30 -48){$\vdots$}
\htext(10 -48){$\vdots$}
\htext(-10 -48){$\vdots$}
\htext(-30 -48){$\vdots$}
\move(1 -1.5)\ravec(0 -4)\rmove(1 2.5)\htext{$0$}
\move(1 -8.5)\ravec(0 -4)\rmove(1 2.5)\htext{$2$}
\move(1 -15.5)\ravec(0 -4)\rmove(1 2.5)\htext{$1$}
\move(-2 -15.5)\ravec(-7 -4)\rmove(2.5 2.5)\htext{$3$}
\move(4 -15.5)\ravec(7 -4)\rmove(-2.5 2.5)\htext{$4$}
\move(-12 -22.5)\ravec(0 -4)\rmove(-1 2.5)\htext{$1$}
\move(14 -22.5)\ravec(0 -4)\rmove(1 2.5)\htext{$1$}
\move(-9 -22.5)\ravec(7 -4)\rmove(-0.9 1.7)\htext{$4$}
\move(11 -22.5)\ravec(-7 -4)\rmove(0.9 1.7)\htext{$3$}
\move(4 -22.5)\ravec(7 -4)\rmove(-0.9 1.7)\htext{$4$}
\move(-2 -22.5)\ravec(-7 -4)\rmove(0.9 1.7)\htext{$3$}
\move(15 -29.5)\ravec(3 -4)\rmove(-0.6 2.5)\htext{$2$}
\move(11 -29.5)\ravec(-14 -4)\rmove(3.2 2)\htext{$3$}
\move(2 -29.5)\ravec(3.5 -4)\rmove(0.5 1)\htext{$2$}
\move(0 -29.5)\ravec(-3.5 -4)\rmove(0.8 2.5)\htext{$1$}
\move(-11 -29.5)\ravec(3 -4)\rmove(-0.8 2.5)\htext{$4$}
\move(-13 -29.5)\ravec(-4 -4)\rmove(0.9 2.5)\htext{$2$}
\move(-21 -36.5)\ravec(-7 -6)\rmove(2.7 3.5)\htext{$0$}
\move(-19 -36.5)\ravec(0 -4)\rmove(1 2.5)\htext{$4$}
\move(-6 -36.5)\ravec(-2 -6)\rmove(0.3 3.5)\htext{$2$}
\move(6 -36.5)\ravec(-3 -4)\rmove(0.7 2.5)\htext{$0$}
\move(8 -36.5)\ravec(3 -6)\rmove(-0.8 3.5)\htext{$1$}
\move(19 -36.5)\ravec(1 -4)\rmove(0.2 2.5)\htext{$3$}
\move(21 -36.5)\ravec(7 -6)\rmove(-2.7 3.5)\htext{$0$}
\end{texdraw}
\end{center}

\end{example}

\vskip 1cm 

\begin{example}
$A_{2n}^{(2)}$ $(n\ge 2)$

\vskip 3mm

(a) Ground-state paths

\begin{equation*}
\text{\bf p}_{\La_0} = (\text{\bf p}_{\La_0}(k))_{k=0}^{\infty} =
(\,\cdots\,,\emptyset, \emptyset, \emptyset, \emptyset, 
\emptyset, \emptyset)
\end{equation*}

\vskip 3mm

(b) Crystal graph $\B(\La_0)$ for $A_4^{(2)}$

\vskip 3mm

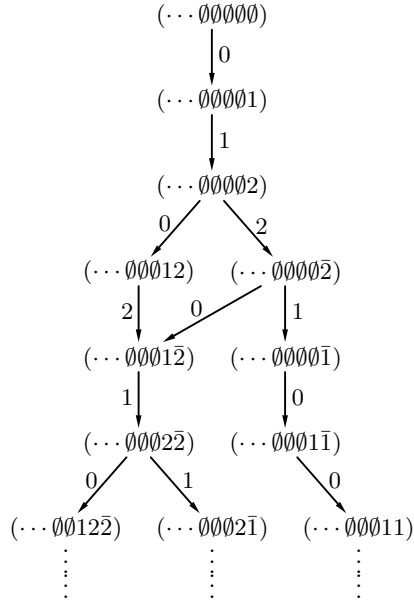
\begin{center}
\begin{texdraw}
%\small
\fontsize{9}{9}\selectfont
\drawdim em
\setunitscale 0.5
\htext(0 35){$(\cdots \emptyset \emptyset \emptyset \emptyset \emptyset)$}
\htext(0 28){$(\cdots \emptyset \emptyset \emptyset \emptyset 1)$}
\htext(0 21){$(\cdots \emptyset \emptyset \emptyset \emptyset 2)$}
\htext(-6 14){$(\cdots \emptyset \emptyset \emptyset 1 2)$}
\htext(6 14){$(\cdots \emptyset \emptyset \emptyset \emptyset \bar{2})$}
\htext(-6 7){$(\cdots \emptyset \emptyset \emptyset 1 \bar{2})$}
\htext(6 7){$(\cdots \emptyset \emptyset \emptyset \emptyset \bar{1}$)}
\htext(-6 0){$(\cdots \emptyset \emptyset \emptyset 2 \bar{2})$}
\htext(6 0){$(\cdots \emptyset \emptyset \emptyset 1 \bar{1})$}
\htext(-12 -7){$(\cdots \emptyset \emptyset 1 2 \bar{2})$}
\htext(0 -7){$(\cdots \emptyset \emptyset \emptyset 2 \bar{1})$}
\htext(12 -7){$(\cdots \emptyset \emptyset \emptyset 1 1)$}
\move(0 33.8)\ravec(0 -4.6)\rmove(1 2.6)\htext{$0$}
\move(0 26.8)\ravec(0 -4.6)\rmove(1 2.6)\htext{$1$}
\move(-1 19.8)\ravec(-4 -4.6)\rmove(1.1 2.8)\htext{$0$}
\move(1 19.8)\ravec(4 -4.6)\rmove(-1 2.6)\htext{$2$}
\move(-6 12.8)\ravec(0 -4.6)\rmove(-1 2.6)\htext{$2$}
\move(4 12.8)\ravec(-8 -4.6)\rmove(2.7 2.8)\htext{$0$}
\move(6 12.8)\ravec(0 -4.6)\rmove(1 2.6)\htext{$1$}
\move(-6 5.8)\ravec(0 -4.6)\rmove(-1 2.6)\htext{$1$}
\move(6 5.8)\ravec(0 -4.6)\rmove(1 2.6)\htext{$0$}
\move(-5 -1.2)\ravec(4 -4.6)\rmove(-1 2.8)\htext{$1$}
\move(7 -1.2)\ravec(4 -4.6)\rmove(-1 2.8)\htext{$0$}
\move(-7 -1.2)\ravec(-4 -4.6)\rmove(1.1 2.8)\htext{$0$}
\htext(-12 -9){$\vdots$}
\htext(0 -9){$\vdots$}
\htext(12 -9){$\vdots$}
\htext(-12 -11){$\vdots$}
\htext(0 -11){$\vdots$}
\htext(12 -11){$\vdots$}
\move(17 -10.3)\move(-17 36.2)
\end{texdraw}
\end{center}

\end{example}

\vskip 5mm

\begin{example}
$D_{n+1}^{(2)}$ ($n\geq2$).

\vskip 2mm

(a) Ground-state paths

\begin{equation*}
\begin{aligned}\mbox{}
& \text{\bf p}_{\La_0} = (\text{\bf p}_{\La_0}(k))_{k=0}^{\infty} =
(\,\cdots\,,\emptyset, \emptyset, \emptyset, \emptyset, 
\emptyset, \emptyset) \\
& \text{\bf p}_{\La_n} = (\text{\bf p}_{\La_n}(k))_{k=0}^{\infty} =
(\,\cdots\,,0, 0, 0, 0, 0, 0) \\
\end{aligned}
\end{equation*}

\vskip 3mm
(b) Crystal graph $\B(\La_0)$ for $D_3^{(2)}$

\vskip 3mm

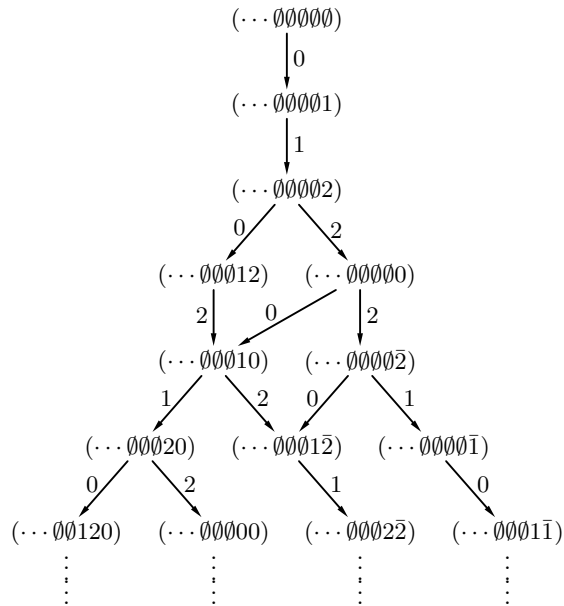
\begin{center}
\begin{texdraw}
%\small
\fontsize{9}{9}\selectfont
\drawdim em
\setunitscale 0.5
\htext(0 35){$(\cdots \emptyset \emptyset \emptyset \emptyset \emptyset)$}
\htext(0 28){$(\cdots \emptyset \emptyset \emptyset \emptyset 1)$}
\htext(0 21){$(\cdots \emptyset \emptyset \emptyset \emptyset 2)$}
\htext(-6 14){$(\cdots \emptyset \emptyset \emptyset 1 2)$}
\htext(6 14){$(\cdots \emptyset \emptyset \emptyset \emptyset 0)$}
\htext(-6 7){$(\cdots \emptyset \emptyset \emptyset 1 0)$}
\htext(6 7){$(\cdots \emptyset \emptyset \emptyset \emptyset \bar{2})$}
\htext(-12 0){$(\cdots \emptyset \emptyset \emptyset 2 0)$}
\htext(0 0){$(\cdots \emptyset \emptyset \emptyset 1 \bar{2})$}
\htext(12 0){$(\cdots \emptyset \emptyset \emptyset \emptyset \bar{1})$}
\htext(-18 -7){$(\cdots \emptyset \emptyset  1 2 0)$}
\htext(-6 -7){$(\cdots \emptyset \emptyset \emptyset 0 0)$}
\htext(6 -7){$(\cdots \emptyset \emptyset \emptyset 2 \bar{2})$}
\htext(18 -7){$(\cdots \emptyset \emptyset \emptyset 1 \bar{1})$}
\move(0 33.8)\ravec(0 -4.6)\rmove(1 2.6)\htext{$0$}
\move(0 26.8)\ravec(0 -4.6)\rmove(1 2.6)\htext{$1$}
\move(-1 19.8)\ravec(-4 -4.6)\rmove(1.1 2.8)\htext{$0$}
\move(1 19.8)\ravec(4 -4.6)\rmove(-1 2.6)\htext{$2$}
\move(-6 12.8)\ravec(0 -4.6)\rmove(-1 2.6)\htext{$2$}
\move(4 12.8)\ravec(-8 -4.6)\rmove(2.7 2.8)\htext{$0$}
\move(6 12.8)\ravec(0 -4.6)\rmove(1 2.6)\htext{$2$}
\move(-7 5.8)\ravec(-4 -4.6)\rmove(1.1 2.8)\htext{$1$}
\move(-5 5.8)\ravec(4 -4.6)\rmove(-1 2.8)\htext{$2$}
\move(5 5.8)\ravec(-4 -4.6)\rmove(1.1 2.8)\htext{$0$}
\move(7 5.8)\ravec(4 -4.6)\rmove(-1 2.8)\htext{$1$}
\move(-13 -1.2)\ravec(-4 -4.6)\rmove(1.1 2.8)\htext{$0$}
\move(-11 -1.2)\ravec(4 -4.6)\rmove(-1 2.8)\htext{$2$}
\move(1 -1.2)\ravec(4 -4.6)\rmove(-1 2.8)\htext{$1$}
\move(13 -1.2)\ravec(4 -4.6)\rmove(-1 2.8)\htext{$0$}
\htext(-18 -9.2){$\vdots$}
\htext(-6 -9.2){$\vdots$}
\htext(6 -9.2){$\vdots$}
\htext(18 -9.2){$\vdots$}
\htext(-18 -11.2){$\vdots$}
\htext(-6 -11.2){$\vdots$}
\htext(6 -11.2){$\vdots$}
\htext(18 -11.2){$\vdots$}
\move(-22.8 -10.6)\move(22.8 36.2)
\end{texdraw}
\end{center}

\end{example}

\newpage

\begin{example}
$B_{n}^{(1)}$ $(n\ge 3)$

\vskip 3mm

(a) Ground-state paths

\begin{equation*}
\begin{aligned}\mbox{}
& \text{\bf p}_{\La_0} = (\text{\bf p}_{\La_0}(k))_{k=0}^{\infty} =
(\,\cdots\,, 1, \bar 1, 1, \bar 1, 1, \bar 1), \\
& \text{\bf p}_{\La_1} = (\text{\bf p}_{\La_1}(k))_{k=0}^{\infty} =
(\, \cdots\,, \bar 1, 1, \bar 1, 1, \bar 1, 1), \\
& \text{\bf p}_{\La_{n}} = (\text{\bf p}_{\La_{n}}(k))_{k=0}^{\infty} =
(\, \cdots\,, 0,0,0,0,0,0)
\end{aligned}
\end{equation*}

\vskip 5mm

(b) Crystal graph $\B(\La_0)$ for $B_3^{(1)}$

\vskip 5mm

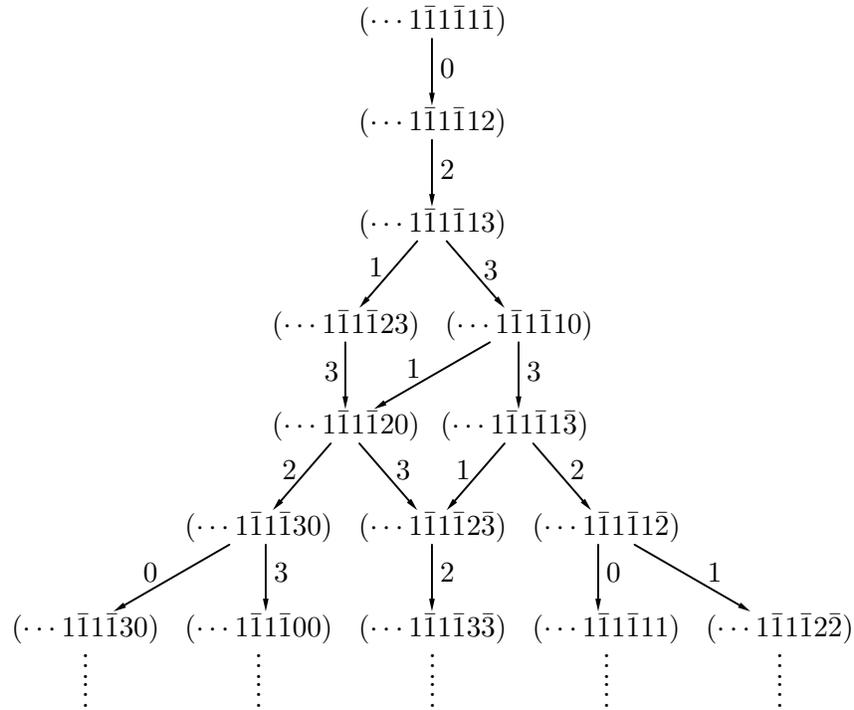
\begin{center}
\begin{texdraw}
\small
\drawdim em
\setunitscale 0.5
\htext(0 35){$(\cdots 1 \bar{1} 1 \bar{1} 1 \bar{1})$}
\htext(0 28){$(\cdots 1 \bar{1} 1 \bar{1} 1 2)$}
\htext(0 21){$(\cdots 1 \bar{1} 1 \bar{1} 1 3)$}
\htext(-6 14){$(\cdots 1 \bar{1} 1 \bar{1} 2 3)$}
\htext(6 14){$(\cdots 1 \bar{1} 1 \bar{1} 1 0)$}
\htext(-6 7){$(\cdots 1 \bar{1} 1 \bar{1} 2 0)$}
\htext(6 7){$(\cdots 1 \bar{1} 1 \bar{1} 1 \bar{3})$ }
\htext(-12 0){$(\cdots 1 \bar{1} 1 \bar{1} 3 0)$}
\htext(0 0){$(\cdots 1 \bar{1} 1 \bar{1} 2 \bar{3})$}
\htext(12 0){$(\cdots 1 \bar{1} 1 \bar{1} 1 \bar{2})$}
\htext(-24 -7){$(\cdots 1 \bar{1} 1 \bar{1} 3 0)$}
\htext(-12 -7){$(\cdots 1 \bar{1} 1 \bar{1} 0 0)$}
\htext(0 -7){$(\cdots 1 \bar{1} 1 \bar{1} 3 \bar{3})$}
\htext(12 -7){$(\cdots 1 \bar{1} 1 \bar{1} 1 1)$}
\htext(24 -7){$(\cdots 1 \bar{1} 1 \bar{1} 2 \bar{2})$}
\move(0 33.8)\ravec(0 -4.6)\rmove(1 2.6)\htext{$0$}
\move(0 26.8)\ravec(0 -4.6)\rmove(1 2.6)\htext{$2$}
\move(-1 19.8)\ravec(-4 -4.6)\rmove(1.1 2.8)\htext{$1$}
\move(1 19.8)\ravec(4 -4.6)\rmove(-1 2.6)\htext{$3$}
\move(-6 12.8)\ravec(0 -4.6)\rmove(-1 2.6)\htext{$3$}
\move(4 12.8)\ravec(-8 -4.6)\rmove(2.7 2.8)\htext{$1$}
\move(6 12.8)\ravec(0 -4.6)\rmove(1 2.6)\htext{$3$}
\move(-7 5.8)\ravec(-4 -4.6)\rmove(1.1 2.8)\htext{$2$}
\move(-5 5.8)\ravec(4 -4.6)\rmove(-1 2.8)\htext{$3$}
\move(5 5.8)\ravec(-4 -4.6)\rmove(1.1 2.8)\htext{$1$}
\move(7 5.8)\ravec(4 -4.6)\rmove(-1 2.8)\htext{$2$}
\move(-14 -1.2)\ravec(-8 -4.6)\rmove(2.5 2.7)\htext{$0$}
\move(-11.5 -1.2)\ravec(0 -4.6)\rmove(1 2.7)\htext{$3$}
\move(0 -1.2)\ravec(0 -4.6)\rmove(1 2.7)\htext{$2$}
\move(11.5 -1.2)\ravec(0 -4.6)\rmove(1 2.7)\htext{$0$}
\move(14 -1.2)\ravec(8 -4.6)\rmove(-2.5 2.7)\htext{$1$}
\htext(-24 -9){$\vdots$}
\htext(-12 -9){$\vdots$}
\htext(0 -9){$\vdots$}
\htext(12 -9){$\vdots$}
\htext(24 -9){$\vdots$}
\htext(-24 -11){$\vdots$}
\htext(-12 -11){$\vdots$}
\htext(0 -11){$\vdots$}
\htext(12 -11){$\vdots$}
\htext(24 -11){$\vdots$}
\move(-29.4 -10.8)\move(29.3 36.2)
\end{texdraw}
\end{center}

\vskip 3mm
\newpage 

(c) Crystal graph $\B(\La_3)$ for $B_3^{(1)}$

\vskip 3mm

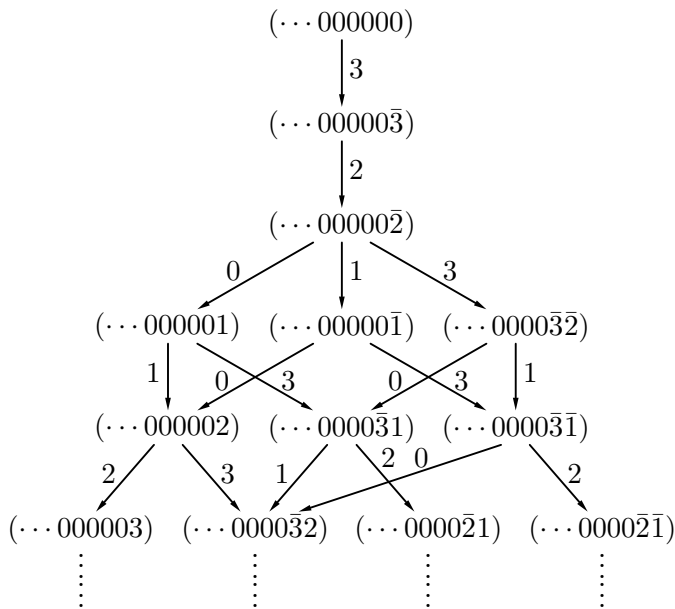
\begin{center}
\begin{texdraw}
\small
\drawdim em
\setunitscale 0.5
\htext(0 35){$(\cdots 0 0 0 0 0 0)$}
\htext(0 28){$(\cdots 0 0 0 0 0 \bar{3})$}
\htext(0 21){$(\cdots 0 0 0 0 0 \bar{2})$}
\htext(-12 14){$(\cdots 0 0 0 0 0 1)$}
\htext(0 14){$(\cdots 0 0 0 0 0 \bar{1})$}
\htext(12 14){$(\cdots 0 0 0 0 \bar{3} \bar{2})$}
\htext(-12 7){$(\cdots 0 0 0 0 0 2)$}
\htext(0 7){$(\cdots 0 0 0 0 \bar{3} 1)$}
\htext(12 7){$(\cdots 0 0 0 0 \bar{3} \bar{1})$}
\htext(-18 0){$(\cdots 0 0 0 0 0 3)$}
\htext(-6 0){$(\cdots 0 0 0 0 \bar{3} 2)$}
\htext(6 0){$(\cdots 0 0 0 0 \bar{2} 1)$}
\htext(18 0){$(\cdots 0 0 0 0 \bar{2} \bar{1})$}
\move(0 33.8)\ravec(0 -4.6)\rmove(1 2.7)\htext{$3$}
\move(0 26.8)\ravec(0 -4.6)\rmove(1 2.7)\htext{$2$}
\move(-2 19.8)\ravec(-8 -4.6)\rmove(2.5 2.7)\htext{$0$}
\move(0 19.8)\ravec(0 -4.6)\rmove(1 2.7)\htext{$1$}
\move(2 19.8)\ravec(8 -4.6)\rmove(-2.5 2.7)\htext{$3$}
\move(-12 12.8)\ravec(0 -4.6)\rmove(-1 2.7)\htext{$1$}
\move(-10 12.8)\ravec(8 -4.6)\rmove(-1.7 2.1)\htext{$3$}
\move(-2 12.8)\ravec(-8 -4.6)\rmove(1.7 2.1)\htext{$0$}
\move(2 12.8)\ravec(8 -4.6)\rmove(-1.7 2.1)\htext{$3$}
\move(10 12.8)\ravec(-8 -4.6)\rmove(1.7 2.1)\htext{$0$}
\move(12 12.8)\ravec(0 -4.6)\rmove(1 2.7)\htext{$1$}
\move(-13 5.8)\ravec(-4 -4.6)\rmove(0.9 2.7)\htext{$2$}
\move(-11 5.8)\ravec(4 -4.6)\rmove(-0.9 2.7)\htext{$3$}
\move(-1 5.8)\ravec(-4 -4.6)\rmove(0.9 2.7)\htext{$1$}
\move(1 5.8)\ravec(4 -4.6)\rmove(-1.8 3.8)\htext{$2$}
\move(11 5.8)\ravec(-14 -4.6)\rmove(8.5 3.8)\htext{$0$}
\move(13 5.8)\ravec(4 -4.6)\rmove(-0.9 2.7)\htext{$2$}
\htext(-18 -2){$\vdots$}
\htext(-6 -2){$\vdots$}
\htext(6 -2){$\vdots$}
\htext(18 -2){$\vdots$}
\htext(-18 -4){$\vdots$}
\htext(-6 -4){$\vdots$}
\htext(6 -4){$\vdots$}
\htext(18 -4){$\vdots$}
\move(-23.1 -3.6)\move(23.1 36.2)
\end{texdraw}
\end{center}

\end{example}

\vskip 1cm

%%%%%%%%%%%%%%%%%%%%%%%%%%%%%%%%%%%%%%%%%%%%%%%%%%%%%%%%%%%%%%%%%%%%%%%%%%%

\section{Young walls}

The purpose of this paper is to give a realization of crystal graph
$B(\la)$ for the basic representations of quantum affine algebras 
in terms of new combinatorial objects which we call the 
\emph{Young walls}.
In this section, we explain the notion of Young walls.
The Young walls are built of colored blocks of three
different shapes:

\vskip 5mm
\hskip 1cm 
\begin{tabular}{rcl}
\raisebox{-0.35em}[1.25em][0.25em]{%
\begin{texdraw}
\drawdim em
\setunitscale 0.15
\textref h:C v:C
\move(0 0)\lvec(10 0)\lvec(10 10)\lvec(0 10)\lvec(0 0)
\move(10 0)\lvec(15 5)\lvec(15 15)\lvec(5 15)\lvec(0 10)
\move(10 10)\lvec(15 15)
\end{texdraw}%
} & : &
unit width, unit height, unit thickness,\\

\\

\raisebox{-0.35em}[1.25em][0.25em]{%
\begin{texdraw}
\drawdim em
\setunitscale 0.15
\textref h:C v:C
\move(0 0)\lvec(10 0)\lvec(10 10)\lvec(0 10)\lvec(0 0)
\move(10 0)\lvec(12.5 2.5)\lvec(12.5 12.5)\lvec(2.5 12.5)\lvec(0 10)
\move(10 10)\lvec(12.5 12.5)
\end{texdraw}%
} & : &
unit width, unit height, half-unit thickness,\\

\\

\raisebox{-0.35em}[1.25em][0.25em]{%
\begin{texdraw}
\drawdim em
\setunitscale 0.15
\textref h:C v:C
\move(0 0)\lvec(10 0)\lvec(10 5)\lvec(0 5)\lvec(0 0)
\move(10 0)\lvec(15 5)\lvec(15 10)\lvec(5 10)\lvec(0 5)
\move(10 5)\lvec(15 10)
\end{texdraw}%
} & : &
unit width, half-unit height, unit thickness.
\end{tabular}

\vskip 5mm

With these colored blocks, we will build a wall of thickness
less than or equal to one unit which extends 
infinitely to the left like playing with the LEGO blocks.
Given a dominant integral weight $\la$ of level 1, we fix a frame
called the {\it ground-state wall} of weight $\la$,
and build the Young walls on this frame.
For each type of classical quantum affine algebras, we use different sets of
colored blocks and ground-state walls, as are described in the 
following. 

%\vskip 5mm

\newpage

(a)  $A_n^{(1)}$ ($n\geq1$) 

\vskip 5mm
\hskip 5mm 
\raisebox{-0.33\height}[0.69\height][0.35\height]{%
\begin{texdraw}
\drawdim em
\setunitscale 0.15
\move(0 0)\lvec(10 0)\lvec(10 10)\lvec(0 10)\lvec(0 0)
\move(10 0)\lvec(15 5)\lvec(15 15)\lvec(5 15)\lvec(0 10)
\move(10 10)\lvec(15 15)
\htext(5 5){$j$}
\end{texdraw}%
}
\quad ($j=0,1,\cdots,n$) : unit width, unit height, unit thickness,\\

\vskip 3mm

\quad $Y_{\La_i} =
\raisebox{-0.5\height}{\begin{texdraw}
\drawdim em
\setunitscale 0.15
\move(-52 0)\lvec(0 0)\lvec(5 5)\lvec(-47 5)
\move(-50 0)\rlvec(5 5)
\move(-40 0)\rlvec(5 5)
\move(-30 0)\rlvec(5 5)
\move(-20 0)\rlvec(5 5)
\move(-10 0)\rlvec(5 5)
\end{texdraw}
} $
\quad ($i=0,1,\cdots,n$).

\vskip 7mm

(b) $A_{2n-1}^{(2)}$ ($n\geq3$)

\vskip 5mm 

\hskip 5mm 
\raisebox{-0.33\height}[0.69\height][0.35\height]{%
\begin{texdraw}
\drawdim em
\setunitscale 0.15
\move(0 0)\lvec(10 0)\lvec(10 10)\lvec(0 10)\lvec(0 0)
\move(10 0)\lvec(12.5 2.5)\lvec(12.5 12.5)\lvec(2.5 12.5)\lvec(0 10)
\move(10 10)\lvec(12.5 12.5)
\htext(5 5){$0$}
\end{texdraw}%
}\,,
\raisebox{-0.33\height}[0.69\height][0.35\height]{%
\begin{texdraw}
\drawdim em
\setunitscale 0.15
\move(0 0)\lvec(10 0)\lvec(10 10)\lvec(0 10)\lvec(0 0)
\move(10 0)\lvec(12.5 2.5)\lvec(12.5 12.5)\lvec(2.5 12.5)\lvec(0 10)
\move(10 10)\lvec(12.5 12.5)
\htext(5 5){$1$}
\end{texdraw}%
}\ :
unit width, unit height, half-unit thickness,\\

\vskip 3mm

\hskip 5mm 
\raisebox{-0.33\height}[0.69\height][0.35\height]{%
\begin{texdraw}
\drawdim em
\setunitscale 0.15
\move(0 0)\lvec(10 0)\lvec(10 10)\lvec(0 10)\lvec(0 0)
\move(10 0)\lvec(15 5)\lvec(15 15)\lvec(5 15)\lvec(0 10)
\move(10 10)\lvec(15 15)
\htext(5 5){$j$}
\end{texdraw}%
}
\quad ($j=2,\cdots,n$) : unit width, unit height, unit thickness,\\

\vskip 3mm

\quad $Y_{\La_0} =
\raisebox{-0.3\height}{\begin{texdraw}
\drawdim em
\setunitscale 0.15
\move(-42 0)\lvec(0 0)\lvec(2.5 2.5)\lvec(2.5 12.5)\lvec(-39.5 12.5)
\move(-42 10)\lvec(0 10)
\move(0 0)\lvec(0 10)\lvec(2.5 12.5)
\move(-10 0)\lvec(-10 10)\lvec(-7.5 12.5)
\move(-20 0)\lvec(-20 10)\lvec(-17.5 12.5)
\move(-30 0)\lvec(-30 10)\lvec(-27.5 12.5)
\move(-40 0)\lvec(-40 10)\lvec(-37.5 12.5)
\move(0 0)\lvec(-2.5 -2.5)\lvec(-44.5 -2.5)
\move(-10 0)\rlvec(-2.5 -2.5)
\move(-20 0)\rlvec(-2.5 -2.5)
\move(-30 0)\rlvec(-2.5 -2.5)
\move(-40 0)\rlvec(-2.5 -2.5)
\htext(-5 5){$1$}
\htext(-15 5){$0$}
\htext(-25 5){$1$}
\htext(-35 5){$0$}
\end{texdraw}}$\,, \\

\vskip 2mm

\quad $Y_{\La_1} =
\raisebox{-0.3\height}{\begin{texdraw}
\drawdim em
\setunitscale 0.15
\move(-42 0)\lvec(0 0)\lvec(2.5 2.5)\lvec(2.5 12.5)\lvec(-39.5 12.5)
\move(-42 10)\lvec(0 10)
\move(0 0)\lvec(0 10)\lvec(2.5 12.5)
\move(-10 0)\lvec(-10 10)\lvec(-7.5 12.5)
\move(-20 0)\lvec(-20 10)\lvec(-17.5 12.5)
\move(-30 0)\lvec(-30 10)\lvec(-27.5 12.5)
\move(-40 0)\lvec(-40 10)\lvec(-37.5 12.5)
\move(0 0)\lvec(-2.5 -2.5)\lvec(-44.5 -2.5)
\move(-10 0)\rlvec(-2.5 -2.5)
\move(-20 0)\rlvec(-2.5 -2.5)
\move(-30 0)\rlvec(-2.5 -2.5)
\move(-40 0)\rlvec(-2.5 -2.5)
\htext(-5 5){$0$}
\htext(-15 5){$1$}
\htext(-25 5){$0$}
\htext(-35 5){$1$}
\end{texdraw}}\,.$

\vskip 7mm

(c) $D_n^{(1)}$ ($n\geq4$)

\vskip 5mm

\hskip 5mm 
\raisebox{-0.33\height}[0.69\height][0.35\height]{%
\begin{texdraw}
\drawdim em
\setunitscale 0.15
\move(0 0)\lvec(10 0)\lvec(10 10)\lvec(0 10)\lvec(0 0)
\move(10 0)\lvec(12.5 2.5)\lvec(12.5 12.5)\lvec(2.5 12.5)\lvec(0 10)
\move(10 10)\lvec(12.5 12.5)
\htext(5 5){$0$}
\end{texdraw}%
}\,,
\raisebox{-0.33\height}[0.69\height][0.35\height]{%
\begin{texdraw}
\drawdim em
\setunitscale 0.15
\move(0 0)\lvec(10 0)\lvec(10 10)\lvec(0 10)\lvec(0 0)
\move(10 0)\lvec(12.5 2.5)\lvec(12.5 12.5)\lvec(2.5 12.5)\lvec(0 10)
\move(10 10)\lvec(12.5 12.5)
\htext(5 5){$1$}
\end{texdraw}%
}\,, 
\raisebox{-0.33\height}[0.69\height][0.35\height]{%
\begin{texdraw}
\drawdim em
\setunitscale 0.15
\move(0 0)\lvec(10 0)\lvec(10 10)\lvec(0 10)\lvec(0 0)
\move(10 0)\lvec(12.5 2.5)\lvec(12.5 12.5)\lvec(2.5 12.5)\lvec(0 10)
\move(10 10)\lvec(12.5 12.5)
\htext(5 5){$_{n-1}$}
\end{texdraw}%
}\,,
\raisebox{-0.33\height}[0.69\height][0.35\height]{%
\begin{texdraw}
\drawdim em
\setunitscale 0.15
\move(0 0)\lvec(10 0)\lvec(10 10)\lvec(0 10)\lvec(0 0)
\move(10 0)\lvec(12.5 2.5)\lvec(12.5 12.5)\lvec(2.5 12.5)\lvec(0 10)
\move(10 10)\lvec(12.5 12.5)
\htext(5 5){$n$}
\end{texdraw}%
} \ :
unit width, unit height, half-unit thickness,\\

\vskip 3mm

\hskip 5mm 
\raisebox{-0.33\height}[0.69\height][0.35\height]{%
\begin{texdraw}
\drawdim em
\setunitscale 0.15
\move(0 0)\lvec(10 0)\lvec(10 10)\lvec(0 10)\lvec(0 0)
\move(10 0)\lvec(15 5)\lvec(15 15)\lvec(5 15)\lvec(0 10)
\move(10 10)\lvec(15 15)
\htext(5 5){$j$}
\end{texdraw}%
}
\quad ($j=2,\cdots,n-2$) : unit width, unit height, unit thickness,\\

\vskip 3mm

\quad $Y_{\La_0} =
\raisebox{-0.3\height}{\begin{texdraw}
\drawdim em
\setunitscale 0.15
\move(-42 0)\lvec(0 0)\lvec(2.5 2.5)\lvec(2.5 12.5)\lvec(-39.5 12.5)
\move(-42 10)\lvec(0 10)
\move(0 0)\lvec(0 10)\lvec(2.5 12.5)
\move(-10 0)\lvec(-10 10)\lvec(-7.5 12.5)
\move(-20 0)\lvec(-20 10)\lvec(-17.5 12.5)
\move(-30 0)\lvec(-30 10)\lvec(-27.5 12.5)
\move(-40 0)\lvec(-40 10)\lvec(-37.5 12.5)
\move(0 0)\lvec(-2.5 -2.5)\lvec(-44.5 -2.5)
\move(-10 0)\rlvec(-2.5 -2.5)
\move(-20 0)\rlvec(-2.5 -2.5)
\move(-30 0)\rlvec(-2.5 -2.5)
\move(-40 0)\rlvec(-2.5 -2.5)
\htext(-5 5){$1$}
\htext(-15 5){$0$}
\htext(-25 5){$1$}
\htext(-35 5){$0$}
\end{texdraw}}$\,,\\

\vskip 2mm

\quad $Y_{\La_1} =
\raisebox{-0.3\height}{\begin{texdraw}
\drawdim em
\setunitscale 0.15
\move(-42 0)\lvec(0 0)\lvec(2.5 2.5)\lvec(2.5 12.5)\lvec(-39.5 12.5)
\move(-42 10)\lvec(0 10)
\move(0 0)\lvec(0 10)\lvec(2.5 12.5)
\move(-10 0)\lvec(-10 10)\lvec(-7.5 12.5)
\move(-20 0)\lvec(-20 10)\lvec(-17.5 12.5)
\move(-30 0)\lvec(-30 10)\lvec(-27.5 12.5)
\move(-40 0)\lvec(-40 10)\lvec(-37.5 12.5)
\move(0 0)\lvec(-2.5 -2.5)\lvec(-44.5 -2.5)
\move(-10 0)\rlvec(-2.5 -2.5)
\move(-20 0)\rlvec(-2.5 -2.5)
\move(-30 0)\rlvec(-2.5 -2.5)
\move(-40 0)\rlvec(-2.5 -2.5)
\htext(-5 5){$0$}
\htext(-15 5){$1$}
\htext(-25 5){$0$}
\htext(-35 5){$1$}
\end{texdraw}}$\,,\\

\vskip 2mm

\quad $Y_{\La_{n-1}} =
\raisebox{-0.3\height}{\begin{texdraw}
\drawdim em
\setunitscale 0.15
\move(-42 0)\lvec(0 0)\lvec(2.5 2.5)\lvec(2.5 12.5)\lvec(-39.5 12.5)
\move(-42 10)\lvec(0 10)
\move(0 0)\lvec(0 10)\lvec(2.5 12.5)
\move(-10 0)\lvec(-10 10)\lvec(-7.5 12.5)
\move(-20 0)\lvec(-20 10)\lvec(-17.5 12.5)
\move(-30 0)\lvec(-30 10)\lvec(-27.5 12.5)
\move(-40 0)\lvec(-40 10)\lvec(-37.5 12.5)
\move(0 0)\lvec(-2.5 -2.5)\lvec(-44.5 -2.5)
\move(-10 0)\rlvec(-2.5 -2.5)
\move(-20 0)\rlvec(-2.5 -2.5)
\move(-30 0)\rlvec(-2.5 -2.5)
\move(-40 0)\rlvec(-2.5 -2.5)
\htext(-5 5){$n$}
\htext(-15 5){$_{n-1}$}
\htext(-25 5){$n$}
\htext(-35 5){$_{n-1}$}
\end{texdraw}}$\,,\\

\vskip 2mm

\quad $Y_{\La_{n}} =
\raisebox{-0.3\height}{\begin{texdraw}
\drawdim em
\setunitscale 0.15
\move(-42 0)\lvec(0 0)\lvec(2.5 2.5)\lvec(2.5 12.5)\lvec(-39.5 12.5)
\move(-42 10)\lvec(0 10)
\move(0 0)\lvec(0 10)\lvec(2.5 12.5)
\move(-10 0)\lvec(-10 10)\lvec(-7.5 12.5)
\move(-20 0)\lvec(-20 10)\lvec(-17.5 12.5)
\move(-30 0)\lvec(-30 10)\lvec(-27.5 12.5)
\move(-40 0)\lvec(-40 10)\lvec(-37.5 12.5)
\move(0 0)\lvec(-2.5 -2.5)\lvec(-44.5 -2.5)
\move(-10 0)\rlvec(-2.5 -2.5)
\move(-20 0)\rlvec(-2.5 -2.5)
\move(-30 0)\rlvec(-2.5 -2.5)
\move(-40 0)\rlvec(-2.5 -2.5)
\htext(-5 5){$_{n-1}$}
\htext(-15 5){$n$}
\htext(-25 5){$_{n-1}$}
\htext(-35 5){$n$}
\end{texdraw}}$\,.\\

\vskip 5mm

\newpage
(d) $A_{2n}^{(2)}$ ($n\geq2$)

\vskip 3mm
\hskip 5mm 
\raisebox{-0.33\height}[0.69\height][0.35\height]{%
\begin{texdraw}
\drawdim em
\setunitscale 0.15
\move(0 0)\lvec(10 0)\lvec(10 5)\lvec(0 5)\lvec(0 0)
\move(10 0)\lvec(15 5)\lvec(15 10)\lvec(5 10)\lvec(0 5)
\move(10 5)\lvec(15 10)
\htext(5 2.5){$_{0}$}
\end{texdraw}%
} \ :
unit width, half-unit height, unit thickness,\\

\vskip 3mm
\hskip 5mm 
\raisebox{-0.33\height}[0.69\height][0.35\height]{%
\begin{texdraw}
\drawdim em
\setunitscale 0.15
\move(0 0)\lvec(10 0)\lvec(10 10)\lvec(0 10)\lvec(0 0)
\move(10 0)\lvec(15 5)\lvec(15 15)\lvec(5 15)\lvec(0 10)
\move(10 10)\lvec(15 15)
\htext(5 5){$j$}
\end{texdraw}%
}
\quad ($j=1,\cdots,n$) : unit width, unit height, unit thickness,\\

\vskip 3mm

\quad $Y_{\La_0} = \ 
\raisebox{-0.33\height}[0.69\height][0.35\height]{%
\begin{texdraw}
\drawdim em
\setunitscale 0.15
\move(-3 0)\lvec(40 0)
\move(-3 5)\lvec(40 5)
\move(-0.5 10)\lvec(45 10)
\move(40 0)\lvec(40 5)\lvec(45 10)\lvec(45 5)\lvec(40 0)
\move(0 0)\lvec(0 5)\lvec(5 10)
\move(10 0)\lvec(10 5)\lvec(15 10)
\move(20 0)\lvec(20 5)\lvec(25 10)
\move(30 0)\lvec(30 5)\lvec(35 10)
\htext(5 2.5){$_{0}$}
\htext(15 2.5){$_{0}$}
\htext(25 2.5){$_{0}$}
\htext(35 2.5){$_{0}$}
\end{texdraw}%
}\,. $

\vskip 7mm

(e) $D_{n+1}^{(2)}$ ($n\geq2$)

\vskip 3mm
\hskip 5mm 
\raisebox{-0.33\height}[0.69\height][0.35\height]{%
\begin{texdraw}
\drawdim em
\setunitscale 0.15
\move(0 0)\lvec(10 0)\lvec(10 5)\lvec(0 5)\lvec(0 0)
\move(10 0)\lvec(15 5)\lvec(15 10)\lvec(5 10)\lvec(0 5)
\move(10 5)\lvec(15 10)
\htext(5 2.5){$_{0}$}
\end{texdraw}%
}\,, \ 
\raisebox{-0.33\height}[0.69\height][0.35\height]{%
\begin{texdraw}
\drawdim em
\setunitscale 0.15
\move(0 0)\lvec(10 0)\lvec(10 5)\lvec(0 5)\lvec(0 0)
\move(10 0)\lvec(15 5)\lvec(15 10)\lvec(5 10)\lvec(0 5)
\move(10 5)\lvec(15 10)
\htext(5 2.5){$_{n}$}
\end{texdraw}%
} \ :
unit width, half-unit height, unit thickness,\\

\vskip 3mm
\hskip 5mm 
\raisebox{-0.33\height}[0.69\height][0.35\height]{%
\begin{texdraw}
\drawdim em
\setunitscale 0.15
\move(0 0)\lvec(10 0)\lvec(10 10)\lvec(0 10)\lvec(0 0)
\move(10 0)\lvec(15 5)\lvec(15 15)\lvec(5 15)\lvec(0 10)
\move(10 10)\lvec(15 15)
\htext(5 5){$j$}
\end{texdraw}%
}
\quad ($j=1,\cdots,n-1$) : unit width, unit height, unit thickness,\\

\vskip 3mm

\quad $Y_{\La_0} = \ 
\raisebox{-0.33\height}[0.69\height][0.35\height]{%
\begin{texdraw}
\drawdim em
\setunitscale 0.15
\move(-3 0)\lvec(40 0)
\move(-3 5)\lvec(40 5)
\move(-0.5 10)\lvec(45 10)
\move(40 0)\lvec(40 5)\lvec(45 10)\lvec(45 5)\lvec(40 0)
\move(0 0)\lvec(0 5)\lvec(5 10)
\move(10 0)\lvec(10 5)\lvec(15 10)
\move(20 0)\lvec(20 5)\lvec(25 10)
\move(30 0)\lvec(30 5)\lvec(35 10)
\htext(5 2.5){$_{0}$}
\htext(15 2.5){$_{0}$}
\htext(25 2.5){$_{0}$}
\htext(35 2.5){$_{0}$}
\end{texdraw}%
}\,, $

\vskip 3mm

\quad $Y_{\La_n} = \ 
\raisebox{-0.33\height}[0.69\height][0.35\height]{%
\begin{texdraw}
\drawdim em
\setunitscale 0.15
\move(-3 0)\lvec(40 0)
\move(-3 5)\lvec(40 5)
\move(-0.5 10)\lvec(45 10)
\move(40 0)\lvec(40 5)\lvec(45 10)\lvec(45 5)\lvec(40 0)
\move(0 0)\lvec(0 5)\lvec(5 10)
\move(10 0)\lvec(10 5)\lvec(15 10)
\move(20 0)\lvec(20 5)\lvec(25 10)
\move(30 0)\lvec(30 5)\lvec(35 10)
\htext(5 2.5){$_{n}$}
\htext(15 2.5){$_{n}$}
\htext(25 2.5){$_{n}$}
\htext(35 2.5){$_{n}$}
\end{texdraw}%
}\,. $

\vskip 7mm

(f) $B_n^{(1)}$ ($n\geq3$)

\vskip 3mm
\hskip 5mm 
\raisebox{-0.33\height}[0.69\height][0.35\height]{%
\begin{texdraw}
\drawdim em
\setunitscale 0.15
\move(0 0)\lvec(10 0)\lvec(10 10)\lvec(0 10)\lvec(0 0)
\move(10 0)\lvec(12.5 2.5)\lvec(12.5 12.5)\lvec(2.5 12.5)\lvec(0 10)
\move(10 10)\lvec(12.5 12.5)
\htext(5 5){$0$}
\end{texdraw}%
}\,,
\raisebox{-0.33\height}[0.69\height][0.35\height]{%
\begin{texdraw}
\drawdim em
\setunitscale 0.15
\move(0 0)\lvec(10 0)\lvec(10 10)\lvec(0 10)\lvec(0 0)
\move(10 0)\lvec(12.5 2.5)\lvec(12.5 12.5)\lvec(2.5 12.5)\lvec(0 10)
\move(10 10)\lvec(12.5 12.5)
\htext(5 5){$1$}
\end{texdraw}%
}\ :
unit width, unit height, half-unit thickness,\\

\vskip 3mm
\hskip 5mm 
\raisebox{-0.33\height}[0.69\height][0.35\height]{%
\begin{texdraw}
\drawdim em
\setunitscale 0.15
\move(0 0)\lvec(10 0)\lvec(10 5)\lvec(0 5)\lvec(0 0)
\move(10 0)\lvec(15 5)\lvec(15 10)\lvec(5 10)\lvec(0 5)
\move(10 5)\lvec(15 10)
\htext(5 2.5){$_{n}$}
\end{texdraw}%
} \ :
unit width, half-unit height, unit thickness,\\

\vskip 3mm
\hskip 5mm 
\raisebox{-0.33\height}[0.69\height][0.35\height]{%
\begin{texdraw}
\drawdim em
\setunitscale 0.15
\move(0 0)\lvec(10 0)\lvec(10 10)\lvec(0 10)\lvec(0 0)
\move(10 0)\lvec(15 5)\lvec(15 15)\lvec(5 15)\lvec(0 10)
\move(10 10)\lvec(15 15)
\htext(5 5){$j$}
\end{texdraw}%
}
\quad ($j=2,\cdots,n-1$) : unit width, unit height, unit thickness,\\

\vskip 3mm

\quad $Y_{\La_0} =
\raisebox{-0.3\height}{\begin{texdraw}
\drawdim em
\setunitscale 0.15
\move(-42 0)\lvec(0 0)\lvec(2.5 2.5)\lvec(2.5 12.5)\lvec(-39.5 12.5)
\move(-42 10)\lvec(0 10)
\move(0 0)\lvec(0 10)\lvec(2.5 12.5)
\move(-10 0)\lvec(-10 10)\lvec(-7.5 12.5)
\move(-20 0)\lvec(-20 10)\lvec(-17.5 12.5)
\move(-30 0)\lvec(-30 10)\lvec(-27.5 12.5)
\move(-40 0)\lvec(-40 10)\lvec(-37.5 12.5)
\move(0 0)\lvec(-2.5 -2.5)\lvec(-44.5 -2.5)
\move(-10 0)\rlvec(-2.5 -2.5)
\move(-20 0)\rlvec(-2.5 -2.5)
\move(-30 0)\rlvec(-2.5 -2.5)
\move(-40 0)\rlvec(-2.5 -2.5)
\htext(-5 5){$1$}
\htext(-15 5){$0$}
\htext(-25 5){$1$}
\htext(-35 5){$0$}
\end{texdraw}}$\,, \\

\vskip 2mm

\quad $Y_{\La_1} =
\raisebox{-0.3\height}{\begin{texdraw}
\drawdim em
\setunitscale 0.15
\move(-42 0)\lvec(0 0)\lvec(2.5 2.5)\lvec(2.5 12.5)\lvec(-39.5 12.5)
\move(-42 10)\lvec(0 10)
\move(0 0)\lvec(0 10)\lvec(2.5 12.5)
\move(-10 0)\lvec(-10 10)\lvec(-7.5 12.5)
\move(-20 0)\lvec(-20 10)\lvec(-17.5 12.5)
\move(-30 0)\lvec(-30 10)\lvec(-27.5 12.5)
\move(-40 0)\lvec(-40 10)\lvec(-37.5 12.5)
\move(0 0)\lvec(-2.5 -2.5)\lvec(-44.5 -2.5)
\move(-10 0)\rlvec(-2.5 -2.5)
\move(-20 0)\rlvec(-2.5 -2.5)
\move(-30 0)\rlvec(-2.5 -2.5)
\move(-40 0)\rlvec(-2.5 -2.5)
\htext(-5 5){$0$}
\htext(-15 5){$1$}
\htext(-25 5){$0$}
\htext(-35 5){$1$}
\end{texdraw}}\,,$

\vskip 3mm

\quad $Y_{\La_n} = \ 
\raisebox{-0.33\height}[0.69\height][0.35\height]{%
\begin{texdraw}
\drawdim em
\setunitscale 0.15
\move(-3 0)\lvec(40 0)
\move(-3 5)\lvec(40 5)
\move(-0.5 10)\lvec(45 10)
\move(40 0)\lvec(40 5)\lvec(45 10)\lvec(45 5)\lvec(40 0)
\move(0 0)\lvec(0 5)\lvec(5 10)
\move(10 0)\lvec(10 5)\lvec(15 10)
\move(20 0)\lvec(20 5)\lvec(25 10)
\move(30 0)\lvec(30 5)\lvec(35 10)
\htext(5 2.5){$_{n}$}
\htext(15 2.5){$_{n}$}
\htext(25 2.5){$_{n}$}
\htext(35 2.5){$_{n}$}
\end{texdraw}%
}\,. $

\vskip 7mm 
\noindent
For convenience, we will use the following notations:

\vskip 3mm

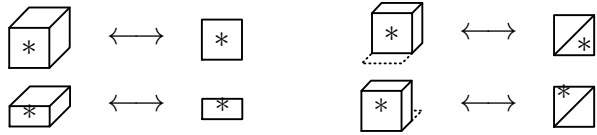
\begin{center}
\begin{tabular}{rcl}
\raisebox{-0.4\height}{
\begin{texdraw}
\drawdim em
\setunitscale 0.13
\move(-10 0)\lvec(0 0)\lvec(0 10)\lvec(-10 10)\lvec(-10 0)
\move(0 0)\lvec(5 5)\lvec(5 15)\lvec(-5 15)\lvec(-10 10)
\move(0 10)\lvec(5 15)
\htext(-5 5){$*$}
\end{texdraw}
}
& $\longleftrightarrow$ &
\raisebox{-0.4\height}{
\begin{texdraw}
\drawdim em
\setunitscale 0.13
\move(-10 0)\lvec(0 0)\lvec(0 10)\lvec(-10 10)\lvec(-10 0)
\htext(-5 5){$*$}
\end{texdraw}
}\\[4mm]
\raisebox{-0.4\height}{
\begin{texdraw}
\drawdim em
\setunitscale 0.13
\move(0 0)\lvec(10 0)\lvec(10 5)\lvec(0 5)\lvec(0 0)
\move(10 0)\lvec(15 5)\lvec(15 10)\lvec(5 10)\lvec(0 5)
\move(10 5)\lvec(15 10)
\htext(5 2.5){*}
\end{texdraw}
}
& $\longleftrightarrow$ &
\raisebox{-0.4\height}{
\begin{texdraw}
\drawdim em
\setunitscale 0.13
\textref h:C v:C
\move(0 0)\lvec(10 0)\lvec(10 5)\lvec(0 5)\lvec(0 0)
\htext(5 2.5){*}
\end{texdraw}
}
\end{tabular}
\qquad
\begin{tabular}{rcl}
\raisebox{-0.4\height}{
\begin{texdraw}
\drawdim em
\setunitscale 0.13
\move(-10 0)\lvec(0 0)\lvec(0 10)\lvec(-10 10)\lvec(-10 0)
\move(0 0)\lvec(2.5 2.5)\lvec(2.5 12.5)\lvec(-7.5 12.5)\lvec(-10 10)
\move(0 10)\lvec(2.5 12.5)
\lpatt(0.3 1)
\move(0 0)\lvec(-2.5 -2.5)\lvec(-12.5 -2.5)\lvec(-10 0)
\htext(-5 5){$*$}
\end{texdraw}
}
& $\longleftrightarrow$ &
\raisebox{-0.4\height}{
\begin{texdraw}
\drawdim em
\setunitscale 0.13
\move(-10 0)\lvec(0 0)\lvec(0 10)\lvec(-10 10)\lvec(-10 0)
\move(0 10)\lvec(-10 0)
\htext(-2.5 2.5){$*$}
\end{texdraw}
}\\[4mm]
\raisebox{-0.4\height}{
\begin{texdraw}
\drawdim em
\setunitscale 0.13
\move(-10 0)\lvec(0 0)\lvec(0 10)\lvec(-10 10)\lvec(-10 0)
\move(0 0)\lvec(2.5 2.5)\lvec(2.5 12.5)\lvec(-7.5 12.5)\lvec(-10 10)
\move(0 10)\lvec(2.5 12.5)
\lpatt(0.3 1)
\move(2.5 2.5)\lvec(5 5)\lvec(2.5 5)
\htext(-5 5){*}
\end{texdraw}
}
& $\longleftrightarrow$ &
\raisebox{-0.4\height}{
\begin{texdraw}
\drawdim em
\setunitscale 0.13
\move(-10 0)\lvec(0 0)\lvec(0 10)\lvec(-10 10)\lvec(-10 0)
\move(0 10)\lvec(-10 0)
\htext(-7.5 7.5){*}
\end{texdraw}
}
\end{tabular}
\end{center}

\vskip 3mm
For example, we have

\vskip 3mm

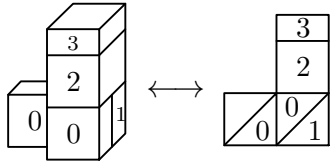
\begin{center}
\raisebox{-0.4\height}{
\begin{texdraw}
\small
\drawdim em
\setunitscale 0.18
\move(0 0)\lvec(0 25)\lvec(-10 25)\lvec(-10 0)\lvec(0 0)
\move(-10 10)\lvec(0 10)\lvec(5 15)
\move(-10 20)\lvec(0 20)\lvec(5 25)
\move(0 0)\lvec(5 5)\lvec(5 30)
\move(0 25)\lvec(5 30)\lvec(-5 30)\lvec(-10 25)
\move(2.5 2.5)\lvec(2.5 12.5)
\move(-10 2.5)\lvec(-17.5 2.5)\lvec(-17.5 12.5)\lvec(-10 12.5)
\move(-17.5 12.5)\lvec(-15 15)\lvec(-10 15)
\htext(-5 5){$0$}
\htext(-5 15){$2$}
\htext(-12.5 7.5){$0$}
\htext(4.2 8.7){$_1$}
\htext(-5 22.5){$_{3}$}
\end{texdraw}
}
$\longleftrightarrow$
\raisebox{-0.37\height}{
\begin{texdraw}
\small
\drawdim em
\setunitscale 0.18
\move(0 0)\lvec(0 10)\lvec(-20 10)\lvec(-20 0)\lvec(0 0)
\move(-20 0)\lvec(-10 10)\lvec(-10 0)\lvec(0 10)
\move(-10 10)\lvec(-10 25)\lvec(0 25)\lvec(0 10)
\move(-10 20)\lvec(0 20)
\htext(-12.8 2.8){$0$}
\htext(-7.2 7.5){$0$}
\htext(-2.8 2.8){$1$}
\htext(-5 15){$2$}
\htext(-5 22.5){$3$}
\end{texdraw}
}
\end{center}

\vskip 3mm

The rules for building the walls are given as follows:
\begin{enumerate}
\item The walls must be built on top of the ground-state walls.
\item The colored blocks should be stacked in columns.
      No block can be placed on top of a column of half-unit thickness.
\item Except for the right-most column, there should be no free space
      to the right of any block.
\item The colored blocks should be stacked in the patterns given below.

\vskip 5mm 

\begin{enumerate}
\savebox{\tmpfiga}{\begin{texdraw}
\fontsize{8}{8}\selectfont
\drawdim em
\setunitscale 1.9
\move(0 0)\rlvec(-7.7 0) \move(0 1)\rlvec(-7.7 0) \move(0 2)\rlvec(-7.7 0)
\move(0 3.5)\rlvec(-7.7 0) \move(0 4.5)\rlvec(-7.7 0)
\move(0 5.5)\rlvec(-7.7 0) \move(0 6.5)\rlvec(-7.7 0)
\move(0 0)\rlvec(0 6.7) \move(-1 0)\rlvec(0 6.7) \move(-2 0)\rlvec(0 6.7)
\move(-3.5 0)\rlvec(0 6.7) \move(-4.5 0)\rlvec(0 6.7)
\move(-5.5 0)\rlvec(0 6.7) \move(-6.5 0)\rlvec(0 6.7)
\move(-7.5 0)\rlvec(0 6.7) \move(-0.5 0.5)
\bsegment
\htext(0 0){$i$} \htext(0 1){$i\!\!+\!\!1$} \vtext(0 2.25){$\cdots$}
\htext(0 3.5){$n$} \htext(0 4.5){$0$} \htext(0 5.5){$1$}
\htext(-1 0){$i\!\!-\!\!1$} \htext(-1 1){$i$} \vtext(-1 2.25){$\cdots$}
\htext(-1 3.5){$n\!\!-\!\!1$} \htext(-1 4.5){$n$} \htext(-1 5.5){$0$}
\htext(-2.25 0){$\cdots$} \htext(-2.25 1){$\cdots$} \htext(-3.5 0){$2$}
\htext(-3.5 1){$3$} \htext(-4.5 0){$1$} \htext(-4.5 1){$2$}
\htext(-5.5 0){$0$} \htext(-5.5 1){$1$} \htext(-6.5 0){$n$}
\htext(-6.5 1){$0$}
\esegment
\end{texdraw}}%
\savebox{\tmpfigb}{\begin{texdraw}
\fontsize{7}{7}\selectfont
\drawdim em
\setunitscale 1.9
\nc{\dtri}{
\bsegment
\move(-1 0)\lvec(0 1)\lvec(0 0)\lvec(-1 0)\ifill f:0.7
\esegment
}
\move(0 0)\dtri \move(-1 0)\dtri \move(-2 0)\dtri \move(-3 0)\dtri
\move(0 0)\rlvec(-4.3 0) \move(0 1)\rlvec(-4.3 0) \move(0 2)\rlvec(-4.3 0)
\move(0 3.5)\rlvec(-4.3 0) \move(0 4.5)\rlvec(-4.3 0) \move(0 6)\rlvec(-4.3 0)
\move(0 7)\rlvec(-4.3 0) \move(0 8)\rlvec(-4.3 0) \move(0 9)\rlvec(-4.3 0)
\move(0 0)\rlvec(0 9.3) \move(-1 0)\rlvec(0 9.3) \move(-2 0)\rlvec(0 9.3)
\move(-3 0)\rlvec(0 9.3) \move(-4 0)\rlvec(0 9.3)
\move(-1 0)\rlvec(1 1) \move(-2 0)\rlvec(1 1) \move(-3 0)\rlvec(1 1)
\move(-4 0)\rlvec(1 1) \move(-1 7)\rlvec(1 1) \move(-2 7)\rlvec(1 1)
\move(-3 7)\rlvec(1 1) \move(-4 7)\rlvec(1 1)
\vtext(-0.5 2.75){$\cdots$} \vtext(-0.5 5.25){$\cdots$}
\vtext(-1.5 2.75){$\cdots$} \vtext(-1.5 5.25){$\cdots$}
\vtext(-2.5 2.75){$\cdots$} \vtext(-2.5 5.25){$\cdots$}
\vtext(-3.5 2.75){$\cdots$} \vtext(-3.5 5.25){$\cdots$}
\htext(-0.25 7.27){$1$} \htext(-0.75 7.75){$0$} \htext(-1.25 7.27){$0$}
\htext(-1.75 7.75){$1$} \htext(-2.25 7.27){$1$} \htext(-2.75 7.75){$0$}
\htext(-3.25 7.27){$0$} \htext(-3.75 7.75){$1$}
\htext(-0.25 0.27){$1$} \htext(-0.75 0.75){$0$} \htext(-1.25 0.27){$0$}
\htext(-1.75 0.75){$1$} \htext(-2.25 0.27){$1$} \htext(-2.75 0.75){$0$}
\htext(-3.25 0.27){$0$} \htext(-3.75 0.75){$1$}
\htext(-0.5 1.5){$2$} \htext(-1.5 1.5){$2$} \htext(-2.5 1.5){$2$}
\htext(-3.5 1.5){$2$}
\htext(-0.5 6.5){$2$} \htext(-1.5 6.5){$2$} \htext(-2.5 6.5){$2$}
\htext(-3.5 6.5){$2$}
\htext(-0.5 8.5){$2$} \htext(-1.5 8.5){$2$} \htext(-2.5 8.5){$2$}
\htext(-3.5 8.5){$2$}
\htext(-0.5 4){$n$} \htext(-1.5 4){$n$} \htext(-2.5 4){$n$}
\htext(-3.5 4){$n$}
\end{texdraw}}%
\savebox{\tmpfigc}{\begin{texdraw}
\fontsize{7}{7}\selectfont
\drawdim em
\setunitscale 1.9
\nc{\dtri}{
\bsegment
\move(-1 0)\lvec(0 1)\lvec(0 0)\lvec(-1 0)\ifill f:0.7
\esegment
}
\move(0 0)\dtri \move(-1 0)\dtri \move(-2 0)\dtri \move(-3 0)\dtri
\move(0 0)\rlvec(-4.3 0) \move(0 1)\rlvec(-4.3 0) \move(0 2)\rlvec(-4.3 0)
\move(0 3.5)\rlvec(-4.3 0) \move(0 4.5)\rlvec(-4.3 0) \move(0 6)\rlvec(-4.3 0)
\move(0 7)\rlvec(-4.3 0) \move(0 8)\rlvec(-4.3 0) \move(0 9)\rlvec(-4.3 0)
\move(0 0)\rlvec(0 9.3) \move(-1 0)\rlvec(0 9.3) \move(-2 0)\rlvec(0 9.3)
\move(-3 0)\rlvec(0 9.3) \move(-4 0)\rlvec(0 9.3)
\move(-1 0)\rlvec(1 1) \move(-2 0)\rlvec(1 1) \move(-3 0)\rlvec(1 1)
\move(-4 0)\rlvec(1 1) \move(-1 7)\rlvec(1 1) \move(-2 7)\rlvec(1 1)
\move(-3 7)\rlvec(1 1) \move(-4 7)\rlvec(1 1)
\vtext(-0.5 2.75){$\cdots$} \vtext(-0.5 5.25){$\cdots$}
\vtext(-1.5 2.75){$\cdots$} \vtext(-1.5 5.25){$\cdots$}
\vtext(-2.5 2.75){$\cdots$} \vtext(-2.5 5.25){$\cdots$}
\vtext(-3.5 2.75){$\cdots$} \vtext(-3.5 5.25){$\cdots$}
\htext(-0.25 7.27){$0$} \htext(-0.75 7.75){$1$} \htext(-1.25 7.27){$1$}
\htext(-1.75 7.75){$0$} \htext(-2.25 7.27){$0$} \htext(-2.75 7.75){$1$}
\htext(-3.25 7.27){$1$} \htext(-3.75 7.75){$0$}
\htext(-0.25 0.27){$0$} \htext(-0.75 0.75){$1$} \htext(-1.25 0.27){$1$}
\htext(-1.75 0.75){$0$} \htext(-2.25 0.27){$0$} \htext(-2.75 0.75){$1$}
\htext(-3.25 0.27){$1$} \htext(-3.75 0.75){$0$}
\htext(-0.5 1.5){$2$} \htext(-1.5 1.5){$2$} \htext(-2.5 1.5){$2$}
\htext(-3.5 1.5){$2$}
\htext(-0.5 6.5){$2$} \htext(-1.5 6.5){$2$} \htext(-2.5 6.5){$2$}
\htext(-3.5 6.5){$2$}
\htext(-0.5 8.5){$2$} \htext(-1.5 8.5){$2$} \htext(-2.5 8.5){$2$}
\htext(-3.5 8.5){$2$}
\htext(-0.5 4){$n$} \htext(-1.5 4){$n$} \htext(-2.5 4){$n$}
\htext(-3.5 4){$n$}
\end{texdraw}}%
\savebox{\tmpfigd}{\begin{texdraw}
\fontsize{10}{10}\selectfont
\drawdim em
\setunitscale 1.9
\nc{\dtri}{
\bsegment
\move(-1 0)\lvec(0 1)\lvec(0 0)\lvec(-1 0)\ifill f:0.7
\esegment
}
\move(0 0)\dtri \move(-1 0)\dtri \move(-2 0)\dtri \move(-3 0)\dtri
\move(0 0)\rlvec(-4.3 0) \move(0 1)\rlvec(-4.3 0) \move(0 2)\rlvec(-4.3 0)
\move(0 3.5)\rlvec(-4.3 0) \move(0 4.5)\rlvec(-4.3 0)
\move(0 5.5)\rlvec(-4.3 0) \move(0 6.5)\rlvec(-4.3 0)
\move(0 8)\rlvec(-4.3 0) \move(0 9)\rlvec(-4.3 0)
\move(0 10)\rlvec(-4.3 0) \move(0 11)\rlvec(-4.3 0)
\move(0 0)\rlvec(0 11.3) \move(-1 0)\rlvec(0 11.3) \move(-2 0)\rlvec(0 11.3)
\move(-3 0)\rlvec(0 11.3) \move(-4 0)\rlvec(0 11.3)
\move(-1 0)\rlvec(1 1) \move(-2 0)\rlvec(1 1) \move(-3 0)\rlvec(1 1)
\move(-4 0)\rlvec(1 1) \move(-1 9)\rlvec(1 1) \move(-2 9)\rlvec(1 1)
\move(-3 9)\rlvec(1 1) \move(-4 9)\rlvec(1 1)
\htext(-0.3 0.25){$1$} \htext(-0.75 0.75){$0$} \htext(-0.5 1.5){$2$}
\vtext(-0.5 2.75){$\cdots$} \htext(-0.5 4){$n\!\!-\!\!2$}
\htext(-0.5 6){$n\!\!-\!\!2$} \htext(-0.5 8.5){$2$}
\htext(-0.3 9.25){$1$} \htext(-0.75 9.75){$0$} \htext(-0.5 10.5){$2$}
\htext(-2.3 0.25){$1$} \htext(-2.75 0.75){$0$} \htext(-2.5 1.5){$2$}
\vtext(-2.5 2.75){$\cdots$} \htext(-2.5 4){$n\!\!-\!\!2$}
\htext(-2.5 6){$n\!\!-\!\!2$} \htext(-2.5 8.5){$2$}
\htext(-2.3 9.25){$1$} \htext(-2.75 9.75){$0$} \htext(-2.5 10.5){$2$}
\htext(-1.3 0.25){$0$} \htext(-1.75 0.75){$1$} \htext(-1.5 1.5){$2$}
\vtext(-1.5 2.75){$\cdots$} \htext(-1.5 4){$n\!\!-\!\!2$}
\htext(-1.5 6){$n\!\!-\!\!2$} \htext(-1.5 8.5){$2$}
\htext(-1.3 9.25){$0$} \htext(-1.75 9.75){$1$} \htext(-1.5 10.5){$2$}
\htext(-3.3 0.25){$0$} \htext(-3.75 0.75){$1$} \htext(-3.5 1.5){$2$}
\vtext(-3.5 2.75){$\cdots$} \htext(-3.5 4){$n\!\!-\!\!2$}
\htext(-3.5 6){$n\!\!-\!\!2$} \htext(-3.5 8.5){$2$}
\htext(-3.3 9.25){$0$} \htext(-3.75 9.75){$1$} \htext(-3.5 10.5){$2$}
\htext(-0.3 4.75){$n$} \htext(-2.3 4.75){$n$} \htext(-1.75 5.25){$n$}
\htext(-3.75 5.25){$n$}
\move(-1 4.5)\rlvec(0.5 0.5)\rmove(0.45 0.45)\rlvec(0.05 0.05)
\move(-2 4.5)\rlvec(0.05 0.05)\rmove(0.45 0.45)\rlvec(0.5 0.5)
\move(-3 4.5)\rlvec(0.5 0.5)\rmove(0.45 0.45)\rlvec(0.05 0.05)
\move(-4 4.5)\rlvec(0.05 0.05)\rmove(0.45 0.45)\rlvec(0.5 0.5)
\htext(-1.5 4.75){$n\!\!-\!\!1$} \htext(-3.5 4.75){$n\!\!-\!\!1$}
\htext(-0.5 5.25){$n\!\!-\!\!1$} \htext(-2.5 5.25){$n\!\!-\!\!1$}
\vtext(-0.5 7.25){$\cdots$} \vtext(-1.5 7.25){$\cdots$}
\vtext(-2.5 7.25){$\cdots$} \vtext(-3.5 7.25){$\cdots$}
\end{texdraw}}%
\savebox{\tmpfige}{\begin{texdraw}
\fontsize{10}{10}\selectfont
\drawdim em
\setunitscale 1.9
\nc{\dtri}{
\bsegment
\move(-1 0)\lvec(0 1)\lvec(0 0)\lvec(-1 0)\ifill f:0.7
\esegment
}
\move(0 0)\dtri \move(-1 0)\dtri \move(-2 0)\dtri \move(-3 0)\dtri
\move(0 0)\rlvec(-4.3 0) \move(0 1)\rlvec(-4.3 0) \move(0 2)\rlvec(-4.3 0)
\move(0 3.5)\rlvec(-4.3 0) \move(0 4.5)\rlvec(-4.3 0)
\move(0 5.5)\rlvec(-4.3 0) \move(0 6.5)\rlvec(-4.3 0)
\move(0 8)\rlvec(-4.3 0) \move(0 9)\rlvec(-4.3 0) \move(0 10)\rlvec(-4.3 0)
\move(0 11)\rlvec(-4.3 0) \move(0 0)\rlvec(0 11.3) \move(-1 0)\rlvec(0 11.3)
\move(-2 0)\rlvec(0 11.3) \move(-3 0)\rlvec(0 11.3) \move(-4 0)\rlvec(0 11.3)
\move(-1 0)\rlvec(1 1) \move(-2 0)\rlvec(1 1)
\move(-3 0)\rlvec(1 1) \move(-4 0)\rlvec(1 1)
\move(-1 9)\rlvec(1 1) \move(-2 9)\rlvec(1 1)
\move(-3 9)\rlvec(1 1) \move(-4 9)\rlvec(1 1)
\htext(-0.3 0.25){$0$} \htext(-0.75 0.75){$1$} \htext(-0.5 1.5){$2$}
\vtext(-0.5 2.75){$\cdots$} \htext(-0.5 4){$n\!\!-\!\!2$}
\htext(-0.5 6){$n\!\!-\!\!2$} \htext(-0.5 8.5){$2$}
\htext(-0.3 9.25){$0$} \htext(-0.75 9.75){$1$} \htext(-0.5 10.5){$2$}
\htext(-2.3 0.25){$0$} \htext(-2.75 0.75){$1$} \htext(-2.5 1.5){$2$}
\vtext(-2.5 2.75){$\cdots$} \htext(-2.5 4){$n\!\!-\!\!2$}
\htext(-2.5 6){$n\!\!-\!\!2$} \htext(-2.5 8.5){$2$}
\htext(-2.3 9.25){$0$} \htext(-2.75 9.75){$1$} \htext(-2.5 10.5){$2$}
\htext(-1.3 0.25){$1$} \htext(-1.75 0.75){$0$} \htext(-1.5 1.5){$2$}
\vtext(-1.5 2.75){$\cdots$} \htext(-1.5 4){$n\!\!-\!\!2$}
\htext(-1.5 6){$n\!\!-\!\!2$} \htext(-1.5 8.5){$2$}
\htext(-1.3 9.25){$1$} \htext(-1.75 9.75){$0$} \htext(-1.5 10.5){$2$}
\htext(-3.3 0.25){$1$} \htext(-3.75 0.75){$0$} \htext(-3.5 1.5){$2$}
\vtext(-3.5 2.75){$\cdots$} \htext(-3.5 4){$n\!\!-\!\!2$}
\htext(-3.5 6){$n\!\!-\!\!2$} \htext(-3.5 8.5){$2$}
\htext(-3.3 9.25){$1$} \htext(-3.75 9.75){$0$} \htext(-3.5 10.5){$2$}
\htext(-0.3 4.75){$n$} \htext(-2.3 4.75){$n$} \htext(-1.75 5.25){$n$}
\htext(-3.75 5.25){$n$}
\move(-1 4.5)\rlvec(0.5 0.5)\rmove(0.45 0.45)\rlvec(0.05 0.05)
\move(-2 4.5)\rlvec(0.05 0.05)\rmove(0.45 0.45)\rlvec(0.5 0.5)
\move(-3 4.5)\rlvec(0.5 0.5)\rmove(0.45 0.45)\rlvec(0.05 0.05)
\move(-4 4.5)\rlvec(0.05 0.05)\rmove(0.45 0.45)\rlvec(0.5 0.5)
\htext(-1.5 4.75){$n\!\!-\!\!1$} \htext(-3.5 4.75){$n\!\!-\!\!1$}
\htext(-0.5 5.25){$n\!\!-\!\!1$} \htext(-2.5 5.25){$n\!\!-\!\!1$}
\vtext(-0.5 7.25){$\cdots$} \vtext(-1.5 7.25){$\cdots$}
\vtext(-2.5 7.25){$\cdots$} \vtext(-3.5 7.25){$\cdots$}
\end{texdraw}}%
\savebox{\tmpfigf}{\begin{texdraw}
\fontsize{10}{10}\selectfont
\drawdim em
\setunitscale 1.9
\nc{\dtri}{
\bsegment
\move(-1 0)\lvec(0 1)\lvec(0 0)\lvec(-1 0)\ifill f:0.7
\esegment
}
\move(0 0)\dtri \move(-1 0)\dtri \move(-2 0)\dtri \move(-3 0)\dtri
\move(0 0)\rlvec(-4.3 0) \move(0 1)\rlvec(-4.3 0) \move(0 2)\rlvec(-4.3 0)
\move(0 3.5)\rlvec(-4.3 0) \move(0 4.5)\rlvec(-4.3 0) \move(0 5.5)\rlvec(-4.3 0)
\move(0 6.5)\rlvec(-4.3 0) \move(0 8)\rlvec(-4.3 0) \move(0 9)\rlvec(-4.3 0)
\move(0 10)\rlvec(-4.3 0) \move(0 11)\rlvec(-4.3 0)
\move(0 0)\rlvec(0 11.3) \move(-1 0)\rlvec(0 11.3) \move(-2 0)\rlvec(0 11.3)
\move(-3 0)\rlvec(0 11.3) \move(-4 0)\rlvec(0 11.3)
\move(-1 4.5)\rlvec(1 1) \move(-2 4.5)\rlvec(1 1) \move(-3 4.5)\rlvec(1 1)
\move(-4 4.5)\rlvec(1 1)
\htext(-0.3 0.25){$n$}
\htext(-0.5 1.5){$n\!\!-\!\!2$} \vtext(-0.5 2.75){$\cdots$} \htext(-0.5 4){$2$}
\htext(-0.3 4.75){$1$} \htext(-0.75 5.25){$0$} \htext(-0.5 6){$2$}
\vtext(-0.5 7.25){$\cdots$} \htext(-0.5 8.5){$n\!\!-\!\!2$}
\htext(-0.3 9.25){$n$}
\htext(-0.5 10.5){$n\!\!-\!\!2$}
\htext(-2.3 0.25){$n$}
\htext(-2.5 1.5){$n\!\!-\!\!2$} \vtext(-2.5 2.75){$\cdots$} \htext(-2.5 4){$2$}
\htext(-2.3 4.75){$1$} \htext(-2.75 5.25){$0$} \htext(-2.5 6){$2$}
\vtext(-2.5 7.25){$\cdots$} \htext(-2.5 8.5){$n\!\!-\!\!2$}
\htext(-2.3 9.25){$n$}
\htext(-2.5 10.5){$n\!\!-\!\!2$}
\htext(-1.75 0.75){$n$}
\htext(-1.5 1.5){$n\!\!-\!\!2$} \vtext(-1.5 2.75){$\cdots$} \htext(-1.5 4){$2$}
\htext(-1.3 4.75){$0$} \htext(-1.75 5.25){$1$} \htext(-1.5 6){$2$}
\vtext(-1.5 7.25){$\cdots$} \htext(-1.5 8.5){$n\!\!-\!\!2$}
\htext(-1.75 9.75){$n$}
\htext(-1.5 10.5){$n\!\!-\!\!2$}
\htext(-3.75 0.75){$n$}
\htext(-3.5 1.5){$n\!\!-\!\!2$} \vtext(-3.5 2.75){$\cdots$} \htext(-3.5 4){$2$}
\htext(-3.3 4.75){$0$} \htext(-3.75 5.25){$1$} \htext(-3.5 6){$2$}
\vtext(-3.5 7.25){$\cdots$} \htext(-3.5 8.5){$n\!\!-\!\!2$}
\htext(-3.75 9.75){$n$}
\htext(-3.5 10.5){$n\!\!-\!\!2$}
\move(-1 0)\rlvec(0.5 0.5)\rmove(0.45 0.45)\rlvec(0.05 0.05)
\move(-2 0)\rlvec(0.05 0.05)\rmove(0.45 0.45)\rlvec(0.5 0.5)
\move(-3 0)\rlvec(0.5 0.5)\rmove(0.45 0.45)\rlvec(0.05 0.05)
\move(-4 0)\rlvec(0.05 0.05)\rmove(0.45 0.45)\rlvec(0.5 0.5)
\move(-1 9)\rlvec(0.5 0.5)\rmove(0.45 0.45)\rlvec(0.05 0.05)
\move(-2 9)\rlvec(0.05 0.05)\rmove(0.45 0.45)\rlvec(0.5 0.5)
\move(-3 9)\rlvec(0.5 0.5)\rmove(0.45 0.45)\rlvec(0.05 0.05)
\move(-4 9)\rlvec(0.05 0.05)\rmove(0.45 0.45)\rlvec(0.5 0.5)
\htext(-3.5 9.25){$n\!\!-\!\!1$} \htext(-3.5 0.25){$n\!\!-\!\!1$}
\htext(-2.5 9.75){$n\!\!-\!\!1$} \htext(-2.5 0.75){$n\!\!-\!\!1$}
\htext(-1.5 0.25){$n\!\!-\!\!1$} \htext(-1.5 9.25){$n\!\!-\!\!1$}
\htext(-0.5 9.75){$n\!\!-\!\!1$} \htext(-0.5 0.75){$n\!\!-\!\!1$}
\end{texdraw}}%
\savebox{\tmpfigg}{\begin{texdraw}
\fontsize{10}{10}\selectfont
\drawdim em
\setunitscale 1.9
\nc{\dtri}{
\bsegment
\move(-1 0)\lvec(0 1)\lvec(0 0)\lvec(-1 0)\ifill f:0.7
\esegment
}
\move(0 0)\dtri \move(-1 0)\dtri \move(-2 0)\dtri \move(-3 0)\dtri
\move(0 0)\rlvec(-4.3 0) \move(0 1)\rlvec(-4.3 0) \move(0 2)\rlvec(-4.3 0)
\move(0 3.5)\rlvec(-4.3 0) \move(0 4.5)\rlvec(-4.3 0)
\move(0 5.5)\rlvec(-4.3 0) \move(0 6.5)\rlvec(-4.3 0)
\move(0 8)\rlvec(-4.3 0) \move(0 9)\rlvec(-4.3 0)
\move(0 10)\rlvec(-4.3 0) \move(0 11)\rlvec(-4.3 0)
\move(0 0)\rlvec(0 11.3) \move(-1 0)\rlvec(0 11.3)
\move(-2 0)\rlvec(0 11.3) \move(-3 0)\rlvec(0 11.3) \move(-4 0)\rlvec(0 11.3)
\move(-1 4.5)\rlvec(1 1) \move(-2 4.5)\rlvec(1 1) \move(-3 4.5)\rlvec(1 1)
\move(-4 4.5)\rlvec(1 1)
\htext(-0.75 0.75){$n$}
\htext(-0.5 1.5){$n\!\!-\!\!2$} \vtext(-0.5 2.75){$\cdots$} \htext(-0.5 4){$2$}
\htext(-0.3 4.75){$1$} \htext(-0.75 5.25){$0$} \htext(-0.5 6){$2$}
\vtext(-0.5 7.25){$\cdots$} \htext(-0.5 8.5){$n\!\!-\!\!2$}
\htext(-0.75 9.75){$n$}
\htext(-0.5 10.5){$n\!\!-\!\!2$}
\htext(-2.75 0.75){$n$}
\htext(-2.5 1.5){$n\!\!-\!\!2$} \vtext(-2.5 2.75){$\cdots$} \htext(-2.5 4){$2$}
\htext(-2.3 4.75){$1$} \htext(-2.75 5.25){$0$} \htext(-2.5 6){$2$}
\vtext(-2.5 7.25){$\cdots$} \htext(-2.5 8.5){$n\!\!-\!\!2$}
\htext(-2.75 9.75){$n$}
\htext(-2.5 10.5){$n\!\!-\!\!2$}
\htext(-1.3 0.25){$n$}
\htext(-1.5 1.5){$n\!\!-\!\!2$} \vtext(-1.5 2.75){$\cdots$} \htext(-1.5 4){$2$}
\htext(-1.3 4.75){$0$} \htext(-1.75 5.25){$1$} \htext(-1.5 6){$2$}
\vtext(-1.5 7.25){$\cdots$} \htext(-1.5 8.5){$n\!\!-\!\!2$}
\htext(-1.3 9.25){$n$}
\htext(-1.5 10.5){$n\!\!-\!\!2$}
\htext(-3.3 0.25){$n$}
\htext(-3.5 1.5){$n\!\!-\!\!2$} \vtext(-3.5 2.75){$\cdots$} \htext(-3.5 4){$2$}
\htext(-3.3 4.75){$0$} \htext(-3.75 5.25){$1$} \htext(-3.5 6){$2$}
\vtext(-3.5 7.25){$\cdots$} \htext(-3.5 8.5){$n\!\!-\!\!2$}
\htext(-3.3 9.25){$n$}
\htext(-3.5 10.5){$n\!\!-\!\!2$}
\move(-2 0)\rlvec(0.5 0.5)\rmove(0.45 0.45)\rlvec(0.05 0.05)
\move(-1 0)\rlvec(0.05 0.05)\rmove(0.45 0.45)\rlvec(0.5 0.5)
\move(-4 0)\rlvec(0.5 0.5)\rmove(0.45 0.45)\rlvec(0.05 0.05)
\move(-3 0)\rlvec(0.05 0.05)\rmove(0.45 0.45)\rlvec(0.5 0.5)
\move(-2 9)\rlvec(0.5 0.5)\rmove(0.45 0.45)\rlvec(0.05 0.05)
\move(-1 9)\rlvec(0.05 0.05)\rmove(0.45 0.45)\rlvec(0.5 0.5)
\move(-4 9)\rlvec(0.5 0.5)\rmove(0.45 0.45)\rlvec(0.05 0.05)
\move(-3 9)\rlvec(0.05 0.05)\rmove(0.45 0.45)\rlvec(0.5 0.5)
\htext(-0.5 0.25){$n\!\!-\!\!1$} \htext(-0.5 9.25){$n\!\!-\!\!1$}
\htext(-2.5 0.25){$n\!\!-\!\!1$} \htext(-2.5 9.25){$n\!\!-\!\!1$}
\htext(-1.5 0.75){$n\!\!-\!\!1$} \htext(-1.5 9.75){$n\!\!-\!\!1$}
\htext(-3.5 0.75){$n\!\!-\!\!1$} \htext(-3.5 9.75){$n\!\!-\!\!1$}
\end{texdraw}}%
\savebox{\tmpfigh}{\begin{texdraw}
\fontsize{10}{10}\selectfont
\drawdim em
\setunitscale 1.9
\move(0 0)\lvec(-4 0)\lvec(-4 0.5)\lvec(0 0.5)\ifill f:0.7
\move(0 0)\rlvec(-4.3 0) \move(0 0.5)\rlvec(-4.3 0) \move(0 1)\rlvec(-4.3 0)
\move(0 2)\rlvec(-4.3 0) \move(0 3.5)\rlvec(-4.3 0) \move(0 4.5)\rlvec(-4.3 0)
\move(0 6)\rlvec(-4.3 0) \move(0 7)\rlvec(-4.3 0) \move(0 7.5)\rlvec(-4.3 0)
\move(0 8)\rlvec(-4.3 0) \move(0 9)\rlvec(-4.3 0) \move(0 0)\rlvec(0 9.3)
\move(-1 0)\rlvec(0 9.3) \move(-2 0)\rlvec(0 9.3) \move(-3 0)\rlvec(0 9.3)
\move(-4 0)\rlvec(0 9.3)
\htext(-0.5 0.25){$0$} \htext(-1.5 0.25){$0$} \htext(-2.5 0.25){$0$}
\htext(-3.5 0.25){$0$} \htext(-0.5 0.75){$0$} \htext(-1.5 0.75){$0$}
\htext(-2.5 0.75){$0$} \htext(-3.5 0.75){$0$}
\htext(-0.5 1.5){$1$} \htext(-1.5 1.5){$1$} \htext(-2.5 1.5){$1$}
\htext(-3.5 1.5){$1$}
\htext(-0.5 4){$n$} \htext(-1.5 4){$n$} \htext(-2.5 4){$n$}
\htext(-3.5 4){$n$}
\htext(-0.5 6.5){$1$} \htext(-1.5 6.5){$1$} \htext(-2.5 6.5){$1$}
\htext(-3.5 6.5){$1$}
\htext(-0.5 7.25){$0$} \htext(-1.5 7.25){$0$} \htext(-2.5 7.25){$0$}
\htext(-3.5 7.25){$0$} \htext(-0.5 7.75){$0$} \htext(-1.5 7.75){$0$}
\htext(-2.5 7.75){$0$} \htext(-3.5 7.75){$0$}
\htext(-0.5 8.5){$1$} \htext(-1.5 8.5){$1$} \htext(-2.5 8.5){$1$}
\htext(-3.5 8.5){$1$}
\vtext(-0.5 2.75){$\cdots$} \vtext(-1.5 2.75){$\cdots$}
\vtext(-2.5 2.75){$\cdots$} \vtext(-3.5 2.75){$\cdots$}
\vtext(-0.5 5.25){$\cdots$} \vtext(-1.5 5.25){$\cdots$}
\vtext(-2.5 5.25){$\cdots$} \vtext(-3.5 5.25){$\cdots$}
\end{texdraw}}%
\savebox{\tmpfigi}{\begin{texdraw}
\fontsize{10}{10}\selectfont
\drawdim em
\setunitscale 1.9
\move(0 0)\lvec(-4 0)\lvec(-4 0.5)\lvec(0 0.5)\ifill f:0.7
\move(0 0)\rlvec(-4.3 0) \move(0 0.5)\rlvec(-4.3 0) \move(0 1)\rlvec(-4.3 0)
\move(0 2)\rlvec(-4.3 0) \move(0 3.5)\rlvec(-4.3 0) \move(0 4.5)\rlvec(-4.3 0)
\move(0 4.5)\rlvec(-4.3 0) \move(0 5)\rlvec(-4.3 0) \move(0 5.5)\rlvec(-4.3 0)
\move(0 6.5)\rlvec(-4.3 0) \move(0 8)\rlvec(-4.3 0) \move(0 9)\rlvec(-4.3 0)
\move(0 9.5)\rlvec(-4.3 0) \move(0 10)\rlvec(-4.3 0) \move(0 11)\rlvec(-4.3 0)
\move(0 0)\rlvec(0 11.3) \move(-1 0)\rlvec(0 11.3) \move(-2 0)\rlvec(0 11.3)
\move(-3 0)\rlvec(0 11.3) \move(-4 0)\rlvec(0 11.3)
\htext(-0.5 0.25){$0$} \htext(-1.5 0.25){$0$} \htext(-2.5 0.25){$0$}
\htext(-3.5 0.25){$0$} \htext(-0.5 0.75){$0$} \htext(-1.5 0.75){$0$}
\htext(-2.5 0.75){$0$} \htext(-3.5 0.75){$0$}
\htext(-0.5 1.5){$1$} \htext(-1.5 1.5){$1$} \htext(-2.5 1.5){$1$}
\htext(-3.5 1.5){$1$}
\vtext(-0.5 2.75){$\cdots$} \vtext(-1.5 2.75){$\cdots$}
\vtext(-2.5 2.75){$\cdots$} \vtext(-3.5 2.75){$\cdots$}
\htext(-0.5 4){$n\!\!-\!\!1$} \htext(-1.5 4){$n\!\!-\!\!1$}
\htext(-2.5 4){$n\!\!-\!\!1$} \htext(-3.5 4){$n\!\!-\!\!1$}
\htext(-0.5 4.75){$n$} \htext(-1.5 4.75){$n$} \htext(-2.5 4.75){$n$}
\htext(-3.5 4.75){$n$} \htext(-0.5 5.25){$n$} \htext(-1.5 5.25){$n$}
\htext(-2.5 5.25){$n$} \htext(-3.5 5.25){$n$}
\htext(-0.5 6){$n\!\!-\!\!1$} \htext(-1.5 6){$n\!\!-\!\!1$}
\htext(-2.5 6){$n\!\!-\!\!1$} \htext(-3.5 6){$n\!\!-\!\!1$}
\vtext(-0.5 7.25){$\cdots$} \vtext(-1.5 7.25){$\cdots$}
\vtext(-2.5 7.25){$\cdots$} \vtext(-3.5 7.25){$\cdots$}
\htext(-0.5 8.5){$1$} \htext(-1.5 8.5){$1$} \htext(-2.5 8.5){$1$}
\htext(-3.5 8.5){$1$}
\htext(-0.5 9.25){$0$} \htext(-1.5 9.25){$0$} \htext(-2.5 9.25){$0$}
\htext(-3.5 9.25){$0$} \htext(-0.5 9.75){$0$} \htext(-1.5 9.75){$0$}
\htext(-2.5 9.75){$0$} \htext(-3.5 9.75){$0$}
\htext(-0.5 10.5){$1$} \htext(-1.5 10.5){$1$} \htext(-2.5 10.5){$1$}
\htext(-3.5 10.5){$1$}
\end{texdraw}}%
\savebox{\tmpfigj}{\begin{texdraw}
\fontsize{10}{10}\selectfont
\drawdim em
\setunitscale 1.9
\move(0 0)\lvec(-4 0)\lvec(-4 0.5)\lvec(0 0.5)\ifill f:0.7
\move(0 0)\rlvec(-4.3 0) \move(0 0.5)\rlvec(-4.3 0) \move(0 1)\rlvec(-4.3 0)
\move(0 2)\rlvec(-4.3 0) \move(0 3.5)\rlvec(-4.3 0) \move(0 4.5)\rlvec(-4.3 0)
\move(0 4.5)\rlvec(-4.3 0) \move(0 5)\rlvec(-4.3 0) \move(0 5.5)\rlvec(-4.3 0)
\move(0 6.5)\rlvec(-4.3 0) \move(0 8)\rlvec(-4.3 0) \move(0 9)\rlvec(-4.3 0)
\move(0 9.5)\rlvec(-4.3 0) \move(0 10)\rlvec(-4.3 0) \move(0 11)\rlvec(-4.3 0)
\move(0 0)\rlvec(0 11.3) \move(-1 0)\rlvec(0 11.3) \move(-2 0)\rlvec(0 11.3)
\move(-3 0)\rlvec(0 11.3) \move(-4 0)\rlvec(0 11.3)
\htext(-0.5 0.25){$n$} \htext(-1.5 0.25){$n$} \htext(-2.5 0.25){$n$}
\htext(-3.5 0.25){$n$} \htext(-0.5 0.75){$n$} \htext(-1.5 0.75){$n$}
\htext(-2.5 0.75){$n$} \htext(-3.5 0.75){$n$}
\htext(-0.5 1.5){$n\!\!-\!\!1$} \htext(-1.5 1.5){$n\!\!-\!\!1$}
\htext(-2.5 1.5){$n\!\!-\!\!1$} \htext(-3.5 1.5){$n\!\!-\!\!1$}
\vtext(-0.5 2.75){$\cdots$} \vtext(-1.5 2.75){$\cdots$}
\vtext(-2.5 2.75){$\cdots$} \vtext(-3.5 2.75){$\cdots$}
\htext(-0.5 4){$1$} \htext(-1.5 4){$1$}
\htext(-2.5 4){$1$} \htext(-3.5 4){$1$}
\htext(-0.5 4.75){$0$} \htext(-1.5 4.75){$0$} \htext(-2.5 4.75){$0$}
\htext(-3.5 4.75){$0$} \htext(-0.5 5.25){$0$} \htext(-1.5 5.25){$0$}
\htext(-2.5 5.25){$0$} \htext(-3.5 5.25){$0$}
\htext(-0.5 6){$1$} \htext(-1.5 6){$1$}
\htext(-2.5 6){$1$} \htext(-3.5 6){$1$}
\vtext(-0.5 7.25){$\cdots$} \vtext(-1.5 7.25){$\cdots$}
\vtext(-2.5 7.25){$\cdots$} \vtext(-3.5 7.25){$\cdots$}
\htext(-0.5 8.5){$n\!\!-\!\!1$} \htext(-1.5 8.5){$n\!\!-\!\!1$}
\htext(-2.5 8.5){$n\!\!-\!\!1$} \htext(-3.5 8.5){$n\!\!-\!\!1$}
\htext(-0.5 9.25){$n$} \htext(-1.5 9.25){$n$} \htext(-2.5 9.25){$n$}
\htext(-3.5 9.25){$n$} \htext(-0.5 9.75){$n$} \htext(-1.5 9.75){$n$}
\htext(-2.5 9.75){$n$} \htext(-3.5 9.75){$n$}
\htext(-0.5 10.5){$n\!\!-\!\!1$} \htext(-1.5 10.5){$n\!\!-\!\!1$}
\htext(-2.5 10.5){$n\!\!-\!\!1$} \htext(-3.5 10.5){$n\!\!-\!\!1$}
\end{texdraw}}%
\savebox{\tmpfigk}{\begin{texdraw}
\fontsize{10}{10}\selectfont
\drawdim em
\setunitscale 1.9
\nc{\dtri}{
\bsegment
\move(-1 0)\lvec(0 1)\lvec(0 0)\lvec(-1 0)\ifill f:0.7
\esegment
}
\move(0 0)\dtri \move(-1 0)\dtri \move(-2 0)\dtri \move(-3 0)\dtri
\move(0 0)\rlvec(-4.3 0) \move(0 1)\rlvec(-4.3 0) \move(0 2)\rlvec(-4.3 0)
\move(0 3.5)\rlvec(-4.3 0) \move(0 4.5)\rlvec(-4.3 0)
\move(0 5.5)\rlvec(-4.3 0) \move(0 6.5)\rlvec(-4.3 0)
\move(0 8)\rlvec(-4.3 0) \move(0 9)\rlvec(-4.3 0)
\move(0 10)\rlvec(-4.3 0) \move(0 11)\rlvec(-4.3 0)
\move(0 0)\rlvec(0 11.3) \move(-1 0)\rlvec(0 11.3) \move(-2 0)\rlvec(0 11.3)
\move(-3 0)\rlvec(0 11.3) \move(-4 0)\rlvec(0 11.3)
\move(-1 0)\rlvec(1 1) \move(-2 0)\rlvec(1 1) \move(-3 0)\rlvec(1 1)
\move(-4 0)\rlvec(1 1) \move(-1 9)\rlvec(1 1) \move(-2 9)\rlvec(1 1)
\move(-3 9)\rlvec(1 1) \move(-4 9)\rlvec(1 1)
\move(0 5)\rlvec(-4.3 0)
\htext(-0.3 0.25){$1$} \htext(-0.75 0.75){$0$} \htext(-0.5 1.5){$2$}
\vtext(-0.5 2.75){$\cdots$} \htext(-0.5 4){$n\!\!-\!\!1$}
\htext(-0.5 6){$n\!\!-\!\!1$} \htext(-0.5 8.5){$2$}
\htext(-0.3 9.25){$1$} \htext(-0.75 9.75){$0$} \htext(-0.5 10.5){$2$}
\htext(-2.3 0.25){$1$} \htext(-2.75 0.75){$0$} \htext(-2.5 1.5){$2$}
\vtext(-2.5 2.75){$\cdots$} \htext(-2.5 4){$n\!\!-\!\!1$}
\htext(-2.5 6){$n\!\!-\!\!1$} \htext(-2.5 8.5){$2$}
\htext(-2.3 9.25){$1$} \htext(-2.75 9.75){$0$} \htext(-2.5 10.5){$2$}
\htext(-1.3 0.25){$0$} \htext(-1.75 0.75){$1$} \htext(-1.5 1.5){$2$}
\vtext(-1.5 2.75){$\cdots$} \htext(-1.5 4){$n\!\!-\!\!1$}
\htext(-1.5 6){$n\!\!-\!\!1$} \htext(-1.5 8.5){$2$}
\htext(-1.3 9.25){$0$} \htext(-1.75 9.75){$1$} \htext(-1.5 10.5){$2$}
\htext(-3.3 0.25){$0$} \htext(-3.75 0.75){$1$} \htext(-3.5 1.5){$2$}
\vtext(-3.5 2.75){$\cdots$} \htext(-3.5 4){$n\!\!-\!\!1$}
\htext(-3.5 6){$n\!\!-\!\!1$} \htext(-3.5 8.5){$2$}
\htext(-3.3 9.25){$0$} \htext(-3.75 9.75){$1$} \htext(-3.5 10.5){$2$}
\htext(-0.5 4.75){$n$} \htext(-2.5 4.75){$n$} \htext(-1.5 5.25){$n$}
\htext(-3.5 5.25){$n$} \htext(-1.5 4.75){$n$} \htext(-3.5 4.75){$n$}
\htext(-0.5 5.25){$n$} \htext(-2.5 5.25){$n$}
\vtext(-0.5 7.25){$\cdots$} \vtext(-1.5 7.25){$\cdots$}
\vtext(-2.5 7.25){$\cdots$} \vtext(-3.5 7.25){$\cdots$}
\end{texdraw}}%
\savebox{\tmpfigl}{\begin{texdraw}
\fontsize{10}{10}\selectfont
\drawdim em
\setunitscale 1.9
\nc{\dtri}{
\bsegment
\move(-1 0)\lvec(0 1)\lvec(0 0)\lvec(-1 0)\ifill f:0.7
\esegment
}
\move(0 0)\dtri \move(-1 0)\dtri \move(-2 0)\dtri \move(-3 0)\dtri
\move(0 0)\rlvec(-4.3 0) \move(0 1)\rlvec(-4.3 0) \move(0 2)\rlvec(-4.3 0)
\move(0 3.5)\rlvec(-4.3 0) \move(0 4.5)\rlvec(-4.3 0)
\move(0 5.5)\rlvec(-4.3 0) \move(0 6.5)\rlvec(-4.3 0)
\move(0 8)\rlvec(-4.3 0) \move(0 9)\rlvec(-4.3 0)
\move(0 10)\rlvec(-4.3 0) \move(0 11)\rlvec(-4.3 0)
\move(0 0)\rlvec(0 11.3) \move(-1 0)\rlvec(0 11.3) \move(-2 0)\rlvec(0 11.3)
\move(-3 0)\rlvec(0 11.3) \move(-4 0)\rlvec(0 11.3)
\move(-1 0)\rlvec(1 1) \move(-2 0)\rlvec(1 1) \move(-3 0)\rlvec(1 1)
\move(-4 0)\rlvec(1 1) \move(-1 9)\rlvec(1 1) \move(-2 9)\rlvec(1 1)
\move(-3 9)\rlvec(1 1) \move(-4 9)\rlvec(1 1)
\move(0 5)\rlvec(-4.3 0)
\htext(-0.3 0.25){$0$} \htext(-0.75 0.75){$1$} \htext(-0.5 1.5){$2$}
\vtext(-0.5 2.75){$\cdots$} \htext(-0.5 4){$n\!\!-\!\!1$}
\htext(-0.5 6){$n\!\!-\!\!1$} \htext(-0.5 8.5){$2$}
\htext(-0.3 9.25){$0$} \htext(-0.75 9.75){$1$} \htext(-0.5 10.5){$2$}
\htext(-2.3 0.25){$0$} \htext(-2.75 0.75){$1$} \htext(-2.5 1.5){$2$}
\vtext(-2.5 2.75){$\cdots$} \htext(-2.5 4){$n\!\!-\!\!1$}
\htext(-2.5 6){$n\!\!-\!\!1$} \htext(-2.5 8.5){$2$}
\htext(-2.3 9.25){$0$} \htext(-2.75 9.75){$1$} \htext(-2.5 10.5){$2$}
\htext(-1.3 0.25){$1$} \htext(-1.75 0.75){$0$} \htext(-1.5 1.5){$2$}
\vtext(-1.5 2.75){$\cdots$} \htext(-1.5 4){$n\!\!-\!\!1$}
\htext(-1.5 6){$n\!\!-\!\!1$} \htext(-1.5 8.5){$2$}
\htext(-1.3 9.25){$1$} \htext(-1.75 9.75){$0$} \htext(-1.5 10.5){$2$}
\htext(-3.3 0.25){$1$} \htext(-3.75 0.75){$0$} \htext(-3.5 1.5){$2$}
\vtext(-3.5 2.75){$\cdots$} \htext(-3.5 4){$n\!\!-\!\!1$}
\htext(-3.5 6){$n\!\!-\!\!1$} \htext(-3.5 8.5){$2$}
\htext(-3.3 9.25){$1$} \htext(-3.75 9.75){$0$} \htext(-3.5 10.5){$2$}
\htext(-0.5 4.75){$n$} \htext(-2.5 4.75){$n$} \htext(-1.5 5.25){$n$}
\htext(-3.5 5.25){$n$} \htext(-1.5 4.75){$n$} \htext(-3.5 4.75){$n$}
\htext(-0.5 5.25){$n$} \htext(-2.5 5.25){$n$}
\vtext(-0.5 7.25){$\cdots$} \vtext(-1.5 7.25){$\cdots$}
\vtext(-2.5 7.25){$\cdots$} \vtext(-3.5 7.25){$\cdots$}
\end{texdraw}}%
\savebox{\tmpfigm}{\begin{texdraw}
\fontsize{10}{10}\selectfont
\drawdim em
\setunitscale 1.9
\move(0 0)\lvec(-4 0)\lvec(-4 0.5)\lvec(0 0.5)\ifill f:0.7
\move(0 0)\rlvec(-4.3 0) \move(0 1)\rlvec(-4.3 0) \move(0 2)\rlvec(-4.3 0)
\move(0 3.5)\rlvec(-4.3 0) \move(0 4.5)\rlvec(-4.3 0)
\move(0 5.5)\rlvec(-4.3 0) \move(0 6.5)\rlvec(-4.3 0) \move(0 8)\rlvec(-4.3 0)
\move(0 9)\rlvec(-4.3 0) \move(0 10)\rlvec(-4.3 0) \move(0 11)\rlvec(-4.3 0)
\move(0 0)\rlvec(0 11.3) \move(-1 0)\rlvec(0 11.3) \move(-2 0)\rlvec(0 11.3)
\move(-3 0)\rlvec(0 11.3) \move(-4 0)\rlvec(0 11.3)
\move(-1 4.5)\rlvec(1 1) \move(-2 4.5)\rlvec(1 1) \move(-3 4.5)\rlvec(1 1)
\move(-4 4.5)\rlvec(1 1)
\htext(-0.5 0.25){$n$}
\htext(-0.5 1.5){$n\!\!-\!\!1$} \vtext(-0.5 2.75){$\cdots$} \htext(-0.5 4){$2$}
\htext(-0.3 4.75){$1$} \htext(-0.75 5.25){$0$} \htext(-0.5 6){$2$}
\vtext(-0.5 7.25){$\cdots$} \htext(-0.5 8.5){$n\!\!-\!\!1$}
\htext(-0.5 9.25){$n$}
\htext(-0.5 10.5){$n\!\!-\!\!1$}
\htext(-2.5 0.25){$n$}
\htext(-2.5 1.5){$n\!\!-\!\!1$} \vtext(-2.5 2.75){$\cdots$} \htext(-2.5 4){$2$}
\htext(-2.3 4.75){$1$} \htext(-2.75 5.25){$0$} \htext(-2.5 6){$2$}
\vtext(-2.5 7.25){$\cdots$} \htext(-2.5 8.5){$n\!\!-\!\!1$}
\htext(-2.5 9.25){$n$}
\htext(-2.5 10.5){$n\!\!-\!\!1$}
\htext(-1.5 0.75){$n$}
\htext(-1.5 1.5){$n\!\!-\!\!1$} \vtext(-1.5 2.75){$\cdots$} \htext(-1.5 4){$2$}
\htext(-1.3 4.75){$0$} \htext(-1.75 5.25){$1$} \htext(-1.5 6){$2$}
\vtext(-1.5 7.25){$\cdots$} \htext(-1.5 8.5){$n\!\!-\!\!1$}
\htext(-1.5 9.75){$n$}
\htext(-1.5 10.5){$n\!\!-\!\!1$}
\htext(-3.5 0.75){$n$}
\htext(-3.5 1.5){$n\!\!-\!\!1$} \vtext(-3.5 2.75){$\cdots$} \htext(-3.5 4){$2$}
\htext(-3.3 4.75){$0$} \htext(-3.75 5.25){$1$} \htext(-3.5 6){$2$}
\vtext(-3.5 7.25){$\cdots$} \htext(-3.5 8.5){$n\!\!-\!\!1$}
\htext(-3.5 9.75){$n$}
\htext(-3.5 10.5){$n\!\!-\!\!1$}
\move(0 0.5)\rlvec(-4.3 0)
\move(0 9.5)\rlvec(-4.3 0)
\htext(-3.5 9.25){$n$} \htext(-3.5 0.25){$n$}
\htext(-2.5 9.75){$n$} \htext(-2.5 0.75){$n$}
\htext(-1.5 0.25){$n$} \htext(-1.5 9.25){$n$}
\htext(-0.5 9.75){$n$} \htext(-0.5 0.75){$n$}
\end{texdraw}}%

\item $A_n^{(1)}$ ($n\geq1$)\hfill
\vskip 3mm

      \begin{center}
      On $Y_{\La_i}$ : \raisebox{-\height}{\usebox{\tmpfiga}}
      \end{center}

\vskip 5mm

\item $A_{2n-1}^{(2)}$ ($n\geq3$)\hfill
\vskip 3mm

      \begin{center}
      On $Y_{\La_0}$ : \raisebox{-\height}{\usebox{\tmpfigb}}
      \qquad
      On $Y_{\La_1}$ : \raisebox{-\height}{\usebox{\tmpfigc}}
      \end{center}

\vskip 5mm

\item $D_n^{(1)}$ ($n\geq4$)\hfill
\vskip 5mm

      \begin{center}
      On $Y_{\La_0}$ : \raisebox{-\height}{\usebox{\tmpfigd}}
      \qquad\qquad
      On $Y_{\La_1}$ : \raisebox{-\height}{\usebox{\tmpfige}}
      \end{center}
      \vspace{7em}
      \begin{center}
      On $Y_{\La_{n-1}}$ : \raisebox{-\height}{\usebox{\tmpfigf}}
      \qquad\qquad
      On $Y_{\La_n}$ : \raisebox{-\height}{\usebox{\tmpfigg}}
      \end{center}

\vskip 5mm
\newpage
\item $A_{2n}^{(2)}$ ($n\geq2$)\hfill
\vskip 7mm

      \begin{center}
      On $Y_{\La_0}$ : \raisebox{-\height}{\usebox{\tmpfigh}}
      \end{center}

\vskip 1.5cm

\item $D_{n+1}^{(2)}$ ($n\geq2$)\hfill
\vskip 3mm

      \begin{center}
      On $Y_{\La_0}$ : \raisebox{-\height}{\usebox{\tmpfigi}}
      \qquad\qquad 
      On $Y_{\La_n}$ : \raisebox{-\height}{\usebox{\tmpfigj}}
      \end{center}

\vskip 5mm
\newpage

\item $B_{n}^{(1)}$ ($n\geq3$)\hfill
\vskip 3mm

      \begin{center}
      On $Y_{\La_0}$ : \raisebox{-\height}{\usebox{\tmpfigk}}
      \qquad\qquad 
      On $Y_{\La_1}$ : \raisebox{-\height}{\usebox{\tmpfigl}}
      \end{center}
      \vspace{1cm}
      \begin{center}
      On $Y_{\La_{n}}$ : \raisebox{-\height}{\usebox{\tmpfigm}}
      \qquad
      \phantom{On $Y_{\La_n}$ : \raisebox{-\height}{\usebox{\tmpfigl}}}
      \end{center}
\end{enumerate}
\end{enumerate}

\vskip 5mm
\newpage

\begin{defi}\hfill

(a) A wall built on the ground-state wall $Y_{\La_i}$ following the rules
listed above is called a {\it Young wall on the ground-state $\La_i$} if the
heights of its columns are weakly decreasing as we proceed from the right
to the left.

(b) A column in a Young wall is called a {\it full column} if its
height is a multiple of the unit length and its top is of unit
thickness.

(c) For the classical quantum affine algebras of type 
$A_{2n-1}^{(2)}$ $(n\ge 3)$, 
$D_n^{(1)}$ $(n\ge 4)$, $A_{2n}^{(2)}$ $(n\ge 2)$, 
$D_{n+1}^{(2)}$ $(n\ge 2)$ and $B_n^{(1)}$ $(n\ge 3)$, 
a Young wall is said to be {\it proper} if none of the full columns
have the same height. 

(d) For the quantum affine algebras of type $A_n^{(1)}$ $(n\ge 1)$,
every Young wall is defined to be {\it proper}. 

\end{defi}

\vskip 3mm 

\begin{example}
Consider the Young wall
$Y = (y_k)_{k=0}^{\infty} = (\cdots,y_k,\cdots,y_1,y_0)$
for $B_3^{(1)}$ built on the ground-state wall $Y_{\La_0}$:

\vskip 3mm

\begin{center}
\savebox{\tmppic}{\begin{texdraw}
\drawdim em
\setunitscale 1.7
\nc{\dtri}{
\bsegment
\move(-1 0)\lvec(0 1)\lvec(0 0)\lvec(-1 0)\ifill f:0.7
\esegment
}
\move(0 0)\dtri
\move(-1 0)\dtri
\move(-2 0)\dtri
\move(-3 0)\dtri
\move(-4 0)\dtri
\move(-5 0)\dtri
\move(0 0)\rlvec(-6 0)
\move(0 1)\rlvec(-6 0)
\move(0 2)\rlvec(-5 0)
\move(0 2.5)\rlvec(-5 0)
\move(0 3)\rlvec(-4 0)
\move(0 4)\rlvec(-3 0)
\move(0 5)\rlvec(-3 0)
\move(0 6)\rlvec(-1 0)
\move(0 0)\rlvec(0 6)
\move(-1 0)\rlvec(0 6)
\move(-2 0)\rlvec(0 5)
\move(-3 0)\rlvec(0 5)
\move(-4 0)\rlvec(0 3)
\move(-5 0)\rlvec(0 2.5)
\move(-6 0)\rlvec(0 1)
\move(-1 0)\rlvec(1 1)
\move(-2 0)\rlvec(1 1)
\move(-3 0)\rlvec(1 1)
\move(-4 0)\rlvec(1 1)
\move(-5 0)\rlvec(1 1)
\move(-6 0)\rlvec(1 1)
\move(-1 4)\rlvec(1 1)
\move(-2 4)\rlvec(1 1)
\move(-3 4)\rlvec(1 1)
\htext(-0.25 0.27){$1$}
\htext(-0.75 0.75){$0$}
\htext(-1.25 0.27){$0$}
\htext(-1.75 0.75){$1$}
\htext(-2.25 0.27){$1$}
\htext(-2.75 0.75){$0$}
\htext(-3.25 0.27){$0$}
\htext(-3.75 0.75){$1$}
\htext(-4.25 0.27){$1$}
\htext(-4.75 0.75){$0$}
\htext(-5.25 0.27){$0$}
\htext(-5.75 0.75){$1$}
\htext(-0.25 4.27){$1$}
\htext(-0.75 4.75){$0$}
\htext(-1.25 4.27){$0$}
\htext(-1.75 4.75){$1$}
%\htext(-2.25 4.27){$1$}
\htext(-2.75 4.75){$0$}
\htext(-0.5 1.5){$2$}
\htext(-1.5 1.5){$2$}
\htext(-2.5 1.5){$2$}
\htext(-3.5 1.5){$2$}
\htext(-4.5 1.5){$2$}
\htext(-0.5 3.5){$2$}
\htext(-1.5 3.5){$2$}
\htext(-2.5 3.5){$2$}
\htext(-0.5 2.25){$3$}
\htext(-1.5 2.25){$3$}
\htext(-2.5 2.25){$3$}
\htext(-3.5 2.25){$3$}
\htext(-4.5 2.25){$3$}
\htext(-0.5 2.75){$3$}
\htext(-1.5 2.75){$3$}
\htext(-2.5 2.75){$3$}
\htext(-3.5 2.75){$3$}
\htext(-0.5 5.5){$2$}
\end{texdraw}}
Y = \raisebox{-0.5\height}{\usebox{\tmppic}}
\end{center}

\vskip 3mm
\noindent
The columns $y_0$, $y_1$, $y_3$ $y_5$ are full columns and
$y_2$, $y_4$ are not full.
Hence $Y$ is a proper Young wall.

\vskip 3mm

On the other hand, 
\begin{equation*}
\savebox{\tmppic}{\begin{texdraw}
\drawdim em
\setunitscale 1.7
\nc{\dtri}{
\bsegment
\move(-1 0)\lvec(0 1)\lvec(0 0)\lvec(-1 0)\ifill f:0.7
\esegment
}
\move(0 0)\dtri
\move(-1 0)\dtri
\move(-2 0)\dtri
\move(-3 0)\dtri
\move(-4 0)\dtri
\move(-5 0)\dtri
\move(0 0)\rlvec(-6 0)
\move(0 1)\rlvec(-6 0)
\move(0 2)\rlvec(-5 0)
\move(0 2.5)\rlvec(-5 0)
\move(0 3)\rlvec(-5 0)
\move(0 4)\rlvec(-3 0)
\move(0 5)\rlvec(-3 0)
\move(0 0)\rlvec(0 5)
\move(-1 0)\rlvec(0 5)
\move(-2 0)\rlvec(0 5)
\move(-3 0)\rlvec(0 5)
\move(-4 0)\rlvec(0 3)
\move(-5 0)\rlvec(0 3)
\move(-6 0)\rlvec(0 1)
\move(-1 0)\rlvec(1 1)
\move(-2 0)\rlvec(1 1)
\move(-3 0)\rlvec(1 1)
\move(-4 0)\rlvec(1 1)
\move(-5 0)\rlvec(1 1)
\move(-6 0)\rlvec(1 1)
\move(-1 4)\rlvec(1 1)
\move(-2 4)\rlvec(1 1)
\move(-3 4)\rlvec(1 1)
\htext(-0.25 0.27){$1$}
\htext(-0.75 0.75){$0$}
\htext(-1.25 0.27){$0$}
\htext(-1.75 0.75){$1$}
\htext(-2.25 0.27){$1$}
\htext(-2.75 0.75){$0$}
\htext(-3.25 0.27){$0$}
\htext(-3.75 0.75){$1$}
\htext(-4.25 0.27){$1$}
\htext(-4.75 0.75){$0$}
\htext(-5.25 0.27){$0$}
\htext(-5.75 0.75){$1$}
\htext(-0.25 4.27){$1$}
\htext(-0.75 4.75){$0$}
\htext(-1.25 4.27){$0$}
\htext(-1.75 4.75){$1$}
%\htext(-2.25 4.27){$1$}
\htext(-2.75 4.75){$0$}
\htext(-0.5 1.5){$2$}
\htext(-1.5 1.5){$2$}
\htext(-2.5 1.5){$2$}
\htext(-3.5 1.5){$2$}
\htext(-4.5 1.5){$2$}
\htext(-0.5 3.5){$2$}
\htext(-1.5 3.5){$2$}
\htext(-2.5 3.5){$2$}
\htext(-0.5 2.25){$3$}
\htext(-1.5 2.25){$3$}
\htext(-2.5 2.25){$3$}
\htext(-3.5 2.25){$3$}
\htext(-4.5 2.25){$3$}
\htext(-0.5 2.75){$3$}
\htext(-1.5 2.75){$3$}
\htext(-2.5 2.75){$3$}
\htext(-3.5 2.75){$3$}
\htext(-4.5 2.75){$3$}
\end{texdraw}}
Y' = \raisebox{-0.5\height}{\usebox{\tmppic}}\quad
\end{equation*}
is not a proper Young wall because the full columns $y'_0$ and $y'_1$
(also $y'_3$ and $y'_4$) have the same height.
\end{example}

\vskip 1cm

\section{The crystal structure}

Let $\la$ be a dominant integral weight of level 1
and let $\F(\la)$ denote the set of all proper Young walls built
on the ground-state $Y_\la$.
In this section, we define an affine crystal structure on $\F(\la)$.
The action of Kashiwara operators are defined using the 
$i$-signatures of proper Young walls in a similar way 
to playing the Tetris game. 

\vskip 3mm 

\begin{defi}\hfill

(a) A block of color $i$  in a proper Young wall is called 
a {\it removable $i$-block}
if the wall remains a proper Young wall after removing the block.
A column in a proper Young wall is called {\it $i$-removable}
if the top of that column is a removable $i$-block.

(b) A place in a proper Young wall where one may add an
$i$-block to obtain another proper Young wall is called an
{\it $i$-admissible slot}.
A column in a proper Young wall is called
{\it $i$-admissible} if the top of that column is an
$i$-admissible slot.

\end{defi}

\begin{example}
In the following figure, we 
consider a proper Young wall for $B_3^{(1)}$ built on the
ground-state wall $Y_{\La_0}$
and indicate all the removable blocks and admissible slots.

\vskip 3mm

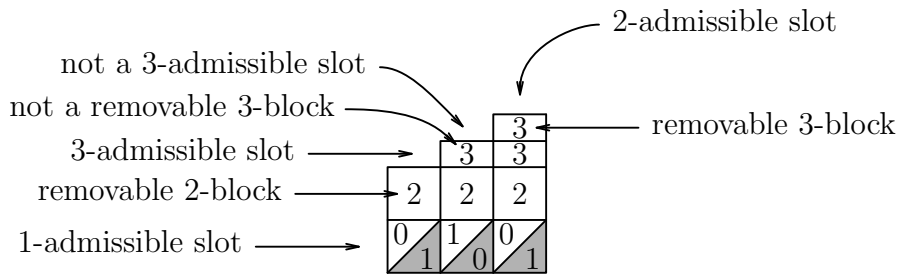
\begin{center}
\begin{texdraw}
\drawdim em
\setunitscale 1.7
\nc{\dtri}{
\bsegment
\move(-1 0)\lvec(0 1)\lvec(0 0)\lvec(-1 0)\ifill f:0.7
\esegment
}
\move(0 0)\dtri
\move(-1 0)\dtri
\move(-2 0)\dtri
\move(-3 0)\lvec(0 0)\lvec(0 3)
\move(-3 0)\lvec(-3 2)\lvec(0 2)
\move(-2 0)\lvec(-2 2.5)\lvec(0 2.5)
\move(-1 0)\lvec(-1 3)\lvec(0 3)
\move(-3 1)\lvec(0 1)
\move(-2 0)\rlvec(1 1)
\move(-1 0)\rlvec(1 1)
\move(-3 0)\rlvec(1 1)
\htext(-0.5 2.75){$3$}
\htext(-0.5 2.25){$3$}
\htext(-1.5 2.25){$3$}
\htext(-0.5 1.5){$2$}
\htext(-1.5 1.5){$2$}
\htext(-2.5 1.5){$2$}
\htext(-0.25 0.27){$1$}
\htext(-0.75 0.75){$0$}
\htext(-1.25 0.27){$0$}
\htext(-1.75 0.75){$1$}
\htext(-2.25 0.27){$1$}
\htext(-2.75 0.75){$0$}
\arrowheadtype t:V
\textref h:L v:C
\htext(2 2.85){removable 3-block}
    \move(1.75 2.75)\avec(-0.18 2.75)
\htext(1.25 4.8){2-admissible slot}
    \move(1 4.7)
    \clvec(0.5 4.7)(-0.5 4.25)(-0.5 3.5)
    \move(-0.5 3.5)\ahead{0}{-1}
\textref h:R v:C
\htext(-3.4 4){not a 3-admissible slot}
    \move(-3.15 3.9)\clvec(-2.65 3.9)(-2.2 3.4)(-1.7 2.9)
    \move(-1.66 2.86)\ahead{1}{-1}
\htext(-4 3.2){not a removable 3-block}
    \move(-3.8 3.1)\clvec(-3.25 3.1)(-2.34 3.06)(-1.84 2.56)
    \move(-1.8 2.52)\ahead{1}{-1}
\htext(-4.75 2.35){3-admissible slot}
    \move(-4.5 2.25)\avec(-2.55 2.25)
\htext(-5 1.6){removable 2-block}
    \move(-4.75 1.5)\avec(-2.82 1.5)
\htext(-5.75 0.6){1-admissible slot}
    \move(-5.5 0.5)\avec(-3.55 0.5)
\move(6.8 5.1)\move(-10.5 -0.1)
\end{texdraw}
\end{center}
\end{example}

\vskip 3mm

We now define the action of Kashiwara operators $\eit$, $\fit$ ($i\in I$)
on $\F(\la)$.
Fix $i\in I$ and let $Y= (y_k)_{k=0}^{\infty}\in\F(\la)$ be a
proper Young wall.

\begin{enumerate}
\item To each column $y_k$ of $Y$, we assign its {\it $i$-signature} as follows:
\begin{enumerate}
\item we assign $-\,-$ if the column $y_k$ is twice $i$-removable;
\item we assign $-$ if the column is once $i$-removable,
      but not $i$-admissible (the $i$-block may be of unit height or
      of half-unit height);
\item we assign $-\,+$ if the column is once $i$-removable and
      once $i$-admissible (the $i$-block will be of half-unit height
      in this case);
\item we assign $+$ if the column is once $i$-admissible, but not
      $i$-removable (the $i$-block may be of unit height or
      of half-unit height);
\item we assign $+\,+$ if the column is twice $i$-admissible
      (the $i$-block will be of half-unit height in this case).
\end{enumerate}
\item From the (infinite) sequence of $+$'s and $-$'s, cancel out every
$(+,-)$-pair to obtain a finite sequence of $-$'s followed by $+$'s,
reading from left to right.
This sequence is called the {\it $i$-signature} 
of the proper Young wall $Y$.
\item We define $\eit Y$ to be the proper Young wall obtained from $Y$
by removing the $i$-block corresponding to the right-most $-$ in the
$i$-signature of $Y$.
We define $\eit Y = 0$ if there exists no $-$ in the $i$-signature of $Y$.
\item We define $\fit Y$ to be the proper Young wall obtained from $Y$
by adding an $i$-block to the column corresponding to the left-most $+$
in the $i$-signature of $Y$.
We define $\fit Y = 0$ if there exists no $+$ in the $i$-signature of $Y$.
\end{enumerate}

We also define the maps
\begin{equation*}
\wt : \F(\la) \longrightarrow P,\quad 
\vep_i : \F(\la) \longrightarrow \Z, \quad 
\vphi_i : \F(\la) \longrightarrow \Z
\end{equation*}
by
\begin{equation*}
\begin{aligned}\mbox{}
\wt(Y) &= \la - \sum_{i\in I} k_i \ali,\\
\vep_i(Y) &= \text{the number of $-$'s in the $i$-signature of $Y$,}\\
\vphi_i(Y) &= \text{the number of $+$'s in the $i$-signature of $Y$,}
\end{aligned}
\end{equation*}
where $k_i$ is the number of $i$-blocks in $Y$ that have been added to
the ground-state wall $Y_\la$.

\vskip 3mm 
Then it is straightforward to verify that the following theorem holds. 

\vskip 3mm 
\begin{thm}\label{thm:933}
The maps $\wt : \F(\la) \ra P$,
$\eit,\fit : \F(\la) \ra \F(\la)\cup \{0\}$ 
and $\vep_i,\vphi_i : \F(\la) \ra \Z$ define an affine crystal
structure on the set $\F(\la)$ of all proper Young walls.
\end{thm}

\vskip 3mm

Let $\delta$ be the null root for the quantum affine algebra $U_q(\g)$
and write
\begin{equation*}
\begin{cases}
\delta = a_0\alpha_0 + a_1\alpha_1 + \cdots + a_n\alpha_n
&\text{for $\g = A_n^{(1)}, \cdots, B_n^{(1)}$},\\
2\delta = a_0\alpha_0 + a_1\alpha_1 + \cdots + a_n\alpha_n
&\text{for $\g = D_{n+1}^{(2)}$}
\end{cases}
\end{equation*}
The part of a column consisting of $a_0$-many 0-blocks,
$a_1$-many 1-blocks, $\cdots$, $a_n$-many $n$-blocks in some cyclic
order is called a {\it $\delta$-column}.

\vskip 3mm

\begin{example}\hfill
\begin{enumerate}
\item The following are $\delta$-columns for $B_3^{(1)}$.

\vskip 3mm

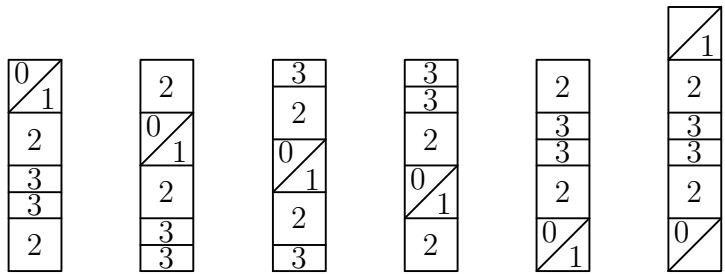
\begin{center}
\begin{texdraw}
\drawdim em
\setunitscale 1.7
\move(0 0)
\bsegment
\move(0 0)\lvec(1 0)\lvec(1 4)\lvec(0 4)\lvec(0 0)
\move(0 1)\lvec(1 1)
\move(0 1.5)\lvec(1 1.5)
\move(0 2)\lvec(1 2)
\move(0 3)\lvec(1 3)
\move(0 3)\lvec(1 4)
\htext(0.5 0.5){$2$}
\htext(0.5 1.25){$3$}
\htext(0.5 1.75){$3$}
\htext(0.5 2.5){$2$}
\htext(0.25 3.75){$0$}
\htext(0.75 3.27){$1$}
\esegment
\move(2.5 0)
\bsegment
\move(0 0)\lvec(1 0)\lvec(1 4)\lvec(0 4)\lvec(0 0)
\move(0 0.5)\lvec(1 0.5)
\move(0 1)\lvec(1 1)
\move(0 2)\lvec(1 2)
\move(0 3)\lvec(1 3)
\move(0 2)\lvec(1 3)
\htext(0.5 3.5){$2$}
\htext(0.5 0.25){$3$}
\htext(0.5 0.75){$3$}
\htext(0.5 1.5){$2$}
\htext(0.25 2.75){$0$}
\htext(0.75 2.27){$1$}
\esegment
\move(5 0)
\bsegment
\move(0 0)\lvec(1 0)\lvec(1 4)\lvec(0 4)\lvec(0 0)
\move(0 0.5)\lvec(1 0.5)
\move(0 1.5)\lvec(1 1.5)
\move(0 2.5)\lvec(1 2.5)
\move(0 3.5)\lvec(1 3.5)
\move(0 1.5)\lvec(1 2.5)
\htext(0.5 1){$2$}
\htext(0.5 0.25){$3$}
\htext(0.5 3.75){$3$}
\htext(0.5 3){$2$}
\htext(0.75 1.72){$1$}
\htext(0.25 2.25){$0$}
\esegment
\move(7.5 0)
\bsegment
\move(0 0)\lvec(1 0)\lvec(1 4)\lvec(0 4)\lvec(0 0)
\move(0 3.5)\lvec(1 3.5)
\move(0 1)\lvec(1 1)
\move(0 2)\lvec(1 2)
\move(0 3)\lvec(1 3)
\move(0 1)\lvec(1 2)
\htext(0.5 2.5){$2$}
\htext(0.5 3.25){$3$}
\htext(0.5 3.75){$3$}
\htext(0.5 0.5){$2$}
\htext(0.25 1.75){$0$}
\htext(0.75 1.27){$1$}
\esegment
\move(10 0)
\bsegment
\move(0 0)\lvec(1 0)\lvec(1 4)\lvec(0 4)\lvec(0 0)
\move(0 2.5)\lvec(1 2.5)
\move(0 1)\lvec(1 1)
\move(0 2)\lvec(1 2)
\move(0 3)\lvec(1 3)
\move(0 0)\lvec(1 1)
\htext(0.5 3.5){$2$}
\htext(0.5 2.25){$3$}
\htext(0.5 2.75){$3$}
\htext(0.5 1.5){$2$}
\htext(0.25 0.75){$0$}
\htext(0.75 0.27){$1$}
\esegment
\move(12.5 0)
\bsegment
\move(0 0)\lvec(1 0)\lvec(1 4)\lvec(0 4)\lvec(0 0)
\move(0 2.5)\lvec(1 2.5)
\move(0 1)\lvec(1 1)
\move(0 2)\lvec(1 2)
\move(0 3)\lvec(1 3)
\move(0 0)\lvec(1 1)
\htext(0.5 3.5){$2$}
\htext(0.5 2.25){$3$}
\htext(0.5 2.75){$3$}
\htext(0.5 1.5){$2$}
\htext(0.25 0.75){$0$}
\htext(0.75 4.27){$1$}
\move(1 4)\lvec(1 5)\lvec(0 5)\lvec(0 4)\lvec(1 5)
\esegment
\end{texdraw}
\end{center}

\vskip 5mm

\item The following are $\delta$-columns for $D_3^{(2)}$. 

\vskip 3mm

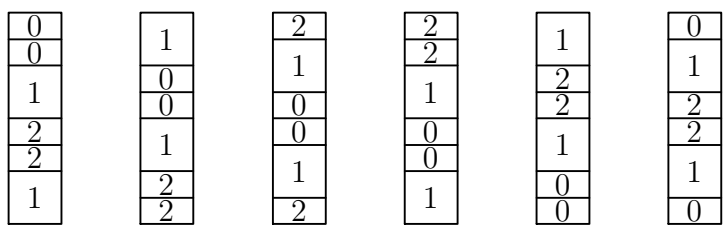
\begin{center}
\begin{texdraw}
\drawdim em
\setunitscale 1.7
\move(0 0)
\bsegment
\move(0 0)\lvec(1 0)\lvec(1 4)\lvec(0 4)\lvec(0 0)
\move(0 1)\rlvec(1 0)
\move(0 2)\rlvec(1 0)
\move(0 3)\rlvec(1 0)
\move(0 1.5)\rlvec(1 0)
\move(0 3.5)\rlvec(1 0)
\htext(0.5 0.5){$1$}
\htext(0.5 1.25){$2$}
\htext(0.5 1.75){$2$}
\htext(0.5 2.5){$1$}
\htext(0.5 3.25){$0$}
\htext(0.5 3.75){$0$}
\esegment
\move(2.5 0)
\bsegment
\move(0 0)\lvec(1 0)\lvec(1 4)\lvec(0 4)\lvec(0 0)
\move(0 1)\rlvec(1 0)
\move(0 2)\rlvec(1 0)
\move(0 3)\rlvec(1 0)
\move(0 0.5)\rlvec(1 0)
\move(0 2.5)\rlvec(1 0)
\htext(0.5 3.5){$1$}
\htext(0.5 0.25){$2$}
\htext(0.5 0.75){$2$}
\htext(0.5 1.5){$1$}
\htext(0.5 2.25){$0$}
\htext(0.5 2.75){$0$}
\esegment
\move(5 0)
\bsegment
\move(0 0)\lvec(1 0)\lvec(1 4)\lvec(0 4)\lvec(0 0)
\move(0 2.5)\rlvec(1 0)
\move(0 2)\rlvec(1 0)
\move(0 3.5)\rlvec(1 0)
\move(0 0.5)\rlvec(1 0)
\move(0 1.5)\rlvec(1 0)
\htext(0.5 3){$1$}
\htext(0.5 3.75){$2$}
\htext(0.5 0.25){$2$}
\htext(0.5 1){$1$}
\htext(0.5 2.25){$0$}
\htext(0.5 1.75){$0$}
\esegment
\move(7.5 0)
\bsegment
\move(0 0)\lvec(1 0)\lvec(1 4)\lvec(0 4)\lvec(0 0)
\move(0 1)\rlvec(1 0)
\move(0 2)\rlvec(1 0)
\move(0 3)\rlvec(1 0)
\move(0 1.5)\rlvec(1 0)
\move(0 3.5)\rlvec(1 0)
\htext(0.5 0.5){$1$}
\htext(0.5 1.25){$0$}
\htext(0.5 1.75){$0$}
\htext(0.5 2.5){$1$}
\htext(0.5 3.25){$2$}
\htext(0.5 3.75){$2$}
\esegment
\move(10 0)
\bsegment
\move(0 0)\lvec(1 0)\lvec(1 4)\lvec(0 4)\lvec(0 0)
\move(0 1)\rlvec(1 0)
\move(0 2)\rlvec(1 0)
\move(0 3)\rlvec(1 0)
\move(0 0.5)\rlvec(1 0)
\move(0 2.5)\rlvec(1 0)
\htext(0.5 3.5){$1$}
\htext(0.5 0.25){$0$}
\htext(0.5 0.75){$0$}
\htext(0.5 1.5){$1$}
\htext(0.5 2.25){$2$}
\htext(0.5 2.75){$2$}
\esegment
\move(12.5 0)
\bsegment
\move(0 0)\lvec(1 0)\lvec(1 4)\lvec(0 4)\lvec(0 0)
\move(0 2.5)\rlvec(1 0)
\move(0 2)\rlvec(1 0)
\move(0 3.5)\rlvec(1 0)
\move(0 0.5)\rlvec(1 0)
\move(0 1.5)\rlvec(1 0)
\htext(0.5 3){$1$}
\htext(0.5 3.75){$0$}
\htext(0.5 0.25){$0$}
\htext(0.5 1){$1$}
\htext(0.5 2.25){$2$}
\htext(0.5 1.75){$2$}
\esegment
\end{texdraw}
\end{center}
\end{enumerate}
\end{example}

\vskip 3mm

\begin{defi}\hfill

(a) A column in a proper Young wall is said to
{\it contain a removable $\delta$} if we may remove a $\delta$-column
from $Y$ and still obtain a proper Young wall.

(b) A proper Young wall is said to be {\it reduced} if none of its
columns contain a removable $\delta$.

\end{defi}

\vskip 3mm
\begin{example}
The first Young wall drawn below is reduced, but the second one is not
reduced.

\vskip 3mm

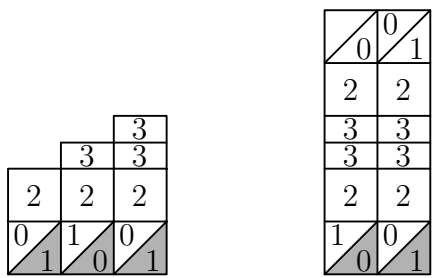
\begin{center}
\begin{texdraw}
\drawdim em
\setunitscale 1.7
\nc{\dtri}{
\bsegment
\move(-1 0)\lvec(0 1)\lvec(0 0)\lvec(-1 0)\ifill f:0.7
\esegment
}
\move(0 0)
\bsegment
\move(1 0)\dtri
\move(2 0)\dtri
\move(3 0)\dtri
\move(3 2)\lvec(0 2)\lvec(0 0)\lvec(3 0)\lvec(3 3)\lvec(2 3)\lvec(2 0)
\lvec(3 1)\lvec(0 1)
\move(1 0)\lvec(1 2.5)\lvec(3 2.5)
\move(0 0)\lvec(1 1)
\move(1 0)\lvec(2 1)
\htext(0.5 1.5){$2$}
\htext(1.5 1.5){$2$}
\htext(2.5 1.5){$2$}
\htext(1.5 2.25){$3$}
\htext(2.5 2.25){$3$}
\htext(2.5 2.75){$3$}
\htext(0.25 0.75){$0$}
\htext(0.75 0.27){$1$}
\htext(1.25 0.75){$1$}
\htext(1.75 0.27){$0$}
\htext(2.25 0.75){$0$}
\htext(2.75 0.27){$1$}
\esegment
\move(6 0)
\bsegment
\move(1 0)\dtri
\move(2 0)\dtri
\move(0 0)\lvec(2 0)\lvec(2 5)\lvec(0 5)\lvec(0 0)\lvec(1 1)
\move(1 0)\lvec(2 1)\lvec(0 1)
\move(1 0)\lvec(1 5)\lvec(0 4)\lvec(2 4)
\move(1 4)\lvec(2 5)
\move(0 2)\rlvec(2 0)
\move(0 2.5)\rlvec(2 0)
\move(0 3)\rlvec(2 0)
\htext(0.25 0.75){$1$}
\htext(0.75 0.27){$0$}
\htext(1.25 0.75){$0$}
\htext(1.75 0.27){$1$}
\htext(0.5 1.5){$2$}
\htext(1.5 1.5){$2$}
\htext(1.5 2.25){$3$}
\htext(1.5 2.75){$3$}
\htext(0.5 2.25){$3$}
\htext(0.5 2.75){$3$}
\htext(0.5 3.5){$2$}
\htext(1.5 3.5){$2$}
\htext(0.75 4.27){$0$}
\htext(1.25 4.75){$0$}
\htext(1.75 4.27){$1$}
\esegment
\end{texdraw}
\end{center}
\end{example}

\vskip 3mm
Let $\Y(\la)$ be the set of all reduced proper Young walls
built on the ground-state wall $Y_\la$.
We first prove that the set $\Y(\la)$ is stable under the 
Kashiwara operators. 

\begin{thm}\label{thm:closed}
For all $i\in I$ and $Y\in \Y(\la)$, we have 
$$\eit Y\in \Y(\la)\cup \{0\} \quad \text{and} \quad 
\fit Y\in \Y(\la)\cup \{0\}.$$
Hence the set $\Y(\la)$ has an affine crystal structure for the
quantum affine algebra $U_q(\g)$.
\end{thm}
\begin{proof}
Let $Y\in \Y(\la)$ and fix an index $i\in I$.
If $\eit Y\neq 0$, by definition of removable blocks, $\eit Y$ is a proper
Young wall.
Suppose that $\eit Y$ is not reduced.
This means that removing an $i$-block for $Y$ has created a removable
$\delta$ in $\eit Y$.
Hence $Y$ must have a column of the form:

\vskip 5mm

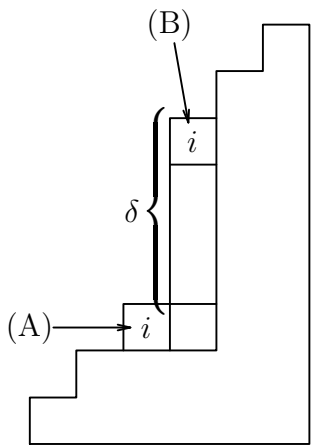
\begin{center}
\begin{texdraw}
\drawdim em
\setunitscale 1.5
%\move(0 0)\lvec(0 1)\lvec(1 1)\lvec(1 2)\lvec(2 2)
\move(2 2)\lvec(2 3)\lvec(3 3)
\lvec(3 4)\lvec(6 4)\lvec(6 10)\lvec(7 10)\lvec(7 11)\lvec(8 11)
%\lvec(8 0)\lvec(0 0)
\lvec(8 2)\lvec(2 2)
\move(4 4)\lvec(4 5)\lvec(6 5)
\move(5 4)\lvec(5 9)\lvec(6 9)
\move(5 8)\lvec(6 8)
\vtext(4.7 7){$\overbrace{\rule{6.6em}{0em}}$}
\htext(4.2 7){$\delta$}
\htext(4.5 4.5){$i$}
\htext(5.5 8.5){$i$}
\htext(5 11){(B)}
\htext(2 4.5){(A)}
\arrowheadtype t:V
\move(2.5 4.5)\avec(4.15 4.5)
\move(5.1 10.6)\avec(5.4 8.9)
\end{texdraw}
\end{center}

\vskip 3mm
Here, the $i$-block at (A) is removed to give $\eit Y$.
Then, in this case, the rules for the action of Kashiwara operators tell us
that the operator $\eit$ would remove the $i$-block at (B), not the one at
(A), which is a contradiction.
Hence $\eit Y$ is reduced.

Similarly, if $\fit Y\neq 0$, one can show that 
$\fit Y$ is a reduced proper Young wall.
\end{proof}

\vskip 1cm 
%%%%%%%%%%%%%%%%%%%%%%%%%%%%%%%%%%%%%%%%%%%%%%%%%%%%%%%%%%%
\section{Crystal graphs for basic representations}

{% This brace is placed here and at the end of this section to limit
%  the scope of \tmpfig's. Do not erase.
%FIGSTART

In this section, we will prove that the crystal graph $B(\la)$ for
the basic representation $V(\la)$ of the quantum affine algebra
$U_q(\g)$ is isomorphic to the affine crystal $\Y(\la)$
consisting of all reduced
proper Young walls built on the ground-state wall $Y_\la$.
Hence the crystal $\Y(\la)$ would be  the connected component of
the crystal $\F(\la)$
contaning the ground-state wall $Y_\la$.

\vskip 3mm

\begin{thm}\label{thm:main}
There exists an isomorphism of $U_q(\g)$-crystals
\begin{equation}
\Y(\la) \stackrel{\sim} \longrightarrow B(\la)
\quad \text{given by} \ \ 
Y_{\la} \longmapsto u_{\la},
\end{equation}
where $u_\la$ is the highest weight vector in $B(\la)$.
\end{thm}

\vskip 5mm 
\noindent
{\it Proof.} \ 
The rest of this section is devoted to proving this theorem.
It suffices to show $\Y(\la)\cong  \bP(\la)$ as $U'_q(\g)$-crystals.
Let us define a map $\Psi : \Y(\la)\longrightarrow  \bP(\la)$
as follows: \ 
We read the top parts of each column in a reduced proper Young wall
and write down a sequence of elements from $\B$, the perfect crystal of
level 1, to obtain a $\la$-path.
The rules for obtaining the paths are given below for each affine type
in question.

\newpage
\savebox{\tmpfiga}{\begin{texdraw}
\fontsize{8}{8}\selectfont
\drawdim em
\setunitscale 1.7
\move(-1 0)\lvec(-1 1)\lvec(0 1)\lvec(-1 0)\lvec(0 0)\lvec(0 1)
\htext(-0.7 0.75){$0$}
\end{texdraw}}%
\savebox{\tmpfigb}{\begin{texdraw}
\fontsize{8}{8}\selectfont
\drawdim em
\setunitscale 1.7
\move(-1 0)\lvec(-1 1)\lvec(0 1)\lvec(-1 0)\lvec(0 0)\lvec(0 1)
\htext(-0.3 0.3){$0$}
\end{texdraw}}%
\savebox{\tmpfigc}{\begin{texdraw}
\fontsize{8}{8}\selectfont
\drawdim em
\setunitscale 1.7
\move(0 0)\lvec(-1 0)\lvec(-1 1)\lvec(0 1)\lvec(0 0)
\htext(-0.5 0.5){$1$}
\end{texdraw}}%
\savebox{\tmpfigd}{\begin{texdraw}
\fontsize{8}{8}\selectfont
\drawdim em
\setunitscale 1.7
\move(-1 0)\lvec(-1 1)\lvec(0 1)\lvec(-1 0)\lvec(0 0)\lvec(0 1)
\htext(-0.7 0.75){$1$}
\end{texdraw}}%
\savebox{\tmpfige}{\begin{texdraw}
\fontsize{8}{8}\selectfont
\drawdim em
\setunitscale 1.7
\move(-1 0)\lvec(-1 1)\lvec(0 1)\lvec(-1 0)\lvec(0 0)\lvec(0 1)
\htext(-0.3 0.3){$1$}
\end{texdraw}}%
\savebox{\tmpfigf}{\begin{texdraw}
\fontsize{8}{8}\selectfont
\drawdim em
\setunitscale 1.7
\move(0 0)\lvec(-1 0)\lvec(-1 1)\lvec(0 1)\lvec(0 0)
\htext(-0.5 0.5){$\bar{1}$}
\end{texdraw}}%
\savebox{\tmpfigg}{\begin{texdraw}
\fontsize{8}{8}\selectfont
\drawdim em
\setunitscale 1.7
\lpatt(0.1 0.15)
\move(0 1)\lvec(0 2)\lvec(-1 2)\lvec(-1 1)
\move(-1 1)\lvec(0 2)
\lpatt(1 0)
\move(0 0)\lvec(-1 0)\lvec(-1 1)\lvec(0 1)\lvec(0 0)
\htext(-0.5 0.5){$2$}
\end{texdraw}}%
\savebox{\tmpfigh}{\begin{texdraw}
\fontsize{8}{8}\selectfont
\drawdim em
\setunitscale 1.7
\move(-1 0)\lvec(-1 2)\lvec(0 2)\lvec(0 0)\lvec(-1 0)
\move(-1 0)\lvec(0 1)\lvec(-1 1)
\htext(-0.5 1.5){$2$}
\end{texdraw}}%
\savebox{\tmpfigi}{\begin{texdraw}
\fontsize{8}{8}\selectfont
\drawdim em
\setunitscale 1.7
\move(0 0)\lvec(-1 0)\lvec(-1 1)\lvec(0 1)\lvec(0 0)
\htext(-0.5 0.5){$2$}
\end{texdraw}}%
\savebox{\tmpfigj}{\begin{texdraw}
\fontsize{8}{8}\selectfont
\drawdim em
\setunitscale 1.7
\lpatt(0.1 0.15)
\move(-1 1)\lvec(-1 2)\lvec(0 2)\lvec(0 1)
\lpatt(1 0)
\move(0 0)\lvec(-1 0)\lvec(-1 1)\lvec(0 1)\lvec(0 0)
\htext(-0.5 0.5){$3$}
\htext(-0.5 1.5){$2$}
\end{texdraw}}%
\savebox{\tmpfigk}{\begin{texdraw}
\fontsize{8}{8}\selectfont
\drawdim em
\setunitscale 1.7
\move(0 0)\lvec(-1 0)\lvec(-1 1)\lvec(0 1)\lvec(0 0)
\htext(-0.5 0.5){$\bar{2}$}
\end{texdraw}}%
\savebox{\tmpfigl}{\begin{texdraw}
\fontsize{8}{8}\selectfont
\drawdim em
\setunitscale 1.7
\move(0 0)\lvec(-1 0)
\move(0 1)\lvec(-1 1)\move(0 0)\lvec(0 2)\lvec(-1 2)\lvec(-1 0)
\htext(-0.5 0.5){$j\!\!-\!\!2$}
\htext(-0.5 1.5){$j\!\!-\!\!1$}
\end{texdraw}}%
\savebox{\tmpfigm}{\begin{texdraw}
\fontsize{8}{8}\selectfont
\drawdim em
\setunitscale 1.7
\move(0 0)\lvec(-1 0)\lvec(-1 1)\lvec(0 1)\lvec(0 0)
\htext(-0.5 0.5){$j$}
\end{texdraw}}%
\savebox{\tmpfign}{\begin{texdraw}
\fontsize{8}{8}\selectfont
\drawdim em
\setunitscale 1.7
\move(0 0)\lvec(-1 0)
\move(0 1)\lvec(-1 1)\move(0 0)\lvec(0 2)\lvec(-1 2)\lvec(-1 0)
\htext(-0.5 1.5){$j$}
\htext(-0.5 0.5){$j\!\!+\!\!1$}
\end{texdraw}}%
\savebox{\tmpfigo}{\begin{texdraw}
\fontsize{8}{8}\selectfont
\drawdim em
\setunitscale 1.7
\move(0 0)\lvec(-1 0)\lvec(-1 1)\lvec(0 1)\lvec(0 0)
\htext(-0.5 0.5){$\bar{j}$}
\end{texdraw}}%
\savebox{\tmpfigp}{\begin{texdraw}
\fontsize{8}{8}\selectfont
\drawdim em
\setunitscale 1.7
\move(0 0)\lvec(-1 0)\move(0 1)\lvec(-1 1)
\move(0 0)\lvec(0 2)\lvec(-1 2)\lvec(-1 0)
\htext(-0.5 0.5){$n\!\!-\!\!1$}
\htext(-0.5 1.5){$n$}
\end{texdraw}}%
\savebox{\tmpfigq}{\begin{texdraw}
\fontsize{8}{8}\selectfont
\drawdim em
\setunitscale 1.7
\move(0 0)\lvec(-1 0)\lvec(-1 1)\lvec(0 1)\lvec(0 0)
\htext(-0.5 0.5){$0$}
\end{texdraw}}%
\savebox{\tmpfigr}{\begin{texdraw}
\fontsize{8}{8}\selectfont
\drawdim em
\setunitscale 1.7
\move(0 0)\lvec(-1 0)\lvec(-1 1)\lvec(0 1)\lvec(0 0)
\move(-1 0)\lvec(-0.5 0.5)\move(-0.05 0.95)\lvec(0 1)
\htext(-0.5 0.75){$n\!\!-\!\!1$}
\end{texdraw}}%
\savebox{\tmpfigs}{\begin{texdraw}
\fontsize{8}{8}\selectfont
\drawdim em
\setunitscale 1.7
\move(0 0)\lvec(-1 0)\lvec(-1 1)\lvec(0 1)\lvec(0 0)
\move(-1 0)\lvec(-0.95 0.05)\move(-0.5 0.5)\lvec(0 1)
\htext(-0.5 0.25){$n\!\!-\!\!1$}
\end{texdraw}}%
\savebox{\tmpfigt}{\begin{texdraw}
\fontsize{8}{8}\selectfont
\drawdim em
\setunitscale 1.7
\move(0 0)\lvec(-1 0)\lvec(-1 1)\lvec(0 1)\lvec(0 0)
\htext(-0.5 0.5){$n$}
\end{texdraw}}%
\savebox{\tmpfigu}{\begin{texdraw}
\fontsize{8}{8}\selectfont
\drawdim em
\setunitscale 1.7
\move(-1 0)\lvec(-1 1)\lvec(0 1)\lvec(-1 0)\lvec(0 0)\lvec(0 1)
\htext(-0.7 0.75){$n$}
\end{texdraw}}%
\savebox{\tmpfigv}{\begin{texdraw}
\fontsize{8}{8}\selectfont
\drawdim em
\setunitscale 1.7
\move(-1 0)\lvec(-1 1)\lvec(0 1)\lvec(-1 0)\lvec(0 0)\lvec(0 1)
\htext(-0.3 0.3){$n$}
\end{texdraw}}%
\savebox{\tmpfigw}{\begin{texdraw}
\fontsize{8}{8}\selectfont
\drawdim em
\setunitscale 1.7
\move(0 0)\lvec(-1 0)\lvec(-1 1)\lvec(0 1)\lvec(0 0)
\htext(-0.5 0.5){$\bar{n}$}
\end{texdraw}}%
\savebox{\tmpfigx}{\begin{texdraw}
\fontsize{8}{8}\selectfont
\drawdim em
\setunitscale 1.7
\move(0 0)\lvec(0 1)\lvec(-1 1)\lvec(-1 0)
\move(0 0)\lvec(-1 0)\lvec(0 1)
\htext(-0.3 0.3){$1$}
\htext(-0.7 0.75){$0$}
\end{texdraw}}%
\savebox{\tmpfigy}{\begin{texdraw}
\fontsize{8}{8}\selectfont
\drawdim em
\setunitscale 1.7
\move(0 0)\lvec(0 1)\lvec(-1 1)\lvec(-1 0)
\move(0 0)\lvec(-1 0)\lvec(0 1)
\htext(-0.3 0.3){$0$}
\htext(-0.7 0.75){$1$}
\end{texdraw}}%
\savebox{\tmpfigz}{\begin{texdraw}
\fontsize{8}{8}\selectfont
\drawdim em
\setunitscale 1.7
\lpatt(0.1 0.15)
\move(-1 1)\lvec(-1 2)\lvec(0 2)\lvec(0 1)
\lpatt(1 0)
\move(0 0)\lvec(-1 0)\lvec(-1 1)\lvec(0 1)\lvec(0 0)
\htext(-0.5 0.5){$n\!\!-\!\!2$}
\htext(-0.5 1.5){$n\!\!-\!\!1$}
\end{texdraw}}%
\savebox{\tmpfigaa}{\begin{texdraw}
\fontsize{8}{8}\selectfont
\drawdim em
\setunitscale 1.7
\move(0 0)\lvec(-1 0)\lvec(-1 1)\lvec(0 1)\lvec(0 0)
\htext(-0.5 0.5){$n\!\!-\!\!1$}
\end{texdraw}}%
\savebox{\tmpfigad}{\begin{texdraw}
\fontsize{8}{8}\selectfont
\drawdim em
\setunitscale 1.7
\move(0 0)\lvec(-1 0)\lvec(-1 1)\lvec(0 1)\lvec(0 0)
\htext(-0.5 0.5){$\overline{n\!\!-\!\!1}$}
\end{texdraw}}%
\savebox{\tmpfigab}{\begin{texdraw}
\fontsize{8}{8}\selectfont
\drawdim em
\setunitscale 1.7
\move(0 0)\lvec(-1 0)
\move(0 1)\lvec(-1 1)\lvec(-0.5 1.5)
\move(-0.05 1.95)\lvec(0 2)
\move(0 0)\lvec(0 2)\lvec(-1 2)\lvec(-1 0)
\htext(-0.5 0.5){$n\!\!-\!\!2$}
\htext(-0.3 1.3){$n$}
\htext(-0.5 1.75){$n\!\!-\!\!1$}
\end{texdraw}}%
\savebox{\tmpfigac}{\begin{texdraw}
\fontsize{8}{8}\selectfont
\drawdim em
\setunitscale 1.7
\move(0 0)\lvec(-1 0)
\move(0 1)\lvec(-1 1)\lvec(-0.95 1.05)
\move(-0.5 1.5)\lvec(0 2)
\move(0 0)\lvec(0 2)\lvec(-1 2)\lvec(-1 0)
\htext(-0.5 0.5){$n\!\!-\!\!2$}
\htext(-0.7 1.75){$n$}
\htext(-0.5 1.25){$n\!\!-\!\!1$}
\end{texdraw}}%
\savebox{\tmpfigag}{\begin{texdraw}
\fontsize{8}{8}\selectfont
\drawdim em
\setunitscale 1.7
\move(0 0)\lvec(-1 0)\lvec(-1 1)\lvec(0 1)\lvec(0 0)
\htext(-0.5 0.5){$\overline{n\!\!-\!\!2}$}
\end{texdraw}}%
\savebox{\tmpfigae}{\begin{texdraw}
\fontsize{8}{8}\selectfont
\drawdim em
\setunitscale 1.7
\move(0 0)\lvec(-1 0)\lvec(-0.5 0.5)
\move(-0.05 0.95)\lvec(0 1)
\move(0 1)\lvec(-1 1)
\move(0 0)\lvec(0 2)\lvec(-1 2)\lvec(-1 0)
\htext(-0.5 1.5){$n\!\!-\!\!2$}
\htext(-0.5 0.75){$n\!\!-\!\!1$}
\htext(-0.3 0.3){$n$}
\end{texdraw}}%
\savebox{\tmpfigaf}{\begin{texdraw}
\fontsize{8}{8}\selectfont
\drawdim em
\setunitscale 1.7
\move(0 0)\lvec(-1 0)\lvec(-0.95 0.05)
\move(-0.5 0.5)\lvec(0 1)
\move(0 1)\lvec(-1 1)
\move(0 0)\lvec(0 2)\lvec(-1 2)\lvec(-1 0)
\htext(-0.5 1.5){$n\!\!-\!\!2$}
\htext(-0.5 0.25){$n\!\!-\!\!1$}
\htext(-0.7 0.75){$n$}
\end{texdraw}}%
\savebox{\tmpfigah}{\begin{texdraw}
\fontsize{8}{8}\selectfont
\drawdim em
\setunitscale 1.7
\move(0 0)\lvec(0 0.5)\lvec(-1 0.5)\lvec(-1 0)\lvec(0 0)
\htext(-0.5 0.25){$0$}
\end{texdraw}}%
\savebox{\tmpfigai}{\begin{texdraw}
\fontsize{8}{8}\selectfont
\drawdim em
\setunitscale 1.7
\move(0 0)\lvec(0 1)\lvec(-1 1)\lvec(-1 0)\lvec(0 0)
\move(0 0.5)\lvec(-1 0.5)
\htext(-0.5 0.25){$0$}
\htext(-0.5 0.75){$0$}
\end{texdraw}}%
\savebox{\tmpfigaj}{\begin{texdraw}
\fontsize{8}{8}\selectfont
\drawdim em
\setunitscale 1.7
\move(0 0)\lvec(0 2)\lvec(-1 2)\lvec(-1 0)
\move(0 0)\lvec(-1 0)\move(0 1)\lvec(-1 1)
\htext(-0.5 0.5){$3$}
\htext(-0.5 1.5){$2$}
\end{texdraw}}%
\savebox{\tmpfigak}{\begin{texdraw}
\fontsize{8}{8}\selectfont
\drawdim em
\setunitscale 1.7
\move(0 0)\lvec(0 2)\lvec(-1 2)\lvec(-1 0)
\move(0 0)\lvec(-1 0)\move(0 1)\lvec(-1 1)
\htext(-0.5 0.5){$n\!\!-\!\!2$}
\htext(-0.5 1.5){$n\!\!-\!\!1$}
\end{texdraw}}%
\savebox{\tmpfigal}{\begin{texdraw}
\fontsize{8}{8}\selectfont
\drawdim em
\setunitscale 1.7
\move(0 0)\lvec(0 0.5)\lvec(-1 0.5)\lvec(-1 0)\lvec(0 0)
\htext(-0.5 0.26){$n$}
\end{texdraw}}%
\savebox{\tmpfigam}{\begin{texdraw}
\fontsize{8}{8}\selectfont
\drawdim em
\setunitscale 1.7
\move(0 0)\lvec(0 1)\lvec(-1 1)\lvec(-1 0)\lvec(0 0)
\move(0 0.5)\lvec(-1 0.5)
\htext(-0.5 0.26){$n$}
\htext(-0.5 0.76){$n$}
\end{texdraw}}%
\savebox{\tmpfigan}{\begin{texdraw}
\fontsize{8}{8}\selectfont
\drawdim em
\setunitscale 1.7
\move(0 0)\lvec(-1 0)\lvec(-1 1)\lvec(0 1)\lvec(0 0)
\htext(-0.5 0.5){$\bar{3}$}
\end{texdraw}}%
\savebox{\tmpfigao}{\begin{texdraw}
\fontsize{8}{8}\selectfont
\drawdim em
\setunitscale 1.7
\move(0 0)\lvec(-1 0)\lvec(-1 1)\lvec(0 1)\lvec(0 0)
\htext(-0.5 0.5){$j\!\!-\!\!1$}
\end{texdraw}}%
\savebox{\tmpfigap}{\begin{texdraw}
\fontsize{8}{8}\selectfont
\drawdim em
\setunitscale 1.7
\move(0 0)\lvec(0 1)\lvec(-1 1)\lvec(-1 0)
\move(0 0)\lvec(-1 0)\lvec(0 1)
\move(-1 1)\lvec(-1 2)\lvec(0 2)\lvec(0 1)
\htext(-0.3 0.3){$1$}
\htext(-0.7 0.75){$0$}
\htext(-0.5 1.5){$2$}
\end{texdraw}}%
\savebox{\tmpfigaq}{\begin{texdraw}
\fontsize{8}{8}\selectfont
\drawdim em
\setunitscale 1.7
\move(0 0)\lvec(0 1)\lvec(-1 1)\lvec(-1 0)
\move(0 0)\lvec(-1 0)\lvec(0 1)
\move(-1 1)\lvec(-1 2)\lvec(0 2)\lvec(0 1)
\htext(-0.3 0.3){$0$}
\htext(-0.7 0.75){$1$}
\htext(-0.5 1.5){$2$}
\end{texdraw}}%
\savebox{\tmpfigar}{\begin{texdraw}
\fontsize{8}{8}\selectfont
\drawdim em
\setunitscale 1.7
\move(0 0)\lvec(-1 0)\lvec(-1 1)\lvec(0 1)\lvec(0 0)
\htext(-0.5 0.5){$3$}
\end{texdraw}}%
\savebox{\tmpfigas}{\begin{texdraw}
\fontsize{8}{8}\selectfont
\drawdim em
\setunitscale 1.7
\move(0 0)\lvec(-1 0)\lvec(-1 1)\lvec(0 1)\lvec(0 0)
\htext(-0.5 0.5){$\bar{n}$}
\end{texdraw}}%
\savebox{\tmpfigat}{\begin{texdraw}
\fontsize{8}{8}\selectfont
\drawdim em
\setunitscale 1.7
\move(0 0)\lvec(-1 0)\lvec(-1 1)\lvec(0 1)\lvec(0 0)
\move(-1 0)\lvec(-0.5 0.5)\move(-0.05 0.95)\lvec(0 1)
\htext(-0.5 0.75){$n\!\!-\!\!1$}
\htext(-0.3 0.3){$n$}
\end{texdraw}}%
\savebox{\tmpfigau}{\begin{texdraw}
\fontsize{8}{8}\selectfont
\drawdim em
\setunitscale 1.7
\move(0 0)\lvec(-1 0)\lvec(-1 1)\lvec(0 1)\lvec(0 0)
\move(-1 0)\lvec(-0.95 0.05)\move(-0.5 0.5)\lvec(0 1)
\htext(-0.5 0.25){$n\!\!-\!\!1$}
\htext(-0.7 0.75){$n$}
\end{texdraw}}%
\savebox{\tmpfigav}{\begin{texdraw}
\fontsize{8}{8}\selectfont
\drawdim em
\setunitscale 1.7
\move(0 0)\lvec(-1 0)\lvec(-1 2)\lvec(0 2)\lvec(0 0)
\move(-1 0.5)\lvec(0 0.5)
\move(-1 1)\lvec(0 1)
\htext(-0.5 0.26){$0$}
\htext(-0.5 0.76){$0$}
\htext(-0.5 1.5){$1$}
\end{texdraw}}%
\savebox{\tmpfigaw}{\begin{texdraw}
\fontsize{8}{8}\selectfont
\drawdim em
\setunitscale 1.7
\move(0 0)\lvec(-1 0)\lvec(-1 2)\lvec(0 2)\lvec(0 0)
\move(-1 0.5)\lvec(0 0.5)
\move(-1 1)\lvec(0 1)
\htext(-0.5 0.26){$n$}
\htext(-0.5 0.76){$n$}
\htext(-0.5 1.5){$n\!\!-\!\!1$}
\end{texdraw}}%
\savebox{\tmpfigax}{\begin{texdraw}
\fontsize{8}{8}\selectfont
\drawdim em
\setunitscale 1.7
\move(0 0)\lvec(-1 0)\lvec(-1 1)\lvec(0 1)\lvec(0 0)
\htext(-0.5 0.5){$i$}
\end{texdraw}}%
\savebox{\tmpfigay}{\begin{texdraw}
\fontsize{8}{8}\selectfont
\drawdim em
\setunitscale 1.7
\move(0 0)\lvec(-1 0)\lvec(-1 1)\lvec(0 1)\lvec(0 0)
\htext(-0.5 0.5){$i\!\!+\!\!1$}
\end{texdraw}}%
\savebox{\tmpfigaz}{\begin{texdraw}
\fontsize{8}{8}\selectfont
\drawdim em
\setunitscale 1.7
\move(0 0)\lvec(-1 0)\lvec(-1 1)\lvec(0 1)\lvec(0 0)
\htext(-0.5 0.5){$\overline{i\!\!+\!\!1}$}
\end{texdraw}}%
\savebox{\tmpfigba}{\begin{texdraw}
\fontsize{8}{8}\selectfont
\drawdim em
\setunitscale 1.7
\move(0 0)\lvec(-1 0)\lvec(-1 1)\lvec(0 1)\lvec(0 0)
\htext(-0.5 0.5){$\bar{i}$}
\end{texdraw}}%
\savebox{\tmpfigbb}{\begin{texdraw}
\fontsize{8}{8}\selectfont
\drawdim em
\setunitscale 1.7
\lpatt(0.1 0.15)
\move(-1 1)\lvec(-1 2)\lvec(0 2)\lvec(0 1)
\lpatt(1 0)
\move(0 0)\lvec(-1 0)\lvec(-1 1)\lvec(0 1)\lvec(0 0)
\htext(-0.5 1.5){$i$}
\htext(-0.5 0.5){$i\!\!-\!\!1$}
\end{texdraw}}%
\savebox{\tmpfigbc}{\begin{texdraw}
\fontsize{8}{8}\selectfont
\drawdim em
\setunitscale 1.7
\move(0 0)\lvec(0 2)\lvec(-1 2)\lvec(-1 0)
\move(0 0)\lvec(-1 0)\move(0 1)\lvec(-1 1)
\htext(-0.5 1.5){$i$}
\htext(-0.5 0.5){$i\!\!-\!\!1$}
\end{texdraw}}%
\savebox{\tmpfigbd}{\begin{texdraw}
\fontsize{8}{8}\selectfont
\drawdim em
\setunitscale 1.7
\move(0 0)\lvec(0 2)\lvec(-1 2)\lvec(-1 0)
\move(0 0)\lvec(-1 0)\move(0 1)\lvec(-1 1)
\htext(-0.5 1.5){$i$}
\htext(-0.5 0.5){$i\!\!+\!\!1$}
\end{texdraw}}%
\savebox{\tmpfigbe}{\begin{texdraw}
\fontsize{8}{8}\selectfont
\drawdim em
\setunitscale 1.7
\lpatt(0.1 0.15)
\move(-1 1)\lvec(-1 2)\lvec(0 2)\lvec(0 1)
\lpatt(1 0)
\move(0 0)\lvec(-1 0)\lvec(-1 1)\lvec(0 1)\lvec(0 0)
\htext(-0.5 1.5){$i$}
\htext(-0.5 0.5){$i\!\!+\!\!1$}
\end{texdraw}}%
\savebox{\tmpfigbf}{\begin{texdraw}
\fontsize{8}{8}\selectfont
\drawdim em
\setunitscale 1.7
\move(0 0)\lvec(0 2)\lvec(-1 2)\lvec(-1 0)\lvec(0 0)
\move(0 1)\lvec(-1 1)\lvec(0 2)
\htext(-0.3 1.3){$0$}
\htext(-0.5 0.5){$2$}
\end{texdraw}}%
\savebox{\tmpfigbg}{\begin{texdraw}
\fontsize{8}{8}\selectfont
\drawdim em
\setunitscale 1.7
\move(0 0)\lvec(0 2)\lvec(-1 2)\lvec(-1 0)\lvec(0 0)
\move(0 1)\lvec(-1 1)\lvec(0 2)
\htext(-0.7 1.75){$0$}
\htext(-0.5 0.5){$2$}
\end{texdraw}}%
%
%FIGEND

\vskip 3mm 
\begin{multicols}{2}[\text{(a)  $A_n^{(1)}$ ($n\geq1$)}]

\hskip 5mm {\usebox{\tmpfigao}} \araise{$\longmapsto$} {\usebox{\tmpfigm}}
\quad \araise{($1\leq j\leq n$), }\\[4pt]

\vskip 2mm 
{\usebox{\tmpfigt}}  \araise{$\longmapsto$} {\usebox{\tmpfigq}}
\end{multicols}

\vskip 3mm 

\begin{multicols}{2}[\text{(b)  $A_{2n-1}^{(2)}$ ($n\geq3$)}]

\hskip 5mm 
\usebox{\tmpfigb} \araise{ or } \usebox{\tmpfiga}
\araise{$\longmapsto$} \usebox{\tmpfigc}\\[4pt]

\vskip 2mm 
\hskip 5mm 
\usebox{\tmpfigx} \araise{ or } \usebox{\tmpfigy}
\araise{$\longmapsto$} \usebox{\tmpfigi}\\[4pt]

\vskip 2mm
\hskip 5mm 
\usebox{\tmpfigap} \braise{ or } \usebox{\tmpfigaq}
\braise{$\longmapsto$} \craise{\usebox{\tmpfigar}}\\[4pt]

\vskip 2mm
\hskip 5mm 
\usebox{\tmpfigl}
\braise{$\longmapsto$} \craise{\usebox{\tmpfigm}}
\quad \draise{($4\leq j \leq n$)}\\[4pt]

\vskip 2mm 
\usebox{\tmpfigp}
\braise{$\longmapsto$} \craise{\usebox{\tmpfigw}}\\[4pt]

\vskip 4mm 
\usebox{\tmpfign}
\braise{$\longmapsto$} \craise{\usebox{\tmpfigo}}
\quad \draise{($2\leq j \leq n-1$)}\\[4pt]

\vskip 4mm 
\usebox{\tmpfigd} \araise{ or } \usebox{\tmpfige}
\araise{$\longmapsto$} \usebox{\tmpfigf}
\end{multicols}

\vskip 3mm 
\begin{multicols}{2}[\text{(c)  $D_n^{(1)}$ ($n\geq4$)}]

\hskip 5mm
\usebox{\tmpfigb} \araise{ or } \usebox{\tmpfiga}
\araise{$\longmapsto$} \usebox{\tmpfigc}\\[4pt]

\vskip 2mm
\hskip 5mm 
\usebox{\tmpfigx} \araise{ or } \usebox{\tmpfigy}
\araise{$\longmapsto$} \usebox{\tmpfigi}\\[4pt]

\vskip 2mm
\hskip 5mm 
\usebox{\tmpfigap} \braise{ or } \usebox{\tmpfigaq}
\braise{$\longmapsto$} \craise{\usebox{\tmpfigar}}\\[4pt]

\vskip 2mm
\hskip 5mm 
\usebox{\tmpfigl}
\braise{$\longmapsto$} \craise{\usebox{\tmpfigm}}
\quad \draise{($4\leq j \leq n-1$)}\\[4pt]

\vskip 2mm
\hskip 5mm
\usebox{\tmpfigs} \araise{ or } \usebox{\tmpfigr}
\araise{$\longmapsto$} \usebox{\tmpfigt}\\[4pt]

\vskip 2mm
\usebox{\tmpfigv} \araise{ or } \usebox{\tmpfigu}
\araise{$\longmapsto$} \usebox{\tmpfigas}\\[4pt]

\vskip 2mm
\usebox{\tmpfigat} \araise{ or } \usebox{\tmpfigau}
\araise{$\longmapsto$} \usebox{\tmpfigad}\\[4pt]

\vskip 2mm
\usebox{\tmpfigae} \braise{ or } \usebox{\tmpfigaf}
\braise{$\longmapsto$} \craise{\usebox{\tmpfigag}}\\[4pt]

\vskip 2mm 
\usebox{\tmpfign}
\braise{$\longmapsto$} \craise{\usebox{\tmpfigo}}
\quad \draise{($2\leq j \leq n-3$)}\\[4pt]

\vskip 2mm 
\usebox{\tmpfigd} \araise{ or } \usebox{\tmpfige}
\araise{$\longmapsto$} \usebox{\tmpfigf}
\end{multicols}

\vskip 3mm 

\begin{multicols}{2}[\text{(d)  $A_{2n}^{(2)}$ ($n\geq2$)}]

\hskip 5mm 
\usebox{\tmpfigah}
\araise{$\longmapsto$} \raisebox{3.6pt}{$\emptyset$}\\[4pt]

\vskip 2mm
\hskip 5mm 
\usebox{\tmpfigai}
\araise{$\longmapsto$} \usebox{\tmpfigc}\\[4pt]

\vskip 2mm
\hskip 5mm 
\usebox{\tmpfigav}
\braise{$\longmapsto$} \craise{\usebox{\tmpfigi}}\\[4pt]

%\vskip 2mm 
\usebox{\tmpfigl}
\braise{$\longmapsto$} \craise{\usebox{\tmpfigm}}
\quad \draise{($3\leq j \leq n$)}\\[4pt]

%\vskip 2mm 
\usebox{\tmpfigat}
\araise{$\longmapsto$} \usebox{\tmpfigw}\\[4pt]

%\vskip 2mm
\usebox{\tmpfign}
\braise{$\longmapsto$} \craise{\usebox{\tmpfigo}}
\quad \draise{($1\leq j \leq n-1$)}
\end{multicols}

\vskip 3mm
\begin{multicols}{2}[\text{(e) $D_{n+1}^{(2)}$ ($n\geq2$)}]

\hskip 5mm 
\usebox{\tmpfigah}
\araise{$\longmapsto$} \raisebox{3.6pt}{$\emptyset$}\\[4pt]

%\vskip 2mm 
\hskip 5mm 
\usebox{\tmpfigai}
\araise{$\longmapsto$} \usebox{\tmpfigc}\\[4pt]

%\vskip 2mm 
\hskip 5mm 
\usebox{\tmpfigav}
\braise{$\longmapsto$} \craise{\usebox{\tmpfigi}}\\[4pt]

%\vskip 2mm
\hskip 5mm 
\usebox{\tmpfigl}
\braise{$\longmapsto$} \craise{\usebox{\tmpfigm}}
\quad \draise{($3\leq j \leq n$)}\\[4pt]

%\vskip 2mm 
\usebox{\tmpfigal}
\araise{$\longmapsto$} \usebox{\tmpfigq}\\[4pt]

%\vskip 2mm 
\usebox{\tmpfigam}
\araise{$\longmapsto$} \usebox{\tmpfigw}\\[4pt]

%\vskip 2mm 
\usebox{\tmpfigaw}
\braise{$\longmapsto$} \craise{\usebox{\tmpfigad}}\\[4pt]

%\vskip 2mm 
\usebox{\tmpfign}
\braise{$\longmapsto$} \craise{\usebox{\tmpfigo}}
\quad \draise{($1\leq j \leq n-2$)}
\end{multicols}

\vskip 3mm 
\begin{multicols}{2}[\text{(f) $B_n^{(1)}$ ($n\geq3$)}]

\hskip 5mm 
\usebox{\tmpfiga} \araise{ or } \usebox{\tmpfigb}
\araise{$\longmapsto$} \usebox{\tmpfigc}\\[4pt]

%\vskip 2mm
\hskip 5mm 
\usebox{\tmpfigx} \araise{ or } \usebox{\tmpfigy}
\araise{$\longmapsto$} \usebox{\tmpfigi}\\[4pt]

%\vskip 2mm
\hskip 5mm
\usebox{\tmpfigap} \braise{ or } \usebox{\tmpfigaq}
\braise{$\longmapsto$} \craise{\usebox{\tmpfigar}}\\[4pt]

%\vskip 2mm
\hskip 5mm 
\usebox{\tmpfigl}
\braise{$\longmapsto$} \craise{\usebox{\tmpfigm}}
\quad \draise{($4\leq j \leq n$)}\\[4pt]

%\vskip 2mm
\hskip 5mm
\usebox{\tmpfigal}
\araise{$\longmapsto$} \usebox{\tmpfigq}\\[4pt]

%\vskip 2mm 
\usebox{\tmpfigam}
\araise{$\longmapsto$} \usebox{\tmpfigw}\\[4pt]

%\vskip 2mm
\usebox{\tmpfigaw}
\braise{$\longmapsto$} \craise{\usebox{\tmpfigad}}\\[4pt]

%\vskip 2mm 
\usebox{\tmpfign}
\braise{$\longmapsto$} \craise{\usebox{\tmpfigo}}
\quad \draise{($2\leq j \leq n-2$)}\\[4pt]

%\vskip 2mm 
\usebox{\tmpfigd} \araise{ or } \usebox{\tmpfige}
\araise{$\longmapsto$} \usebox{\tmpfigf}
\end{multicols}

\vskip 3mm

It is clear that this map sends the \emph{ground-state walls} to the
appropriate \emph{ground-state paths} and that the image does indeed
lie in the set $\bP(\la)$ of $\la$-paths.
Moreover, the surjectivity of $\Psi$ follows immediately from the
definition, and the injectivity of $\Psi$ follows from the fact that
the proper Young walls in $\Y(\la)$ are {\it reduced}.
Hence we have only to show that the map 
$\Psi$ commutes with Kashiwara operators.
We give a proof of our claim only for the case of $B_n^{(1)}$.
Other cases may be proved in similar manners and are less complicated.

\vskip 3mm 
Let $Y = (y_k)_{k=0}^{\infty} \in \Y(\la)$, where $y_k$ denotes the $k$-th
column of $Y$ (counting from the right to the left)
and let $\text{\bf p} = (\text{\bf p}(k))_{k=0}^{\infty} = \Psi(Y)$
be the image of $Y$ under $\Psi$.
Recall that the action of the Kashiwara operators is determined by the
{\it $i$-signatures} of $Y$ and $\text{\bf p}$.
Fix an index $i\in I$ and consider the Kashiwara operator $\fit$.
%We divide the proof of our claim into the following separate cases. 

\vskip 3mm 
\noindent 
(1) $i=0$: \  
Suppose that the top of the $k$-th column $y_k$ of $Y$ is the cube
\raisebox{-0.3\height}{\usebox{\tmpfige}}\,.
The column $y_k$ is certainly not 0-removable.
If it is 0-admissible, we would assign a $+$ to the column, which is what
we would also do with the corresponding element
$\text{\bf p}(k) = \raisebox{-0.3\height}{\usebox{\tmpfigf}}$\,.
If it is not 0-admissible, then the top of the column $y_{k-1}$ must be
either
\raisebox{-0.3\height}{\usebox{\tmpfigb}}
or
\raisebox{-0.3\height}{\usebox{\tmpfigy}}\,.
In each of the two cases, $y_{k-1}$ is neither 0-admissible nor 0-removable,
so we would assign nothing to $y_k$ and $y_{k-1}$.

\vskip 3mm 
Now, consider the corresponding terms $\text{\bf p}(k)$ and 
$\text{\bf p}(k-1)$ of
the path $\text{\bf p} = (\text{\bf p}(k))_{k=0}^\infty = \Psi(Y)$.
We have
$\text{\bf p}(k) =
\raisebox{-0.3\height}{\usebox{\tmpfigf}}$
and
$\text{\bf p}(k-1) =
\raisebox{-0.3\height}{\usebox{\tmpfigc}}
\text{ or }
\raisebox{-0.3\height}{\usebox{\tmpfigi}}$\,.
We assign  $+$ for
$\text{\bf p}(k) =
\raisebox{-0.3\height}{\usebox{\tmpfigf}}$
and $-$ for
$\text{\bf p}(k-1) =
\raisebox{-0.3\height}{\usebox{\tmpfigc}}
\text{ or }
\raisebox{-0.3\height}{\usebox{\tmpfigi}}$\,,
which cancel out to give nothing.
Thus the resulting $0$-signatures coincide with each other. 

\vskip 3mm 
Similarly, we can verify that we would assign  the same $0$-signatures
to $y_k$ and $\text{\bf p}(k)$ for the other cases of $y_k$.
In Table~\ref{sidezero}, we list all the possible nontrivial
cases for $i=0$.

\vskip 3mm

\noindent 
(2) $i=1$: \ This case is quite similar to the $i=0$ case.

\vskip 3mm

\noindent 
(3) $i=2$: \ We list all the possible nontrivial cases in Table~\ref{sidetwo}.

\vskip 3mm
\noindent 
(4) $3\leq i \leq n-2$: \ 
We list all the possible nontrivial cases in Table~\ref{sidegen}.

\vskip 3mm
\noindent 
(5) $i=n-1$: \ 
We list all the possible nontrivial cases in Table~\ref{sidenmo}.

\vskip 3mm
\noindent
(6) $i=n$: \ Observe that the top parts of the columns of $Y$ that have 
a nontrivial contribution to the $n$-signature of $Y$ have the 
following form. 

\vskip 3mm 
\begin{center}
\begin{texdraw}
\fontsize{11}{11}\selectfont
\drawdim em
\setunitscale 1.7
\move(-1 0)\lvec(-1 1.5)\lvec(0 1.5)\lvec(0 0)\lvec(-2 0)\lvec(-2 1)\lvec(0 1)
\move(2 0)\lvec(2 2)\lvec(3 2)\lvec(3 0)\lvec(1 0)\lvec(1 1.5)\lvec(3 1.5)
\move(1 1)\lvec(3 1)
\htext(0.53 0.5){$\cdots$}
\htext(3.53 0.5){$\cdots$}
\htext(-2.47 0.5){$\cdots$}
\htext(-1.5 -0.5){$\underbrace{\rule{1.5em}{0em}}_{A}$}
\htext(0.5 -0.5){$\underbrace{\rule{4.95em}{0em}}_{B}$}
\htext(2.5 -0.5){$\underbrace{\rule{1.5em}{0em}}_{C}$}
\htext(-1.5 0.5){$n\!\!-\!\!1$}
\htext(-0.5 0.5){$n\!\!-\!\!1$}
\htext(1.5 0.5){$n\!\!-\!\!1$}
\htext(2.5 0.5){$n\!\!-\!\!1$}
\htext(-0.5 1.26){$n$}
\htext(1.5 1.26){$n$}
\htext(2.5 1.26){$n$}
\htext(2.5 1.76){$n$}
\end{texdraw}
\end{center}

\vskip 3mm

\noindent
Here, we do not exclude the possibility of some of
the parts A, B, or C missing in this drawing.
The ground-state wall of weight $\La_n$ could also be considered as a
degenerate case of the above.
It suffices to verify that the $n$-signature for this segment and that
of the corresponding part in the path $\text{\bf p}$ are the same.

\vskip 3mm
In the table given below, we list all the possible combinations for 
the parts  A, B and C, and write down their $n$-signatures.

\begin{center}
\begin{tabular}{ccccccccc}
\multicolumn{3}{c}{presence of}&\quad&
\multicolumn{3}{c}{signature}&\quad&
total\\
A&B&C&&A&B&C&&signature\\
\hline\hline
yes&yes&yes&&$+$&&$-$&&\\
yes&yes&no&&+&+&&&++\\
yes&no&yes&&$+$&&$-$&&\\
yes&no&no&&$++$&&&&$++$\\
no&yes&yes&&&$-$&$-$&&$--$\\
no&yes&no&&&$-+$&&&$-+$\\
no&no&yes&&&&$--$&&$--$\\
\end{tabular}
\end{center}

\noindent
Now, one can easily verify that the $n$-signature for the above
segment of $Y$ coincides with that of the corresponding part in
the path $\text{\bf p}$.

\begin{sidewaystable}
\begin{tabular}{ccccccccccccccc}
\multicolumn{3}{c}{column}&\quad&
\multicolumn{3}{c}{signature}&\quad&
\multicolumn{3}{c}{path element}&\quad&
\multicolumn{3}{c}{signature}\\
$y_{k+1}$&$y_k$&$y_{k-1}$&&
$y_{k+1}$&$y_k$&$y_{k-1}$&&
$\text{\bf p}(k+1)$&$\text{\bf p}(k)$&$\text{\bf p}(k-1)$&&
$\text{\bf p}(k+1)$&$\text{\bf p}(k)$&$\text{\bf p}(k-1)$\\
\hline\hline&&&&&&&&&&&&&&\\[-7pt]
\usebox{\tmpfige}&\usebox{\tmpfigb}&&&
&&&&
\usebox{\tmpfigf}&\usebox{\tmpfigc}&&&
\araise{$+$}&\araise{$-$}&\\
\usebox{\tmpfigd}&\usebox{\tmpfiga}&&&
&&&&
\usebox{\tmpfigf}&\usebox{\tmpfigc}&&&
\araise{$+$}&\araise{$-$}&\\
\usebox{\tmpfigg}&\usebox{\tmpfigbf} \araise{or} \usebox{\tmpfigbg}&&&
&&&&
\usebox{\tmpfigk}&\usebox{\tmpfigc}&&&
\araise{$+$}&\araise{$-$}&\\
\araise{other}&\usebox{\tmpfigb} \araise{or} \usebox{\tmpfiga}&&&
&\araise{$-$}&&&
&\usebox{\tmpfigc}&&&
&\araise{$-$}&\\
\hline&&&&&&&&&&&&&&\\[-7pt]
\usebox{\tmpfigd}&\usebox{\tmpfigx}&&&
&&&&
\usebox{\tmpfigf}&\usebox{\tmpfigi}&&&
$+$&$-$&\\
\usebox{\tmpfige}&\usebox{\tmpfigy}&&&
&&&&
\usebox{\tmpfigf}&\usebox{\tmpfigi}&&&
$+$&$-$&\\
\araise{other}&\usebox{\tmpfigx} \araise{or} \usebox{\tmpfigy}&&&
&\araise{$-$}&&&
&\usebox{\tmpfigi}&&&
&\araise{$-$}&\\
\hline&&&&&&&&&&&&&&\\[-7pt]
&\usebox{\tmpfigg}&\usebox{\tmpfigbf} \araise{or} \usebox{\tmpfigbg}&&
&&&&
&\usebox{\tmpfigk}&\usebox{\tmpfigc}&&
&\araise{$+$}&\araise{$-$}\\
&\usebox{\tmpfigg}&\araise{other}&&
&\araise{$+$}&&&
&\usebox{\tmpfigk}&&&
&\araise{$+$}&\\
\hline&&&&&&&&&&&&&&\\[-7pt]
&\usebox{\tmpfige}&\usebox{\tmpfigb}&&
&&&&
&\usebox{\tmpfigf}&\usebox{\tmpfigc}&&
&\araise{$+$}&\araise{$-$}\\
&\usebox{\tmpfigd}&\usebox{\tmpfiga}&&
&&&&
&\usebox{\tmpfigf}&\usebox{\tmpfigc}&&
&\araise{$+$}&\araise{$-$}\\
&\usebox{\tmpfige}&\usebox{\tmpfigy}&&
&&&&
&\usebox{\tmpfigf}&\usebox{\tmpfigi}&&
&\araise{$+$}&\araise{$-$}\\
&\usebox{\tmpfigd}&\usebox{\tmpfigx}&&
&&&&
&\usebox{\tmpfigf}&\usebox{\tmpfigi}&&
&\araise{$+$}&\araise{$-$}\\
&\usebox{\tmpfige} \araise{or} \usebox{\tmpfigd}&\araise{other}&&
&\araise{$+$}&&&
&\usebox{\tmpfigf}&&&
&\araise{$+$}&\\
\hline&&&&&&&&&&&&&&\\[-7pt]
\end{tabular}
\caption{$i=0$}
\label{sidezero}
\end{sidewaystable}

\begin{sidewaystable}
\begin{tabular}{ccccccccccccccc}
\multicolumn{3}{c}{column}&\quad&
\multicolumn{3}{c}{signature}&\quad&
\multicolumn{3}{c}{path element}&\quad&
\multicolumn{3}{c}{signature}\\
$y_{k+1}$&$y_k$&$y_{k-1}$&&
$y_{k+1}$&$y_k$&$y_{k-1}$&&
$\text{\bf p}(k+1)$&$\text{\bf p}(k)$&$\text{\bf p}(k-1)$&&
$\text{\bf p}(k+1)$&$\text{\bf p}(k)$&$\text{\bf p}(k-1)$\\
\hline\hline&&&&&&&&&&&&&&\\[-7pt]
&\usebox{\tmpfigx} \araise{or} \usebox{\tmpfigy}&\usebox{\tmpfigh}&&
&&&&
&\usebox{\tmpfigi}&\usebox{\tmpfigar}&&
&\araise{$+$}&\araise{$-$}\\
&\usebox{\tmpfigx} \araise{or} \usebox{\tmpfigy}&\araise{other}&&
&\araise{$+$}&&&
&\usebox{\tmpfigi}&&&
&\araise{$+$}&\\
\hline&&&&&&&&&&&&&&\\[-7pt]
\usebox{\tmpfigx} \araise{or} \usebox{\tmpfigy}&\usebox{\tmpfigh}&&&
&&&&
\usebox{\tmpfigi}&\usebox{\tmpfigar}&&&
\araise{$+$}&\araise{$-$}&\\
\araise{other}&\usebox{\tmpfigh}&&&
&\araise{$-$}&&&
&\usebox{\tmpfigar}&&&
&\araise{$-$}&\\
\hline&&&&&&&&&&&&&&\\[-7pt]
&\usebox{\tmpfigj}&\usebox{\tmpfigaj}&&
&&&&
&\usebox{\tmpfigan}&\usebox{\tmpfigk}&&
&\araise{$+$}&\araise{$-$}\\
&\usebox{\tmpfigj}&\araise{other}&&
&\araise{$+$}&&&
&\usebox{\tmpfigan}&&&
&\araise{$+$}&\\
\hline&&&&&&&&&&&&&&\\[-7pt]
\usebox{\tmpfigj}&\usebox{\tmpfigaj}&&&
&&&&
\usebox{\tmpfigan}&\usebox{\tmpfigk}&&&
\araise{$+$}&\araise{$-$}&\\
\araise{other}&\usebox{\tmpfigaj}&&&
&\araise{$-$}&&&
&\usebox{\tmpfigk}&&&
&\araise{$-$}&\\
\hline&&&&&&&&&&&&&&\\[-7pt]
\end{tabular}
\caption{$i=2$}
\label{sidetwo}
\end{sidewaystable}

\vskip 3mm

\begin{sidewaystable}
\begin{tabular}{ccccccccccccccc}
\multicolumn{3}{c}{column}&\quad&
\multicolumn{3}{c}{signature}&\quad&
\multicolumn{3}{c}{path element}&\quad&
\multicolumn{3}{c}{signature}\\
$y_{k+1}$&$y_k$&$y_{k-1}$&&
$y_{k+1}$&$y_k$&$y_{k-1}$&&
$\text{\bf p}(k+1)$&$\text{\bf p}(k)$&$\text{\bf p}(k-1)$&&
$\tmppathp(k+1)$&$\tmppathp(k)$&$\tmppathp(k-1)$\\
\hline\hline&&&&&&&&&&&&&&\\[-7pt]
&\usebox{\tmpfigbb}&\usebox{\tmpfigbc}&&
&&&&
&\usebox{\tmpfigax}&\usebox{\tmpfigay}&&
&\araise{$+$}&\araise{$-$}\\
&\usebox{\tmpfigbb}&\araise{other}&&
&\araise{$+$}&&&
&\usebox{\tmpfigax}&&&
&\araise{$+$}&\\
\hline&&&&&&&&&&&&&&\\[-7pt]
\usebox{\tmpfigbb}&\usebox{\tmpfigbc}&&&
&&&&
\usebox{\tmpfigax}&\usebox{\tmpfigay}&&&
\araise{$+$}&\araise{$-$}&\\
\araise{other}&\usebox{\tmpfigbc}&&&
&\araise{$-$}&&&
&\usebox{\tmpfigay}&&&
&\araise{$-$}&\\
\hline&&&&&&&&&&&&&&\\[-7pt]
&\usebox{\tmpfigbe}&\usebox{\tmpfigbd}&&
&&&&
&\usebox{\tmpfigaz}&\usebox{\tmpfigba}&&
&\araise{$+$}&\araise{$-$}\\
&\usebox{\tmpfigbe}&\araise{other}&&
&\araise{$+$}&&&
&\usebox{\tmpfigba}&&&
&\araise{$+$}&\\
\hline&&&&&&&&&&&&&&\\[-7pt]
\usebox{\tmpfigbe}&\usebox{\tmpfigbd}&&&
&&&&
\usebox{\tmpfigaz}&\usebox{\tmpfigba}&&&
\araise{$+$}&\araise{$-$}&\\
\araise{other}&\usebox{\tmpfigbd}&&&
&\araise{$-$}&&&
&\usebox{\tmpfigba}&&&
&\araise{$-$}&\\
\hline&&&&&&&&&&&&&&\\[-7pt]
\end{tabular}
\caption{$3\leq i \leq n-2$}
\label{sidegen}
\end{sidewaystable}

\vskip 3mm

\begin{sidewaystable}
\centering
\begin{tabular}{ccccccccccccccc}
\multicolumn{3}{c}{column}&\quad&
\multicolumn{3}{c}{signature}&\quad&
\multicolumn{3}{c}{path element}&\quad&
\multicolumn{3}{c}{signature}\\
$y_{k+1}$&$y_k$&$y_{k-1}$&&
$y_{k+1}$&$y_k$&$y_{k-1}$&&
$\tmppathp(k+1)$&$\tmppathp(k)$&$\tmppathp(k-1)$&&
$\tmppathp(k+1)$&$\tmppathp(k)$&$\tmppathp(k-1)$\\
\hline\hline&&&&&&&&&&&&&&\\[-7pt]
&\usebox{\tmpfigz}&\usebox{\tmpfigak}&&
&&&&
&\usebox{\tmpfigaa}&\usebox{\tmpfigt}&&
&\araise{$+$}&\araise{$-$}\\
&\usebox{\tmpfigz}&\araise{other}&&
&\araise{$+$}&&&
&\usebox{\tmpfigaa}&&&
&\araise{$+$}&\\
\hline&&&&&&&&&&&&&&\\[-7pt]
\usebox{\tmpfigz}&\usebox{\tmpfigak}&&&
&&&&
\usebox{\tmpfigaa}&\usebox{\tmpfigt}&&&
\araise{$+$}&\araise{$-$}&\\
\araise{other}&\usebox{\tmpfigak}&&&
&\araise{$-$}&&&
&\usebox{\tmpfigt}&&&
&\araise{$-$}&\\
\hline&&&&&&&&&&&&&&\\[-7pt]
&\usebox{\tmpfigam}&\usebox{\tmpfigaw}&&
&&&&
&\usebox{\tmpfigw}&\usebox{\tmpfigad}&&
&\araise{$+$}&\araise{$-$}\\
&\usebox{\tmpfigam}&\araise{other}&&
&\araise{$+$}&&&
&\usebox{\tmpfigw}&&&
&\araise{$+$}&\\
\hline&&&&&&&&&&&&&&\\[-7pt]
\usebox{\tmpfigam}&\usebox{\tmpfigaw}&&&
&&&&
\usebox{\tmpfigw}&\usebox{\tmpfigad}&&&
\araise{$+$}&\araise{$-$}&\\
\araise{other}&\usebox{\tmpfigaw}&&&
&\araise{$-$}&&&
&\usebox{\tmpfigad}&&&
&\araise{$-$}&\\
\hline&&&&&&&&&&&&&&\\[-7pt]
\end{tabular}
\caption{$i=n-1$}
\label{sidenmo}
\end{sidewaystable}

\begin{center}
\begin{texdraw}
\small
\drawdim em
\setunitscale 1.4
\move(0 0)
\bsegment
\move(0 0)\lvec(-1 0)\lvec(-1 1)\lvec(0 1)\lvec(0 0)
\htext(-0.5 0.5){$n$}
\esegment
\move(2 0)
\bsegment
\move(0 0)\lvec(-1 0)\lvec(-1 1)\lvec(0 1)\lvec(0 0)
\htext(-0.5 0.5){$0$}
\esegment
\move(6 0)
\bsegment
\move(0 0)\lvec(-1 0)\lvec(-1 1)\lvec(0 1)\lvec(0 0)
\htext(-0.5 0.5){$0$}
\esegment
\move(8 0)
\bsegment
\move(0 0)\lvec(-1 0)\lvec(-1 1)\lvec(0 1)\lvec(0 0)
\htext(-0.5 0.5){$\bar{n}$}
\esegment
\htext(-2.5 0.5){$\cdots$}
\htext(-1.5 0.5){$\otimes$}
\htext(0.5 0.5){$\otimes$}
\htext(2.5 0.5){$\otimes$}
\htext(3.6 0.5){$\cdots$}
\htext(4.5 0.5){$\otimes$}
\htext(6.5 0.5){$\otimes$}
\htext(8.5 0.5){$\otimes$}
\htext(9.5 0.5){$\cdots$}
\htext(-0.5 -0.7){$++$}
\htext(1.5 -0.7){$-+$}
\htext(5.5 -0.7){$-+$}
\htext(7.5 -0.7){$--$}
\htext(-0.5 1.8){$\overbrace{\rule{2.5em}{0em}}^{A}$}
\htext(3.5 1.8){$\overbrace{\rule{8.1em}{0em}}^{B}$}
\htext(7.5 1.8){$\overbrace{\rule{2.5em}{0em}}^{C}$}
\end{texdraw}
\end{center}

Thus, we have shown that for all indices $i \in I$, the $i$-signature 
of the reduced proper Young wall $Y$ is identical to 
that of the path $\text{\bf p}=\Psi(Y)$.
It is now straightforward to verify that the action of Kashiwara operators
is compatible with the identification of $Y$ with $\text{\bf p}$. 
Here are some examples.

\vskip 3mm 
\begin{center}
\begin{texdraw}
\fontsize{10}{10}\selectfont
\drawdim em
\setunitscale 1.7
\move(-4 0)
\bsegment
\move(-1 0)
\bsegment
\move(-1 0)\lvec(-1 1)\lvec(0 1)\lvec(-1 0)\lvec(0 0)\lvec(0 1)
\htext(-0.3 0.3){$1$}
\esegment
\move(2 0)
\bsegment
\move(0 0)\lvec(0 1)\lvec(-1 1)\lvec(-1 0)
\move(0 0)\lvec(-1 0)\lvec(0 1)
\htext(-0.3 0.3){$1$}
\htext(-0.7 0.75){$0$}
\esegment
\htext(-0.1 0.9){$0$}
\move(-0.7 0.5)\avec(0.7 0.5)
\esegment
\move(4 0)
\bsegment
\move(-1 0)
\bsegment
\move(0 0)\lvec(-1 0)\lvec(-1 1)\lvec(0 1)\lvec(0 0)
\htext(-0.5 0.5){$\bar{1}$}
\esegment
\move(2 0)
\bsegment
\move(0 0)\lvec(-1 0)\lvec(-1 1)\lvec(0 1)\lvec(0 0)
\htext(-0.5 0.5){$2$}
\esegment
\htext(-0.1 0.9){$0$}
\move(-0.7 0.5)\avec(0.7 0.5)
\esegment
\htext(0 0.5){$\Longleftrightarrow$}
\end{texdraw}

\vskip 2mm

\begin{texdraw}
\fontsize{10}{10}\selectfont
\drawdim em
\setunitscale 1.7
\move(-4 0)
\bsegment
\move(-1 0)
\bsegment
\move(-1 0)\lvec(-1 1)\lvec(0 1)\lvec(-1 0)\lvec(0 0)\lvec(0 1)
\htext(-0.3 0.3){$1$}
\htext(-0.7 0.75){$0$}
\esegment
\move(2 0)
\bsegment
\move(-1 0)\lvec(-1 2)\lvec(0 2)\lvec(0 0)\lvec(-1 0)
\move(-1 0)\lvec(0 1)\lvec(-1 1)
\htext(-0.5 1.5){$2$}
\htext(-0.3 0.3){$1$}
\htext(-0.7 0.75){$0$}
\esegment
\htext(-0.1 0.9){$2$}
\move(-0.7 0.5)\avec(0.7 0.5)
\esegment
\move(4 0)
\bsegment
\move(-1 0)
\bsegment
\move(0 0)\lvec(-1 0)\lvec(-1 1)\lvec(0 1)\lvec(0 0)
\htext(-0.5 0.5){$2$}
\esegment
\move(2 0)
\bsegment
\move(0 0)\lvec(-1 0)\lvec(-1 1)\lvec(0 1)\lvec(0 0)
\htext(-0.5 0.5){$3$}
\esegment
\htext(-0.1 0.9){$2$}
\move(-0.7 0.5)\avec(0.7 0.5)
\esegment
\htext(0 0.5){$\Longleftrightarrow$}
\end{texdraw}

\vskip 2mm

\begin{texdraw}
\fontsize{10}{10}\selectfont
\drawdim em
\setunitscale 1.7
\move(-4 0)
\bsegment
\move(-1 0)
\bsegment
\move(0 0)\lvec(0 1)\lvec(-1 1)\lvec(-1 0)\lvec(0 0)
\htext(-0.5 0.5){$i\!\!-\!\!1$}
\esegment
\move(2 0)
\bsegment
\move(0 0)\lvec(0 2)\lvec(-1 2)\lvec(-1 0)
\move(0 0)\lvec(-1 0)\move(0 1)\lvec(-1 1)
\htext(-0.5 0.5){$i\!\!-\!\!1$}
\htext(-0.5 1.5){$i$}
\esegment
\htext(-0.1 1.1){$i$}
\move(-0.7 0.7)\avec(0.7 0.7)
\esegment
\move(4 0)
\bsegment
\move(-1 0)
\bsegment
\move(0 0)\lvec(-1 0)\lvec(-1 1)\lvec(0 1)\lvec(0 0)
\htext(-0.5 0.5){$i$}
\esegment
\move(2 0)
\bsegment
\move(0 0)\lvec(-1 0)\lvec(-1 1)\lvec(0 1)\lvec(0 0)
\htext(-0.5 0.5){$i\!\!+\!\!1$}
\esegment
\htext(-0.1 0.9){$i$}
\move(-0.7 0.5)\avec(0.7 0.5)
\esegment
\htext(0 0.5){$\Longleftrightarrow$}
\end{texdraw}

\vskip 2mm

\begin{texdraw}
\fontsize{10}{10}\selectfont
\drawdim em
\setunitscale 1.7
\move(-4 0)
\bsegment
\move(-1 0)
\bsegment
\move(0 0)\lvec(0 1)\lvec(-1 1)\lvec(-1 0)\lvec(0 0)
\htext(-0.5 0.5){$i\!\!+\!\!1$}
\esegment
\move(2 0)
\bsegment
\move(0 0)\lvec(0 2)\lvec(-1 2)\lvec(-1 0)
\move(0 0)\lvec(-1 0)\move(0 1)\lvec(-1 1)
\htext(-0.5 0.5){$i\!\!+\!\!1$}
\htext(-0.5 1.5){$i$}
\esegment
\htext(-0.1 1.1){$i$}
\move(-0.7 0.7)\avec(0.7 0.7)
\esegment
\move(4 0)
\bsegment
\move(-1 0)
\bsegment
\move(0 0)\lvec(-1 0)\lvec(-1 1)\lvec(0 1)\lvec(0 0)
\htext(-0.5 0.5){$\overline{i\!\!+\!\!1}$}
\esegment
\move(2 0)
\bsegment
\move(0 0)\lvec(-1 0)\lvec(-1 1)\lvec(0 1)\lvec(0 0)
\htext(-0.5 0.5){$\bar{i}$}
\esegment
\htext(-0.1 0.9){$i$}
\move(-0.7 0.5)\avec(0.7 0.5)
\esegment
\htext(0 0.5){$\Longleftrightarrow$}
\end{texdraw}

\vskip 2mm 

\begin{texdraw}
\fontsize{10}{10}\selectfont
\drawdim em
\setunitscale 1.7
\move(-4 0)
\bsegment
\move(-1 0)
\bsegment
\move(0 0)\lvec(-1 0)\lvec(-1 1)\lvec(0 1)\lvec(0 0)
\htext(-0.5 0.5){$n\!\!-\!\!1$}
\esegment
\move(2 0)
\bsegment
\move(0 0)\lvec(-1 0)\lvec(-1 1.5)\lvec(0 1.5)\lvec(0 0)
\move(0 1)\lvec(-1 1)
\htext(-0.5 0.5){$n\!\!-\!\!1$}
\htext(-0.5 1.26){$n$}
\esegment
\htext(-0.1 1.1){$n$}
\move(-0.7 0.7)\avec(0.7 0.7)
\esegment
\move(4 0)
\bsegment
\move(-1 0)
\bsegment
\move(0 0)\lvec(-1 0)\lvec(-1 1)\lvec(0 1)\lvec(0 0)
\htext(-0.5 0.5){$n$}
\esegment
\move(2 0)
\bsegment
\move(0 0)\lvec(-1 0)\lvec(-1 1)\lvec(0 1)\lvec(0 0)
\htext(-0.5 0.5){$0$}
\esegment
\htext(-0.1 0.9){$n$}
\move(-0.7 0.5)\avec(0.7 0.5)
\esegment
\htext(0 0.5){$\Longleftrightarrow$}
\end{texdraw}

\vskip 2mm
\begin{texdraw}
\fontsize{10}{10}\selectfont
\drawdim em
\setunitscale 1.7
\move(-4 0)
\bsegment
\move(-1 0)
\bsegment
\move(0 0)\lvec(-1 0)\lvec(-1 0.5)\lvec(0 0.5)\lvec(0 0)
\htext(-0.5 0.26){$n$}
\esegment
\move(2 0)
\bsegment
\move(0 0)\lvec(-1 0)\lvec(-1 1)\lvec(0 1)\lvec(0 0)
\move(0 0.5)\lvec(-1 0.5)
\htext(-0.5 0.26){$n$}
\htext(-0.5 0.76){$n$}
\esegment
\htext(-0.1 0.9){$n$}
\move(-0.7 0.5)\avec(0.7 0.5)
\esegment
\move(4 0)
\bsegment
\move(-1 0)
\bsegment
\move(0 0)\lvec(-1 0)\lvec(-1 1)\lvec(0 1)\lvec(0 0)
\htext(-0.5 0.5){$0$}
\esegment
\move(2 0)
\bsegment
\move(0 0)\lvec(-1 0)\lvec(-1 1)\lvec(0 1)\lvec(0 0)
\htext(-0.5 0.5){$\bar{n}$}
\esegment
\htext(-0.1 0.9){$n$}
\move(-0.7 0.5)\avec(0.7 0.5)
\esegment
\htext(0 0.5){$\Longleftrightarrow$}
\end{texdraw}
\end{center}

%\vskip 3mm
\noindent
Therefore, we obtain the desired crystal isomorphism 
$\Psi: \Y(\la) \stackrel{\sim} \longrightarrow \bP(\la)$. 
\qed

}% This brace is placed here and at the beginning of this section to limit
%  the scope of \tmpfig's. Do not erase.
%

\vskip 5mm
\newpage
\begin{example}
In this example, we illustrate the top parts of the affine crystal
$\Y(\la)$ consisting of reduced proper Young walls. 
Compare them with the path realization of crystal graphs 
given in Section 4.

\vskip 3mm

(a) The crystal $\Y(\La_0)$ for $A_2^{(1)}$

\vskip 3mm 
\begin{center}
\begin{texdraw}
\small
\drawdim em
\setunitscale 1.2
%\htext(-9 4){\normalsize$\bullet$ $\Y(\La_0)$ for $A_2^{(1)}$}
\htext(-1.3 5){$Y_{\La_0}$}
\move(-1.5 4.3)\avec(-1.5 2.5)
\htext(-1.1 3.5){$0$}
\move(-1 1)
\bsegment
\move(0 0)\lvec(-1 0)\lvec(-1 1)\lvec(0 1)\lvec(0 0)
\htext(-0.5 0.5){$0$}
\esegment
\move(-5.5 -4)
\bsegment
\move(0 0)\lvec(-1 0)\lvec(-1 2)\lvec(0 2)\lvec(0 0)
\move(0 1)\lvec(-1 1)
\htext(-0.5 1.5){$1$} \htext(-0.5 0.5){$0$}
\esegment
\move(4 -3.5)
\bsegment
\move(0 0)\lvec(-2 0)\lvec(-2 1)\lvec(0 1)\lvec(0 0)
\move(-1 0)\lvec(-1 1)
\htext(-0.5 0.5){$0$} \htext(-1.5 0.5){$2$}
\esegment
\move(-2.7 -9)
\bsegment
\move(-1 0)\lvec(-1 2)\lvec(0 2)\lvec(0 0)\lvec(-2 0)\lvec(-2 1)\lvec(0 1)
\htext(-0.5 0.5){$0$} \htext(-0.5 1.5){$1$} \htext(-1.5 0.5){$2$}
\esegment
\move(2.5 -9)
\bsegment
\move(0 0)\lvec(-3 0)\lvec(-3 1)\lvec(0 1)\lvec(0 0)
\move(-1 0)\lvec(-1 1) \move(-2 0)\lvec(-2 1)
\htext(-0.5 0.5){$0$} \htext(-1.5 0.5){$2$} \htext(-2.5 0.5){$1$}
\esegment
\move(-8 -14)
\bsegment
\move(-1 0)\lvec(-1 2)\lvec(0 2)\lvec(0 0)\lvec(-3 0)\lvec(-3 1)\lvec(0 1)
\move(-2 0)\lvec(-2 1)
\htext(-0.5 0.5){$0$} \htext(-0.5 1.5){$1$} \htext(-1.5 0.5){$2$}
\htext(-2.5 0.5){$1$}
\esegment
\move(-3 -14)
\bsegment
\move(0 0)\lvec(-2 0)\lvec(-2 2)\lvec(0 2)\lvec(0 0)
\move(0 1)\lvec(-2 1)\move(-1 0)\lvec(-1 2)
\htext(-0.5 0.5){$0$}\htext(-0.5 1.5){$1$}
\htext(-1.5 0.5){$2$}\htext(-1.5 1.5){$0$}
\esegment
\move(3 -14)
\bsegment
\move(0 0)\lvec(-4 0)\lvec(-4 1)\lvec(0 1)\lvec(0 0)
\move(-1 0)\lvec(-1 1)\move(-2 0)\lvec(-2 1)\move(-3 0)\lvec(-3 1)
\htext(-0.5 0.5){$0$}\htext(-1.5 0.5){$2$}\htext(-2.5 0.5){$1$}
\htext(-3.5 0.5){$0$}
\esegment
\move(8 -14)
\bsegment
\move(0 1)\lvec(-2 1)\lvec(-2 0)\lvec(0 0)\lvec(0 3)\lvec(-1 3)\lvec(-1 0)
\move(0 2)\lvec(-1 2)
\htext(-0.5 0.5){$0$}\htext(-1.5 0.5){$2$}
\htext(-0.5 1.5){$1$}\htext(-0.5 2.5){$2$}
\esegment
\move(-0.5 0.5)\avec(2.7 -1.7)\htext(1.4 -0.2){$2$}
\move(3 -4.5)\avec(1 -7.5)\htext(1.6 -5.7){$1$}
\move(1 -9.5)\avec(1 -11.5)\htext(1.4 -10){$0$}
\move(0 -9.5)\avec(-7.5 -11.8)\htext(-6 -10.8){$1$}
\move(-2.5 0.5)\avec(-5.8 -1.7)\htext(-4.5 -0.2){$1$}
\move(-6 -4.5)\avec(-4 -7.5)\htext(-4.6 -5.7){$2$}
\move(-4 -9.5)\avec(-4 -11.5)\htext(-3.6 -10){$0$}
\move(-3 -9.5)\avec(6 -11.8)\htext(4 -10.8){$2$}
\vtext(-9.5 -15){$\cdots$}
\vtext(-4 -15){$\cdots$}
\vtext(1 -15){$\cdots$}
\vtext(7 -15){$\cdots$}
\move(0 -15.7)\move(0 5.4)
\end{texdraw}
\end{center}

%%%%%%%%%%%%%%%%%%%%%

\newpage

(b) The crystal $\Y(\La_0)$ for $A_5^{(2)}$

\vskip 3mm

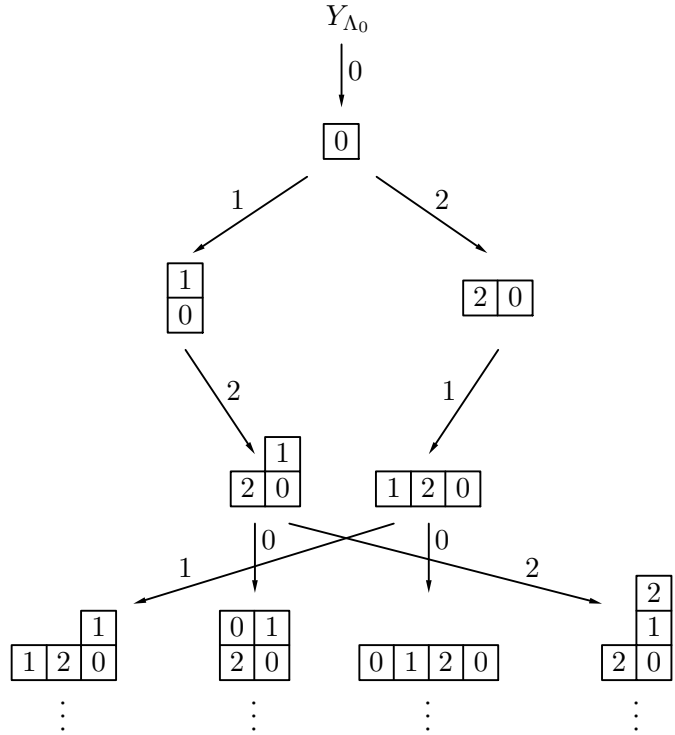
\begin{center}
\begin{texdraw}
\fontsize{7}{7}\selectfont
\drawdim em
\setunitscale 1.7
%\htext(-8 32){\normalsize $\bullet$ $\Y(\La_0)$ for $A_5^{(1)}$}
\nc{\dtri}{
\bsegment
\move(-1 0)\lvec(0 1)\lvec(0 0)\lvec(-1 0)\ifill f:0.7
\esegment
}
\move(-5 0)
\bsegment
\move(0 0)\dtri \move(-1 0)\dtri
\move(0 1)\lvec(-1 0)\lvec(-1 5)\lvec(0 5)\lvec(0 0)\lvec(-2 0)\lvec(-2 1)
\lvec(0 1)
\move(-2 0)\lvec(-1 1)
\move(0 2)\lvec(-1 2)\move(0 3)\lvec(-1 3)\move(0 4)\lvec(-1 4)\lvec(0 5)
\htext(-0.3 0.3){$1$}\htext(-0.7 0.75){$0$}\htext(-1.3 0.3){$0$}
\htext(-1.7 0.75){$1$}\htext(-0.5 1.5){$2$}\htext(-0.5 2.5){$3$}
\htext(-0.5 3.5){$2$}\htext(-0.7 4.75){$0$}
\esegment
\move(0 0)
\bsegment
\move(0 0)\dtri \move(-1 0)\dtri
\move(0 1)\lvec(-1 0)\lvec(-1 5)\lvec(0 5)\lvec(0 0)\lvec(-2 0)\lvec(-2 1)
\lvec(0 1)
\move(-2 0)\lvec(-1 1)
\move(0 2)\lvec(-1 2)\move(0 3)\lvec(-1 3)\move(0 4)\lvec(-1 4)\lvec(0 5)
\htext(-0.3 0.3){$1$}\htext(-0.7 0.75){$0$}\htext(-1.3 0.3){$0$}
\htext(-1.7 0.75){$1$}\htext(-0.5 1.5){$2$}\htext(-0.5 2.5){$3$}
\htext(-0.5 3.5){$2$}\htext(-0.3 4.3){$1$}
\esegment
\move(5 0)
\bsegment
\move(0 0)\dtri \move(-1 0)\dtri
\move(-2 1)\lvec(0 1)\lvec(-1 0)\lvec(-1 4)\lvec(0 4)\lvec(0 0)\lvec(-2 0)
\lvec(-2 2)\lvec(0 2)
\move(0 3)\lvec(-1 3)\move(-2 0)\lvec(-1 1)
\htext(-0.3 0.3){$1$}\htext(-0.7 0.75){$0$}\htext(-1.3 0.3){$0$}
\htext(-1.7 0.75){$1$}\htext(-0.5 1.5){$2$}\htext(-1.5 1.5){$2$}
\htext(-0.5 2.5){$3$}\htext(-0.5 3.5){$2$}
\esegment
\move(-5 8)
\bsegment
\move(0 0)\dtri
\move(-1 1)\lvec(0 1)\lvec(-1 0)\lvec(0 0)\lvec(0 5)\lvec(-1 5)\lvec(-1 0)
\move(0 2)\lvec(-1 2)\move(0 3)\lvec(-1 3)\move(0 4)\lvec(-1 4)\lvec(0 5)
\htext(-0.3 0.3){$1$}\htext(-0.7 0.75){$0$}
\htext(-0.5 1.5){$2$}\htext(-0.5 2.5){$3$}
\htext(-0.5 3.5){$2$}\htext(-0.7 4.75){$0$}
\esegment
\move(0 8)
\bsegment
\move(0 0)\dtri \move(-1 0)\dtri
\move(-1 1)\lvec(-2 0)\lvec(-2 1)\lvec(0 1)\lvec(-1 0)\lvec(-1 4)\lvec(0 4)
\lvec(0 0)\lvec(-2 0)\move(0 2)\lvec(-1 2)\move(0 3)\lvec(-1 3)
\htext(-0.3 0.3){$1$}\htext(-0.7 0.75){$0$}\htext(-1.3 0.3){$0$}
\htext(-1.7 0.75){$1$}\htext(-0.5 1.5){$2$}\htext(-0.5 2.5){$3$}
\htext(-0.5 3.5){$2$}
\esegment
\move(5 8)
\bsegment
\move(0 0)\dtri \move(-1 0)\dtri
\move(-2 1)\lvec(0 1)\lvec(-1 0)\lvec(-1 3)\lvec(0 3)\lvec(0 0)\lvec(-2 0)
\lvec(-2 2)\lvec(0 2)
\move(-2 0)\lvec(-1 1)
\htext(-0.3 0.3){$1$}\htext(-0.7 0.75){$0$}\htext(-1.3 0.3){$0$}
\htext(-1.7 0.75){$1$}\htext(-0.5 1.5){$2$}\htext(-1.5 1.5){$2$}
\htext(-0.5 2.5){$3$}
\esegment
\move(-3 16)
\bsegment
\move(0 0)\dtri
\move(-1 1)\lvec(0 1)\lvec(-1 0)\lvec(-1 4)\lvec(0 4)\lvec(0 0)\lvec(-1 0)
\move(0 2)\lvec(-1 2)\move(0 3)\lvec(-1 3)
\htext(-0.3 0.3){$1$}\htext(-0.7 0.75){$0$}\htext(-0.5 1.5){$2$}
\htext(-0.5 2.5){$3$}\htext(-0.5 3.5){$2$}
\esegment
\move(3 16)
\bsegment
\move(0 0)\dtri \move(-1 0)\dtri
\move(-1 1)\lvec(-2 0)\lvec(-2 1)\lvec(0 1)\lvec(-1 0)\lvec(-1 3)\lvec(0 3)
\lvec(0 0)\lvec(-2 0)\move(0 2)\lvec(-1 2)
\htext(-0.3 0.3){$1$}\htext(-0.7 0.75){$0$}\htext(-1.3 0.3){$0$}
\htext(-1.7 0.75){$1$}\htext(-0.5 1.5){$2$}
\htext(-0.5 2.5){$3$}
\esegment
\move(-3 23)
\bsegment
\move(0 0)\dtri
\move(-1 1)\lvec(0 1)\lvec(-1 0)\lvec(-1 3)\lvec(0 3)\lvec(0 0)\lvec(-1 0)
\move(0 2)\lvec(-1 2)
\htext(-0.3 0.3){$1$}\htext(-0.7 0.75){$0$}\htext(-0.5 1.5){$2$}
\htext(-0.5 2.5){$3$}
\esegment
\move(3 23)
\bsegment
\move(0 0)\dtri \move(-1 0)\dtri
\move(-1 1)\lvec(-2 0)\lvec(-2 1)\lvec(0 1)\lvec(-1 0)\lvec(-1 2)\lvec(0 2)
\lvec(0 0)\lvec(-2 0)
\htext(-0.3 0.3){$1$}\htext(-0.7 0.75){$0$}\htext(-1.3 0.3){$0$}
\htext(-1.7 0.75){$1$}\htext(-0.5 1.5){$2$}
\esegment
\move(0 29)
\bsegment
\move(0 0)\dtri
\move(-1 1)\lvec(0 1)\lvec(-1 0)\lvec(0 0)\lvec(0 2)\lvec(-1 2)\lvec(-1 0)
\htext(-0.3 0.3){$1$}\htext(-0.7 0.75){$0$}\htext(-0.5 1.5){$2$}
\esegment
\move(0 34)
\bsegment
\move(0 0)\dtri
\move(-1 0)\lvec(-1 1)\lvec(0 1)\lvec(-1 0)\lvec(0 0)\lvec(0 1)
\htext(-0.3 0.3){$1$}\htext(-0.7 0.75){$0$}
\esegment
\move(-5.5 7.7)\ravec(0 -2.4)\htext(-5.9 6.6){$1$}
\move(-0.5 7.7)\ravec(0 -2.4)\htext(-0.1 6.6){$1$}
\move(4.3 7.7)\ravec(0 -3.4)\htext(4.7 6.3){$2$}
\move(-1.3 7.7)\ravec(-3.5 -2.4)\htext(-3.4 6.7){$0$}
\move(-3.7 15.7)\ravec(-1.7 -2.4)\htext(-4.9 14.7){$0$}
\move(-3.3 15.7)\ravec(2 -3.4)\htext(-1.9 14.3){$1$}
\move(2.1 15.7)\ravec(1.8 -4.4)\htext(3.4 13.7){$2$}
\move(-3.5 22.7)\ravec(0 -2.4)\htext(-3.9 21.5){$2$}
\move(-3.1 22.7)\ravec(5.1 -3.4)\htext(-0.5 21.5){$1$}
\move(2.5 22.7)\ravec(0 -3.4)\htext(2.9 21.5){$3$}
\move(-0.9 28.7)\ravec(-2.2 -2.4)\htext(-2.6 27.5){$3$}
\move(-0.1 28.7)\ravec(2.2 -3.4)\htext(1.4 27.4){$1$}
\move(-0.5 33.7)\ravec(0 -2.4)\htext(-0.1 32.7){$2$}
\htext(-0.45 37.9){$Y_{\La_0}$}
\vtext(-6 -1){$\cdots$}\vtext(-1 -1){$\cdots$}\vtext(4 -1){$\cdots$}
\move(-0.5 37.3)\ravec(0 -2)\htext(-0.1 36.4){$0$}
\move(0 38.3)\move(0 -1.5)
\end{texdraw}
\end{center}

%%%%%%%%%%%%%%%%%%%%%

\newpage

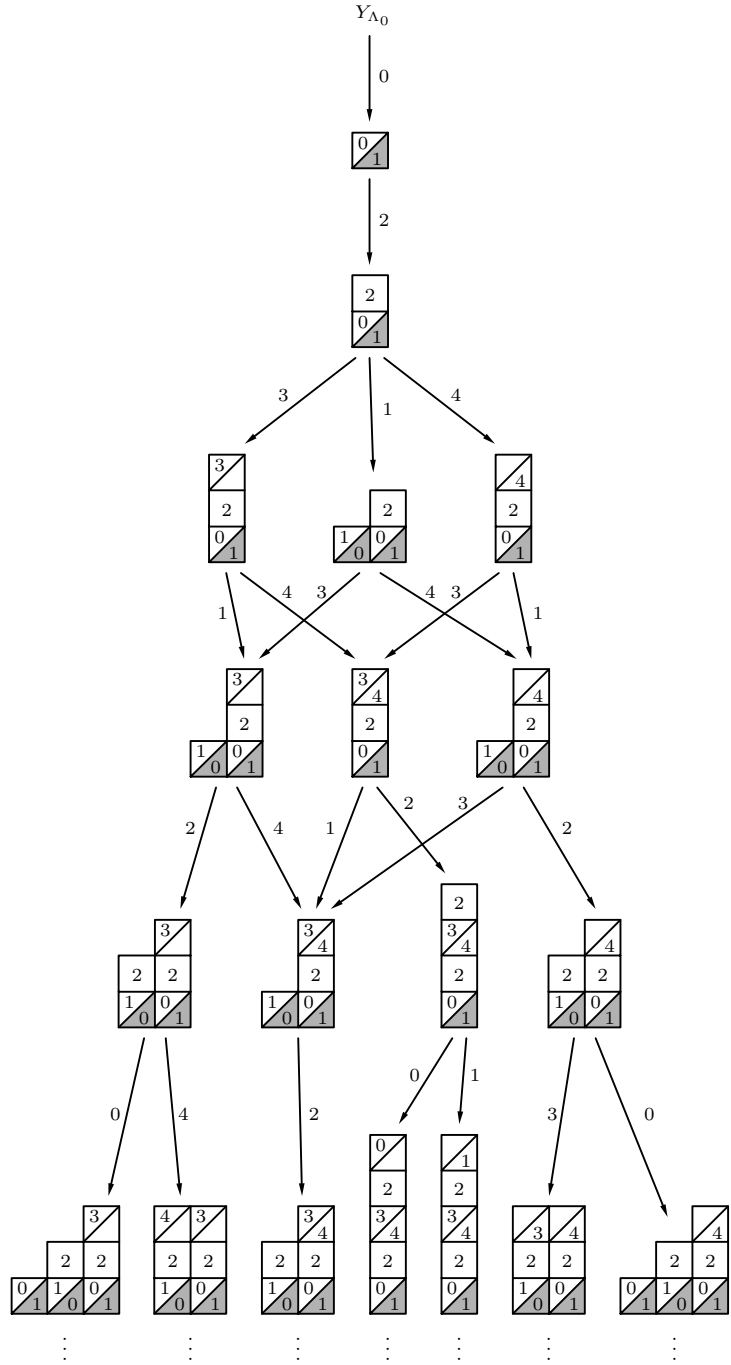
(c) The crystal $\Y(\La_0)$ for $D_4^{(1)}$

\vskip 3mm 

\begin{center}
\begin{texdraw}
\fontsize{7}{7}\selectfont
\drawdim em
\setunitscale 1.7
\nc{\dtri}{
\bsegment
\move(-1 0)\lvec(0 1)\lvec(0 0)\lvec(-1 0)\ifill f:0.7
\esegment
}
%\htext(-8 23){\normalsize$\bullet$ $\Y(\La_0)$ for $D_4^{(1)}$}
\htext(-1.9 28.3){$Y_{\La_0}$}
\move(-1.5 24)
\bsegment
\move(0 0)\dtri
\move(0 1)\lvec(-1 1)\lvec(-1 0)\lvec(0 1)\lvec(0 0)\lvec(-1 0)
\htext(-0.3 0.3){$1$}
\htext(-0.7 0.75){$0$}
\esegment
\move(-1.5 19)
\bsegment
\move(0 0)\dtri
\move(-1 1)\lvec(0 1)\lvec(-1 0)\lvec(-1 2)\lvec(0 2)\lvec(0 0)\lvec(-1 0)
\htext(-0.3 0.3){$1$} \htext(-0.7 0.75){$0$}\htext(-0.5 1.5){$2$}
\esegment
\move(-5.5 13)
\bsegment
\move(0 0)\dtri
\move(-1 1)\lvec(0 1)\lvec(-1 0)\lvec(-1 3)\lvec(0 3)\lvec(0 0)\lvec(-1 0)
\move(0 2)\lvec(-1 2)\lvec(0 3)
\htext(-0.3 0.3){$1$} \htext(-0.7 0.75){$0$}\htext(-0.5 1.5){$2$}
\htext(-0.7 2.75){$3$}
\esegment
\move(-1 13)
\bsegment
\move(0 0)\dtri\move(-1 0)\dtri
\move(-1 1)\lvec(-2 0)\lvec(-2 1)\lvec(0 1)\lvec(-1 0)\lvec(-1 2)\lvec(0 2)
\lvec(0 0)\lvec(-2 0)
\htext(-0.3 0.3){$1$} \htext(-0.7 0.75){$0$}\htext(-1.3 0.3){$0$}
\htext(-1.7 0.75){$1$}\htext(-0.5 1.5){$2$}
\esegment
\move(2.5 13)
\bsegment
\move(0 0)\dtri
\move(-1 1)\lvec(0 1)\lvec(-1 0)\lvec(-1 3)\lvec(0 3)\lvec(0 0)\lvec(-1 0)
\move(0 2)\lvec(-1 2)\lvec(0 3)
\htext(-0.3 0.3){$1$} \htext(-0.7 0.75){$0$}\htext(-0.5 1.5){$2$}
\htext(-0.3 2.3){$4$}
\esegment
\move(-5 7)
\bsegment
\move(0 0)\dtri \move(-1 0)\dtri
\move(-1 1)\lvec(-2 0)\lvec(-2 1)\lvec(0 1)\lvec(-1 0)\lvec(-1 3)\lvec(0 3)
\lvec(0 0)\lvec(-2 0)
\move(0 2)\lvec(-1 2)\lvec(0 3)
\htext(-0.3 0.3){$1$} \htext(-0.7 0.75){$0$} \htext(-1.3 0.3){$0$}
\htext(-1.7 0.75){$1$} \htext(-0.5 1.5){$2$}\htext(-0.7 2.75){$3$}
\esegment
\move(-1.5 7)
\bsegment
\move(0 0)\dtri
\move(-1 1)\lvec(0 1)\lvec(-1 0)\lvec(-1 3)\lvec(0 3)\lvec(0 0)\lvec(-1 0)
\move(0 2)\lvec(-1 2)\lvec(0 3)
\htext(-0.3 0.3){$1$} \htext(-0.7 0.75){$0$}\htext(-0.5 1.5){$2$}
\htext(-0.7 2.75){$3$}\htext(-0.3 2.3){$4$}
\esegment
\move(3 7)
\bsegment
\move(0 0)\dtri \move(-1 0)\dtri
\move(-1 1)\lvec(-2 0)\lvec(-2 1)\lvec(0 1)\lvec(-1 0)\lvec(-1 3)\lvec(0 3)
\lvec(0 0)\lvec(-2 0)
\move(0 2)\lvec(-1 2)\lvec(0 3)
\htext(-0.3 0.3){$1$} \htext(-0.7 0.75){$0$} \htext(-1.3 0.3){$0$}
\htext(-1.7 0.75){$1$} \htext(-0.5 1.5){$2$} \htext(-0.3 2.3){$4$}
\esegment
\move(-7 0)
\bsegment
\move(0 0)\dtri \move(-1 0)\dtri
\move(-2 1)\lvec(0 1)\lvec(-1 0)\lvec(-1 3)\lvec(0 3)\lvec(0 0)\lvec(-2 0)
\lvec(-2 2)\lvec(0 2)
\move(-2 0)\lvec(-1 1)\move(-1 2)\lvec(0 3)
\htext(-0.3 0.3){$1$} \htext(-0.7 0.75){$0$} \htext(-1.3 0.3){$0$}
\htext(-1.7 0.75){$1$} \htext(-0.5 1.5){$2$} \htext(-0.7 2.75){$3$}
\htext(-1.5 1.5){$2$}
\esegment
\move(-3 0)
\bsegment
\move(0 0)\dtri \move(-1 0)\dtri
\move(-1 1)\lvec(-2 0)\lvec(-2 1)\lvec(0 1)\lvec(-1 0)\lvec(-1 3)\lvec(0 3)
\lvec(0 0)\lvec(-2 0)
\move(0 2)\lvec(-1 2)\lvec(0 3)
\htext(-0.3 0.3){$1$} \htext(-0.7 0.75){$0$} \htext(-1.3 0.3){$0$}
\htext(-1.7 0.75){$1$} \htext(-0.5 1.5){$2$} \htext(-0.3 2.3){$4$}
\htext(-0.7 2.75){$3$}
\esegment
\move(1 0)
\bsegment
\move(0 0)\dtri
\move(-1 1)\lvec(0 1)\lvec(-1 0)\lvec(-1 4)\lvec(0 4)\lvec(0 0)\lvec(-1 0)
\move(0 2)\lvec(-1 2)\lvec(0 3)\lvec(-1 3)
\htext(-0.3 0.3){$1$} \htext(-0.7 0.75){$0$} \htext(-0.5 1.5){$2$}
\htext(-0.3 2.3){$4$} \htext(-0.7 2.75){$3$} \htext(-0.5 3.5){$2$}
\esegment
\move(5 0)
\bsegment
\move(0 0)\dtri \move(-1 0)\dtri
\move(-2 1)\lvec(0 1)\lvec(-1 0)\lvec(-1 3)\lvec(0 3)\lvec(0 0)\lvec(-2 0)
\lvec(-2 2)\lvec(0 2)
\move(-2 0)\lvec(-1 1)\move(-1 2)\lvec(0 3)
\htext(-0.3 0.3){$1$} \htext(-0.7 0.75){$0$} \htext(-1.3 0.3){$0$}
\htext(-1.7 0.75){$1$} \htext(-0.5 1.5){$2$} \htext(-0.3 2.3){$4$}
\htext(-1.5 1.5){$2$}
\esegment
\move(8 -8)
\bsegment
\move(0 0)\dtri \move(-1 0)\dtri\move(-2 0)\dtri
\move(-2 1)\lvec(0 1)\lvec(-1 0)\lvec(-1 3)\lvec(0 3)\lvec(0 0)\lvec(-2 0)
\lvec(-2 2)\lvec(0 2)
\move(-2 0)\lvec(-1 1)\move(-1 2)\lvec(0 3)
\move(-2 0)\lvec(-3 0)\lvec(-3 1)\lvec(-2 1)\lvec(-3 0)
\htext(-0.3 0.3){$1$} \htext(-0.7 0.75){$0$} \htext(-1.3 0.3){$0$}
\htext(-1.7 0.75){$1$} \htext(-0.5 1.5){$2$}
\htext(-1.5 1.5){$2$}\htext(-0.3 2.3){$4$}
\htext(-2.3 0.3){$1$}\htext(-2.7 0.75){$0$}
\esegment
\move(4 -8)
\bsegment
\move(0 0)\dtri \move(-1 0)\dtri
\move(-2 1)\lvec(0 1)\lvec(-1 0)\lvec(-1 3)\lvec(0 3)\lvec(0 0)\lvec(-2 0)
\lvec(-2 2)\lvec(0 2)
\move(-2 0)\lvec(-1 1)\move(-1 2)\lvec(0 3)
\move(-2 2)\lvec(-2 3)\lvec(-1 3)\lvec(-2 2)
\htext(-0.3 0.3){$1$} \htext(-0.7 0.75){$0$} \htext(-1.3 0.3){$0$}
\htext(-1.7 0.75){$1$} \htext(-0.5 1.5){$2$} \htext(-0.3 2.3){$4$}
\htext(-1.5 1.5){$2$}\htext(-1.3 2.3){$3$}
\esegment
\move(1 -8)
\bsegment
\move(0 0)\dtri
\move(-1 1)\lvec(0 1)\lvec(-1 0)\lvec(-1 5)\lvec(0 5)\lvec(0 0)\lvec(-1 0)
\move(0 2)\lvec(-1 2)\lvec(0 3)\lvec(-1 3)
\move(0 4)\lvec(-1 4)\lvec(0 5)
\htext(-0.3 0.3){$1$} \htext(-0.7 0.75){$0$} \htext(-0.5 1.5){$2$}
\htext(-0.3 2.3){$4$} \htext(-0.7 2.75){$3$} \htext(-0.5 3.5){$2$}
\htext(-0.3 4.3){$1$}
\esegment
\move(-1 -8)
\bsegment
\move(0 0)\dtri
\move(-1 1)\lvec(0 1)\lvec(-1 0)\lvec(-1 5)\lvec(0 5)\lvec(0 0)\lvec(-1 0)
\move(0 2)\lvec(-1 2)\lvec(0 3)\lvec(-1 3)
\move(0 4)\lvec(-1 4)\lvec(0 5)
\htext(-0.3 0.3){$1$} \htext(-0.7 0.75){$0$} \htext(-0.5 1.5){$2$}
\htext(-0.3 2.3){$4$} \htext(-0.7 2.75){$3$} \htext(-0.5 3.5){$2$}
\htext(-0.7 4.75){$0$}
\esegment
\move(-3 -8)
\bsegment
\move(0 0)\dtri \move(-1 0)\dtri
\move(-2 1)\lvec(0 1)\lvec(-1 0)\lvec(-1 3)\lvec(0 3)\lvec(0 0)\lvec(-2 0)
\lvec(-2 2)\lvec(0 2)
\move(-2 0)\lvec(-1 1)\move(-1 2)\lvec(0 3)
\htext(-0.3 0.3){$1$} \htext(-0.7 0.75){$0$} \htext(-1.3 0.3){$0$}
\htext(-1.7 0.75){$1$} \htext(-0.5 1.5){$2$} \htext(-0.3 2.3){$4$}
\htext(-1.5 1.5){$2$}\htext(-0.7 2.75){$3$}
\esegment
\move(-6 -8)
\bsegment
\move(0 0)\dtri \move(-1 0)\dtri
\move(-2 1)\lvec(0 1)\lvec(-1 0)\lvec(-1 3)\lvec(0 3)\lvec(0 0)\lvec(-2 0)
\lvec(-2 2)\lvec(0 2)
\move(-2 0)\lvec(-1 1)\move(-1 2)\lvec(0 3)
\move(-2 2)\lvec(-2 3)\lvec(-1 3)\lvec(-2 2)
\htext(-0.3 0.3){$1$} \htext(-0.7 0.75){$0$} \htext(-1.3 0.3){$0$}
\htext(-1.7 0.75){$1$} \htext(-0.5 1.5){$2$} \htext(-1.7 2.75){$4$}
\htext(-1.5 1.5){$2$}\htext(-0.7 2.75){$3$}
\esegment
\move(-9 -8)
\bsegment
\move(0 0)\dtri \move(-1 0)\dtri\move(-2 0)\dtri
\move(-2 1)\lvec(0 1)\lvec(-1 0)\lvec(-1 3)\lvec(0 3)\lvec(0 0)\lvec(-2 0)
\lvec(-2 2)\lvec(0 2)
\move(-2 0)\lvec(-1 1)\move(-1 2)\lvec(0 3)
\move(-2 0)\lvec(-3 0)\lvec(-3 1)\lvec(-2 1)\lvec(-3 0)
\htext(-0.3 0.3){$1$} \htext(-0.7 0.75){$0$} \htext(-1.3 0.3){$0$}
\htext(-1.7 0.75){$1$} \htext(-0.5 1.5){$2$}
\htext(-1.5 1.5){$2$}\htext(-0.7 2.75){$3$}
\htext(-2.3 0.3){$1$}\htext(-2.7 0.75){$0$}
\esegment
\move(-2 27.7)\ravec(0 -2.4)\htext(-1.6 26.6){$0$}
\move(-2 23.7)\ravec(0 -2.4)\htext(-1.6 22.6){$2$}
\move(-1.6 18.7)\ravec(3.1 -2.4)\htext(0.4 17.7){$4$}
\move(-2 18.7)\ravec(0.1 -3.2)\htext(-1.5 17.3){$1$}
\move(-2.4 18.7)\ravec(-3.1 -2.4)\htext(-4.4 17.7){$3$}
\move(-6 12.7)\ravec(0.5 -2.4)\htext(-6.1 11.6){$1$}
\move(-5.6 12.7)\ravec(3.2 -2.4)\htext(-4.3 12.2){$4$}
\move(2 12.7)\ravec(0.5 -2.4)\htext(2.7 11.6){$1$}
\move(1.6 12.7)\ravec(-3.2 -2.4)\htext(0.4 12.2){$3$}
\move(-2.3 12.7)\ravec(-2.8 -2.4)\htext(-3.35 12.2){$3$}
\move(-1.7 12.7)\ravec(3.8 -2.4)\htext(-0.3 12.2){$4$}
\move(-6.3 6.7)\ravec(-1 -3.4)\htext(-7 5.6){$2$}
\move(-5.7 6.7)\ravec(1.8 -3.4)\htext(-4.55 5.6){$4$}
\move(-2.2 6.7)\ravec(-1.3 -3.4)\htext(-3.1 5.6){$1$}
\move(-1.8 6.7)\ravec(1.9 -2.4)\htext(-0.9 6.3){$2$}
\move(1.7 6.7)\ravec(-4.8 -3.4)\htext(0.6 6.3){$3$}
\move(2.3 6.7)\ravec(2 -3.4)\htext(3.5 5.6){$2$}
\move(-8.3 -0.3)\ravec(-1 -4.4)\htext(-9.1 -2.4){$0$}
\move(-7.7 -0.3)\ravec(0.4 -4.4)\htext(-7.2 -2.4){$4$}
\move(-4 -0.3)\ravec(0.1 -4.4)\htext(-3.55 -2.4){$2$}
\move(0.3 -0.3)\ravec(-1.5 -2.4)\htext(-0.75 -1.3){$0$}
\move(0.7 -0.3)\ravec(-0.2 -2.4)\htext(0.95 -1.3){$1$}
\move(3.7 -0.3)\ravec(-0.7 -4.4)\htext(3.1 -2.4){$3$}
\move(4.3 -0.3)\ravec(2.0 -5)\htext(5.8 -2.4){$0$}
\vtext(-10.5 -8.9){$\cdots$}
\vtext(-7 -8.9){$\cdots$}
\vtext(-4 -8.9){$\cdots$}
\vtext(-1.5 -8.9){$\cdots$}
\vtext(0.5 -8.9){$\cdots$}
\vtext(3 -8.9){$\cdots$}
\vtext(6.5 -8.9){$\cdots$}
\move(0 -9.6)\move(0 28.8)
\end{texdraw}
\end{center}

\newpage

(d) The crystal $\Y(\La_0)$ for $A_4^{(2)}$ 

\vskip 3mm 

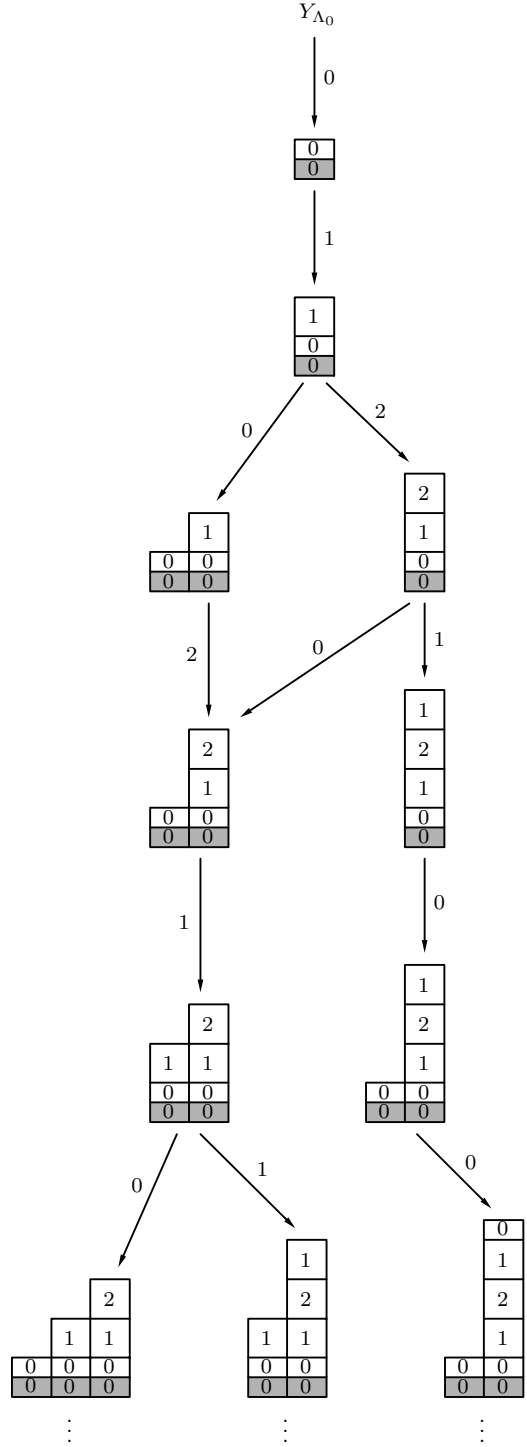
\begin{center}
\begin{texdraw}
\fontsize{8}{8}\selectfont
\drawdim em
\setunitscale 1.75
\htext(-0.2 28.2){$Y_{\La_0}$}
\move(0.2 24)
\bsegment
\move(0 0)\lvec(-1 0)\lvec(-1 0.5)\lvec(0 0.5) \ifill f:0.7
\move(0 0)\lvec(-1 0)\lvec(-1 1)\lvec(0 1)\lvec(0 0)
\move(0 0.5)\lvec(-1 0.5)
\htext(-0.5 0.28){$0$} \htext(-0.5 0.78){$0$}
\esegment
\move(0.2 19)
\bsegment
\move(0 0)\lvec(-1 0)\lvec(-1 0.5)\lvec(0 0.5) \ifill f:0.7
\move(0 0)\lvec(-1 0)\lvec(-1 2)\lvec(0 2)\lvec(0 0)
\move(0 0.5)\lvec(-1 0.5)
\move(0 1)\lvec(-1 1)
\htext(-0.5 0.28){$0$} \htext(-0.5 0.78){$0$}
\htext(-0.5 1.5){$1$}
\esegment
\move(-2.5 13.5)
\bsegment
\move(0 0)\lvec(-2 0)\lvec(-2 0.5)\lvec(0 0.5) \ifill f:0.7
\move(-1 0)\lvec(-1 2)\lvec(0 2)\lvec(0 0)\lvec(-2 0)\lvec(-2 1)\lvec(0 1)
\move(0 0.5)\lvec(-2 0.5)
\htext(-0.5 0.28){$0$} \htext(-0.5 0.78){$0$}
\htext(-1.5 0.28){$0$} \htext(-1.5 0.78){$0$}
\htext(-0.5 1.5){$1$}
\esegment
\move(3 13.5)
\bsegment
\move(0 0)\lvec(-1 0)\lvec(-1 0.5)\lvec(0 0.5) \ifill f:0.7
\move(0 0)\lvec(-1 0)\lvec(-1 3)\lvec(0 3)\lvec(0 0)
\move(0 0.5)\lvec(-1 0.5)
\move(0 1)\lvec(-1 1)
\move(0 2)\lvec(-1 2)
\htext(-0.5 0.28){$0$} \htext(-0.5 0.78){$0$}
\htext(-0.5 1.5){$1$}
\htext(-0.5 2.5){$2$}
\esegment
\move(-2.5 7)
\bsegment
\move(0 0)\lvec(-2 0)\lvec(-2 0.5)\lvec(0 0.5) \ifill f:0.7
\move(-1 0)\lvec(-1 3)\lvec(0 3)\lvec(0 0)\lvec(-2 0)\lvec(-2 1)\lvec(0 1)
\move(0 0.5)\lvec(-2 0.5)
\move(0 2)\lvec(-1 2)
\htext(-0.5 0.28){$0$} \htext(-0.5 0.78){$0$}
\htext(-1.5 0.28){$0$} \htext(-1.5 0.78){$0$}
\htext(-0.5 1.5){$1$}
\htext(-0.5 2.5){$2$}
\esegment
\move(3 7)
\bsegment
\move(0 0)\lvec(-1 0)\lvec(-1 0.5)\lvec(0 0.5) \ifill f:0.7
\move(0 0)\lvec(-1 0)\lvec(-1 4)\lvec(0 4)\lvec(0 0)
\move(0 0.5)\lvec(-1 0.5)
\move(0 1)\lvec(-1 1)
\move(0 2)\lvec(-1 2)
\move(0 3)\lvec(-1 3)
\htext(-0.5 0.28){$0$} \htext(-0.5 0.78){$0$}
\htext(-0.5 1.5){$1$}
\htext(-0.5 2.5){$2$}
\htext(-0.5 3.5){$1$}
\esegment
\move(-2.5 0)
\bsegment
\move(0 0)\lvec(-2 0)\lvec(-2 0.5)\lvec(0 0.5) \ifill f:0.7
\move(-1 0)\lvec(-1 3)\lvec(0 3)\lvec(0 0)\lvec(-2 0)\lvec(-2 2)\lvec(0 2)
\move(0 1)\rlvec(-2 0)
\move(0 0.5)\rlvec(-2 0)
\htext(-0.5 0.28){$0$} \htext(-0.5 0.78){$0$}
\htext(-1.5 0.28){$0$} \htext(-1.5 0.78){$0$}
\htext(-0.5 1.5){$1$} \htext(-1.5 1.5){$1$}
\htext(-0.5 2.5){$2$}
\esegment
\move(3 0)
\bsegment
\move(0 0)\lvec(-2 0)\lvec(-2 0.5)\lvec(0 0.5) \ifill f:0.7
\move(-1 0)\lvec(-1 4)\lvec(0 4)\lvec(0 0)\lvec(-2 0)\lvec(-2 1)\lvec(0 1)
\move(-1 2)\rlvec(1 0)
\move(-1 3)\rlvec(1 0)
\move(-2 0.5)\rlvec(2 0)
\htext(-0.5 0.28){$0$} \htext(-0.5 0.78){$0$}
\htext(-1.5 0.28){$0$} \htext(-1.5 0.78){$0$}
\htext(-0.5 1.5){$1$}
\htext(-0.5 2.5){$2$}
\htext(-0.5 3.5){$1$}
\esegment
\move(-5 -7)
\bsegment
\move(0 0)\lvec(-3 0)\lvec(-3 0.5)\lvec(0 0.5) \ifill f:0.7
\move(0 1)\lvec(-3 1)\lvec(-3 0)\lvec(0 0)\lvec(0 3)\lvec(-1 3)\lvec(-1 0)
\move(0 0.5)\rlvec(-3 0)
\move(-2 0)\lvec(-2 2)\lvec(0 2)
\htext(-0.5 0.28){$0$} \htext(-0.5 0.78){$0$}
\htext(-1.5 0.28){$0$} \htext(-1.5 0.78){$0$}
\htext(-2.5 0.3){$0$} \htext(-2.5 0.8){$0$}
\htext(-0.5 1.5){$1$} \htext(-1.5 1.5){$1$}
\htext(-0.5 2.5){$2$}
\esegment
\move(0 -7)
\bsegment
\move(0 0)\lvec(-2 0)\lvec(-2 0.5)\lvec(0 0.5) \ifill f:0.7
\move(-1 0)\lvec(-1 4)\lvec(0 4)\lvec(0 0)\lvec(-2 0)\lvec(-2 2)\lvec(0 2)
\move(0 1)\lvec(-2 1)
\move(0 0.5)\lvec(-2 0.5)
\move(0 3)\lvec(-1 3)
\htext(-0.5 0.28){$0$} \htext(-0.5 0.78){$0$}
\htext(-1.5 0.28){$0$} \htext(-1.5 0.78){$0$}
\htext(-0.5 1.5){$1$} \htext(-1.5 1.5){$1$}
\htext(-0.5 2.5){$2$}
\htext(-0.5 3.5){$1$}
\esegment
\move(5 -7)
\bsegment
\move(0 0)\lvec(-2 0)\lvec(-2 0.5)\lvec(0 0.5) \ifill f:0.7
\move(-1 0)\lvec(-1 4.5)\lvec(0 4.5)\lvec(0 0)\lvec(-2 0)\lvec(-2 1)\lvec(0 1)
\move(0 0.5)\lvec(-2 0.5)
\move(0 2)\lvec(-1 2)\move(0 3)\lvec(-1 3)\move(0 4)\lvec(-1 4)
\htext(-0.5 0.28){$0$} \htext(-0.5 0.78){$0$}
\htext(-1.5 0.28){$0$} \htext(-1.5 0.78){$0$}
\htext(-0.5 1.5){$1$}
\htext(-0.5 2.5){$2$}
\htext(-0.5 3.5){$1$}
\htext(-0.5 4.3){$0$}
\esegment
\move(-0.3 27.6)\avec(-0.3 25.3)\htext(0.1 26.6){$0$}
\move(-0.3 23.7)\avec(-0.3 21.3)\htext(0.1 22.5){$1$}
\move(0 18.8)\avec(2.1 16.8)\htext(1.4 18.1){$2$}
\move(-0.6 18.8)\avec(-2.8 15.8)\htext(-2 17.6){$0$}
\move(2.5 13.2)\avec(2.5 11.3)\htext(2.9 12.3){$1$}
\move(-3 13.2)\avec(-3 10.3)\htext(-3.4 11.9){$2$}
\move(2.1 13.2)\avec(-2.2 10.3)\htext(-0.2 12.1){$0$}
\move(2.5 6.7)\avec(2.5 4.3)\htext(2.9 5.6){$0$}
\move(-3.2 6.7)\avec(-3.2 3.3)\htext(-3.6 5.1){$1$}
\move(2.3 -0.3)\avec(4.2 -2.2)\htext(3.7 -1){$0$}
\move(-3.2 -0.3)\avec(-0.8 -2.7)\htext(-1.6 -1.2){$1$}
\move(-3.8 -0.3)\avec(-5.3 -3.7)\htext(-4.8 -1.6){$0$}
\vtext(4 -7.8){$\cdots$}
\vtext(-1 -7.8){$\cdots$}
\vtext(-6.5 -7.8){$\cdots$}
%\htext(-5.4 24){\normalsize $\bullet$ $\Y(\La_0)$ for $A_4^{(2)}$}
\move(6 -8.5)\move(0 28.6)
\end{texdraw}
\end{center}

\newpage

(e) The crystal $\Y(\La_0)$ for $D^{(2)}_3$

\vskip 3mm

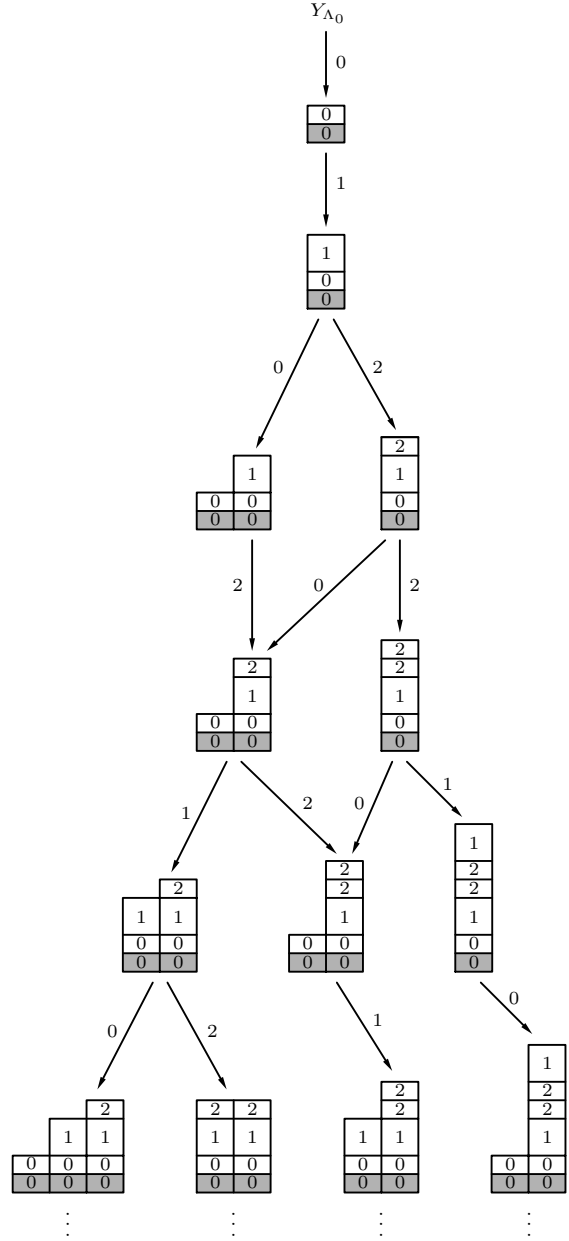
\begin{center}
\begin{texdraw}
\fontsize{7}{7}\selectfont
\drawdim em
\setunitscale 1.75
%\htext(-5.8 21.5){\normalsize $\bullet$ $\Y(\La_0)$ for $D^{(2)}_3$}
\htext(-0.45 26){$Y_{\La_0}$}
\move(0 22.5)
\bsegment
\move(0 0)\lvec(-1 0)\lvec(-1 0.5)\lvec(0 0.5) \ifill f:0.7
\move(-1 0)\lvec(-1 1)\lvec(0 1)\lvec(0 0)\lvec(-1 0)
\move(-1 0.5)\lvec(0 0.5)
\htext(-0.5 0.28){$0$} \htext(-0.5 0.78){$0$}
\esegment
\move(0 18)
\bsegment
\move(0 0)\lvec(-1 0)\lvec(-1 0.5)\lvec(0 0.5) \ifill f:0.7
\move(-1 0)\lvec(-1 2)\lvec(0 2)\lvec(0 0)\lvec(-1 0)
\move(-1 0.5)\lvec(0 0.5)
\move(-1 1)\lvec(0 1)
\htext(-0.5 0.28){$0$} \htext(-0.5 0.78){$0$}
\htext(-0.5 1.5){$1$}
\esegment
\move(-2 12)
\bsegment
\move(0 0)\lvec(-2 0)\lvec(-2 0.5)\lvec(0 0.5) \ifill f:0.7
\move(-1 0)\lvec(-1 2)\lvec(0 2)\lvec(0 0)\lvec(-2 0)\lvec(-2 1)\lvec(0 1)
\move(-2 0.5)\lvec(0 0.5)
\htext(-0.5 0.28){$0$} \htext(-0.5 0.78){$0$}
\htext(-1.5 0.28){$0$} \htext(-1.5 0.78){$0$}
\htext(-0.5 1.5){$1$}
\esegment
\move(2 12)
\bsegment
\move(0 0)\lvec(-1 0)\lvec(-1 0.5)\lvec(0 0.5) \ifill f:0.7
\move(-1 0)\lvec(-1 2.5)\lvec(0 2.5)\lvec(0 0)\lvec(-1 0)
\move(-1 0.5)\lvec(0 0.5)
\move(-1 1)\lvec(0 1)
\move(-1 2)\lvec(0 2)
\move(-1 2.5)\lvec(0 2.5)
\htext(-0.5 0.28){$0$} \htext(-0.5 0.78){$0$}
\htext(-0.5 1.5){$1$}
\htext(-0.5 2.28){$2$}
\esegment
\move(-2 6)
\bsegment
\move(0 0)\lvec(-2 0)\lvec(-2 0.5)\lvec(0 0.5) \ifill f:0.7
\move(-1 0)\lvec(-1 2.5)\lvec(0 2.5)\lvec(0 0)\lvec(-2 0)\lvec(-2 1)\lvec(0 1)
\move(-2 0.5)\lvec(0 0.5)
\move(-1 2)\lvec(0 2)
\htext(-0.5 0.28){$0$} \htext(-0.5 0.78){$0$}
\htext(-1.5 0.28){$0$} \htext(-1.5 0.78){$0$}
\htext(-0.5 1.5){$1$}
\htext(-0.5 2.28){$2$}
\esegment
\move(2 6)
\bsegment
\move(0 0)\lvec(-1 0)\lvec(-1 0.5)\lvec(0 0.5) \ifill f:0.7
\move(-1 0)\lvec(-1 3)\lvec(0 3)\lvec(0 0)\lvec(-1 0)
\move(-1 0.5)\lvec(0 0.5)
\move(-1 1)\lvec(0 1)
\move(-1 2)\lvec(0 2)
\move(-1 2.5)\lvec(0 2.5)
\htext(-0.5 0.28){$0$} \htext(-0.5 0.78){$0$}
\htext(-0.5 1.5){$1$}
\htext(-0.5 2.28){$2$}
\htext(-0.5 2.78){$2$}
\esegment
\move(-4 0)
\bsegment
\move(0 0)\lvec(-2 0)\lvec(-2 0.5)\lvec(0 0.5) \ifill f:0.7
\move(-1 0)\lvec(-1 2.5)\lvec(0 2.5)\lvec(0 0)\lvec(-2 0)\lvec(-2 2)\lvec(0 2)
\move(-2 0.5)\lvec(0 0.5)
\move(-2 1)\lvec(0 1)
\htext(-0.5 0.28){$0$} \htext(-0.5 0.78){$0$}
\htext(-1.5 0.28){$0$} \htext(-1.5 0.78){$0$}
\htext(-0.5 1.5){$1$} \htext(-1.5 1.5){$1$}
\htext(-0.5 2.28){$2$}
\esegment
\move(0.5 0)
\bsegment
\move(0 0)\lvec(-2 0)\lvec(-2 0.5)\lvec(0 0.5) \ifill f:0.7
\move(-1 0)\lvec(-1 3)\lvec(0 3)\lvec(0 0)\lvec(-2 0)\lvec(-2 1)\lvec(0 1)
\move(-2 0.5)\lvec(0 0.5)
\move(-1 2)\lvec(0 2)
\move(-1 2.5)\lvec(0 2.5)
\htext(-0.5 0.28){$0$} \htext(-0.5 0.78){$0$}
\htext(-1.5 0.28){$0$} \htext(-1.5 0.78){$0$}
\htext(-0.5 1.5){$1$}
\htext(-0.5 2.28){$2$} \htext(-0.5 2.78){$2$}
\esegment
\move(4 0)
\bsegment
\move(0 0)\lvec(-1 0)\lvec(-1 0.5)\lvec(0 0.5) \ifill f:0.7
\move(-1 0)\lvec(-1 4)\lvec(0 4)\lvec(0 0)\lvec(-1 0)
\move(-1 0.5)\lvec(0 0.5)
\move(-1 1)\lvec(0 1)
\move(-1 2)\lvec(0 2)
\move(-1 2.5)\lvec(0 2.5)
\move(-1 3)\lvec(0 3)
\htext(-0.5 0.28){$0$} \htext(-0.5 0.78){$0$}
\htext(-0.5 1.5){$1$}
\htext(-0.5 2.28){$2$} \htext(-0.5 2.78){$2$}
\htext(-0.5 3.5){$1$}
\esegment
\move(-6 -6)
\bsegment
\move(0 0)\lvec(-3 0)\lvec(-3 0.5)\lvec(0 0.5) \ifill f:0.7
\move(0 1)\lvec(-3 1)\lvec(-3 0)\lvec(0 0)\lvec(0 2.5)\lvec(-1 2.5)\lvec(-1 0)
\move(0 0.5)\rlvec(-3 0)
\move(-2 0)\lvec(-2 2)\lvec(0 2)
\htext(-0.5 0.28){$0$} \htext(-0.5 0.78){$0$}
\htext(-1.5 0.28){$0$} \htext(-1.5 0.78){$0$}
\htext(-2.5 0.3){$0$} \htext(-2.5 0.8){$0$}
\htext(-0.5 1.5){$1$} \htext(-1.5 1.5){$1$}
\htext(-0.5 2.28){$2$}
\esegment
\move(-2 -6)
\bsegment
\move(0 0)\lvec(-2 0)\lvec(-2 0.5)\lvec(0 0.5) \ifill f:0.7
\move(-2 0)\lvec(-2 2.5)\lvec(0 2.5)\lvec(0 0)\lvec(-2 0)
\move(-2 0.5)\lvec(0 0.5)
\move(-2 1)\lvec(0 1)
\move(-2 2)\lvec(0 2)
\move(-1 0)\lvec(-1 2.5)
\htext(-0.5 0.28){$0$} \htext(-0.5 0.78){$0$}
\htext(-1.5 0.28){$0$} \htext(-1.5 0.78){$0$}
\htext(-0.5 1.5){$1$} \htext(-1.5 1.5){$1$}
\htext(-0.5 2.28){$2$} \htext(-1.5 2.28){$2$}
\esegment
\move(2 -6)
\bsegment
\move(0 0)\lvec(-2 0)\lvec(-2 0.5)\lvec(0 0.5) \ifill f:0.7
\move(-1 0)\lvec(-1 3)\lvec(0 3)\lvec(0 0)\lvec(-2 0)\lvec(-2 2)\lvec(0 2)
\move(-2 0.5)\lvec(0 0.5)
\move(-2 1)\lvec(0 1)
\move(-1 2.5)\lvec(0 2.5)
\htext(-0.5 0.28){$0$} \htext(-0.5 0.78){$0$}
\htext(-1.5 0.28){$0$} \htext(-1.5 0.78){$0$}
\htext(-0.5 1.5){$1$} \htext(-1.5 1.5){$1$}
\htext(-0.5 2.28){$2$} \htext(-0.5 2.78){$2$}
\esegment
\move(6 -6)
\bsegment
\move(0 0)\lvec(-2 0)\lvec(-2 0.5)\lvec(0 0.5) \ifill f:0.7
\move(-1 0)\lvec(-1 4)\lvec(0 4)\lvec(0 0)\lvec(-2 0)\lvec(-2 1)\lvec(0 1)
\move(-2 0.5)\lvec(0 0.5)
\move(-1 2)\lvec(0 2)
\move(-1 3)\lvec(0 3)
\move(-1 2.5)\lvec(0 2.5)
\htext(-0.5 0.28){$0$} \htext(-0.5 0.78){$0$}
\htext(-1.5 0.28){$0$} \htext(-1.5 0.78){$0$}
\htext(-0.5 1.5){$1$}
\htext(-0.5 2.28){$2$} \htext(-0.5 2.78){$2$}
\htext(-0.5 3.5){$1$}
\esegment
\vtext(-7.5 -6.8){$\cdots$}
\vtext(-3 -6.8){$\cdots$}
\vtext(1 -6.8){$\cdots$}
\vtext(5 -6.8){$\cdots$}
\move(-0.5 25.5)\avec(-0.5 23.7)\htext(-0.1 24.7){$0$}
\move(-0.5 22.2)\avec(-0.5 20.2)\htext(-0.1 21.4){$1$}
\move(-0.3 17.7)\avec(1.4 14.7)\htext(0.9 16.4){$2$}
\move(-0.7 17.7)\avec(-2.4 14.2)\htext(-1.8 16.4){$0$}
\move(1.5 11.7)\avec(1.5 9.2)\htext(1.9 10.5){$2$}
\move(1.1 11.7)\avec(-2.1 8.7)\htext(-0.7 10.5){$0$}
\move(-2.5 11.7)\avec(-2.5 8.7)\htext(-2.9 10.5){$2$}
\move(1.7 5.7)\avec(3.2 4.2)\htext(2.8 5.1){$1$}
\move(1.3 5.7)\avec(0.2 3.2)\htext(0.4 4.6){$0$}
\move(-2.8 5.7)\avec(-0.2 3.2)\htext(-1 4.6){$2$}
\move(-3.2 5.7)\avec(-4.7 2.7)\htext(-4.3 4.3){$1$}
\move(3.7 -0.3)\avec(5.1 -1.7)\htext(4.6 -0.7){$0$}
\move(-0.2 -0.3)\avec(1.3 -2.7)\htext(0.9 -1.3){$1$}
\move(-4.8 -0.3)\avec(-3.2 -3.2)\htext(-3.6 -1.6){$2$}
\move(-5.2 -0.3)\avec(-6.7 -3.2)\htext(-6.3 -1.6){$0$}
\move(0 -7.4)\move(0 26.5)
\end{texdraw}
\end{center}

\newpage

(f) The crystal $\Y(\La_0)$ for $B_3^{(1)}$

\vskip 3mm 

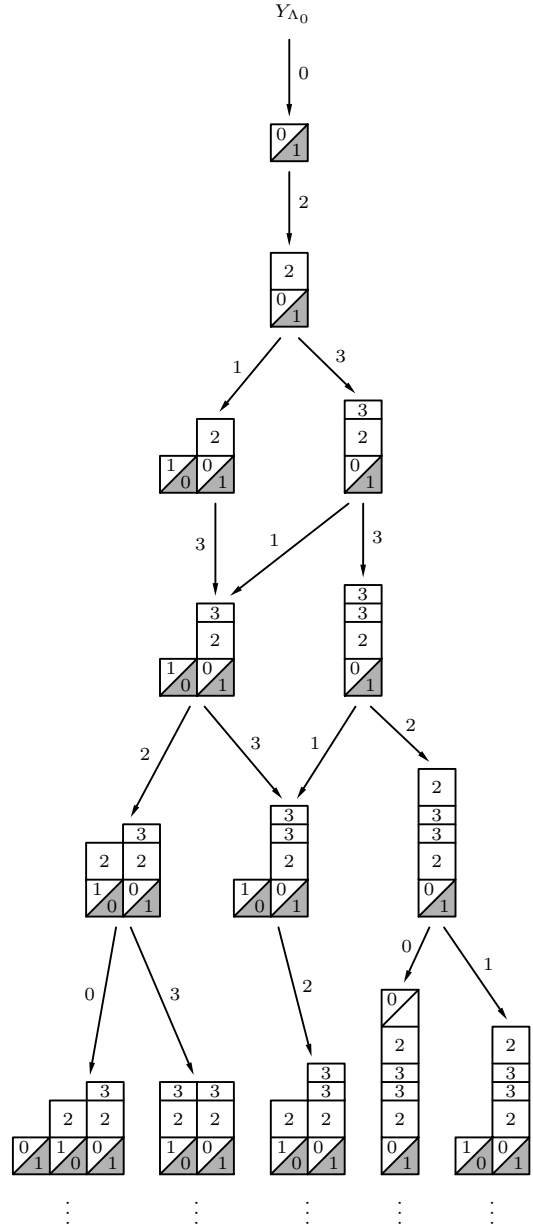
\begin{center}
\begin{texdraw}
\fontsize{7}{7}\selectfont
\drawdim em
\setunitscale 1.75
\nc{\dtri}{
\bsegment
\move(-1 0)\lvec(0 1)\lvec(0 0)\lvec(-1 0)\ifill f:0.7
\esegment
}
%\htext(-7 20){\normalsize $\bullet$ $\Y(\La_0)$ for $B_3^{(1)}$}
\htext(-1.4 24.5){$Y_{\La_0}$}
\move(-1 20.5)
\bsegment
\move(0 0)\dtri
\move(0 1)\lvec(-1 0)\lvec(-1 1)\lvec(0 1)\lvec(0 0)\lvec(-1 0)
\htext(-0.3 0.3){$1$} \htext(-0.7 0.75){$0$}
\esegment
\move(-1 16)
\bsegment
\move(0 0)\dtri
\move(-1 1)\lvec(0 1)\lvec(-1 0)\lvec(-1 2)\lvec(0 2)\lvec(0 0)\lvec(-1 0)
\htext(-0.3 0.3){$1$} \htext(-0.7 0.75){$0$}
\htext(-0.5 1.5){$2$}
\esegment
\move(-3 11.5)
\bsegment
\move(0 0)\dtri \move(-1 0)\dtri
\move(-1 1)\lvec(-2 0)\lvec(-2 1)\lvec(0 1)\lvec(-1 0)\lvec(-1 2)\lvec(0 2)
\lvec(0 0)\lvec(-2 0)
\htext(-0.3 0.3){$1$} \htext(-0.7 0.75){$0$}
\htext(-1.3 0.3){$0$} \htext(-1.7 0.75){$1$}
\htext(-0.5 1.5){$2$}
\esegment
\move(1 11.5)
\bsegment
\move(0 0)\dtri
\move(-1 1)\lvec(0 1)\lvec(-1 0)\lvec(-1 2.5)\lvec(0 2.5)\lvec(0 0)\lvec(-1 0)
\move(-1 2)\lvec(0 2)
\htext(-0.3 0.3){$1$} \htext(-0.7 0.75){$0$}
\htext(-0.5 1.5){$2$}
\htext(-0.5 2.25){$3$}
\esegment
\move(-3 6)
\bsegment
\move(0 0)\dtri \move(-1 0)\dtri
\move(-1 1)\lvec(-2 0)\lvec(-2 1)\lvec(0 1)\lvec(-1 0)\lvec(-1 2.5)\lvec(0 2.5)
\lvec(0 0)\lvec(-2 0)
\move(0 2)\lvec(-1 2)
\htext(-0.3 0.3){$1$} \htext(-0.7 0.75){$0$}
\htext(-1.3 0.3){$0$} \htext(-1.7 0.75){$1$}
\htext(-0.5 1.5){$2$}
\htext(-0.5 2.25){$3$}
\esegment
\move(1 6)
\bsegment
\move(0 0)\dtri
\move(-1 1)\lvec(0 1)\lvec(-1 0)\lvec(-1 3)\lvec(0 3)\lvec(0 0)\lvec(-1 0)
\move(-1 2)\lvec(0 2)
\move(-1 2.5)\lvec(0 2.5)
\htext(-0.3 0.3){$1$} \htext(-0.7 0.75){$0$}
\htext(-0.5 1.5){$2$}
\htext(-0.5 2.25){$3$}
\htext(-0.5 2.75){$3$}
\esegment
\move(-5 0)
\bsegment
\move(0 0)\dtri \move(-1 0)\dtri
\move(-2 1)\lvec(0 1)\lvec(-1 0)\lvec(-1 2.5)\lvec(0 2.5)\lvec(0 0)
\lvec(-2 0)\lvec(-2 2)\lvec(0 2)
\move(-2 0)\lvec(-1 1)
\htext(-0.3 0.3){$1$} \htext(-0.7 0.75){$0$}
\htext(-1.3 0.3){$0$} \htext(-1.7 0.75){$1$}
\htext(-0.5 1.5){$2$}
\htext(-1.5 1.5){$2$}
\htext(-0.5 2.25){$3$}
\esegment
\move(-1 0)
\bsegment
\move(0 0)\dtri \move(-1 0)\dtri
\move(-1 1)\lvec(-2 0)\lvec(-2 1)\lvec(0 1)\lvec(-1 0)\lvec(-1 3)\lvec(0 3)
\lvec(0 0)\lvec(-2 0)
\move(0 2)\lvec(-1 2)
\move(0 2.5)\lvec(-1 2.5)
\htext(-0.3 0.3){$1$} \htext(-0.7 0.75){$0$}
\htext(-1.3 0.3){$0$} \htext(-1.7 0.75){$1$}
\htext(-0.5 1.5){$2$}
\htext(-0.5 2.25){$3$}
\htext(-0.5 2.75){$3$}
\esegment
\move(3 0)
\bsegment
\move(0 0)\dtri
\move(-1 1)\lvec(0 1)\lvec(-1 0)\lvec(-1 4)\lvec(0 4)\lvec(0 0)\lvec(-1 0)
\move(-1 2)\lvec(0 2)
\move(-1 2.5)\lvec(0 2.5)
\move(-1 3)\lvec(0 3)
\htext(-0.3 0.3){$1$} \htext(-0.7 0.75){$0$}
\htext(-0.5 1.5){$2$}
\htext(-0.5 2.25){$3$}
\htext(-0.5 2.75){$3$}
\htext(-0.5 3.5){$2$}
\esegment
\move(-6 -7)
\bsegment
\move(0 0)\dtri \move(-1 0)\dtri \move(-2 0)\dtri
\move(-2 1)\lvec(-3 0)\lvec(-3 1)\lvec(0 1)\lvec(-1 0)\lvec(-1 2.5)
\lvec(0 2.5)\lvec(0 0)\lvec(-3 0)
\move(-1 1)\lvec(-2 0)\lvec(-2 2)\lvec(0 2)
\htext(-0.3 0.3){$1$} \htext(-0.7 0.75){$0$}
\htext(-1.3 0.3){$0$} \htext(-1.7 0.75){$1$}
\htext(-2.3 0.3){$1$} \htext(-2.7 0.75){$0$}
\htext(-0.5 1.5){$2$}
\htext(-1.5 1.5){$2$}
\htext(-0.5 2.25){$3$}
\esegment
\move(-3 -7)
\bsegment
\move(0 0)\dtri \move(-1 0)\dtri
\move(0 0)\lvec(-2 0)\lvec(-2 2.5)\lvec(0 2.5)\lvec(0 0)
\move(-2 0)\lvec(-1 1)
\move(-1 0)\lvec(0 1)\lvec(-2 1)
\move(-1 0)\lvec(-1 2.5)
\move(0 2)\lvec(-2 2)
\htext(-0.3 0.3){$1$} \htext(-0.7 0.75){$0$}
\htext(-1.3 0.3){$0$} \htext(-1.7 0.75){$1$}
\htext(-0.5 1.5){$2$}
\htext(-1.5 1.5){$2$}
\htext(-0.5 2.25){$3$}
\htext(-1.5 2.25){$3$}
\esegment
\move(0 -7)
\bsegment
\move(0 0)\dtri \move(-1 0)\dtri
\move(-1 0)\lvec(-1 3)\lvec(0 3)\lvec(0 0)\lvec(-2 0)\lvec(-2 2)\lvec(0 2)
\move(-2 1)\lvec(0 1)\lvec(-1 0)
\move(-2 0)\lvec(-1 1)
\move(0 2.5)\lvec(-1 2.5)
\htext(-0.3 0.3){$1$} \htext(-0.7 0.75){$0$}
\htext(-1.3 0.3){$0$} \htext(-1.7 0.75){$1$}
\htext(-0.5 1.5){$2$}
\htext(-1.5 1.5){$2$}
\htext(-0.5 2.25){$3$}
\htext(-0.5 2.75){$3$}
\esegment
\move(2 -7)
\bsegment
\move(0 0)\dtri
\move(-1 1)\lvec(0 1)\lvec(-1 0)\lvec(-1 5)\lvec(0 5)\lvec(0 0)\lvec(-1 0)
\move(0 2)\lvec(-1 2)
\move(0 2.5)\lvec(-1 2.5)
\move(0 3)\lvec(-1 3)
\move(0 4)\lvec(-1 4)\lvec(0 5)
\htext(-0.3 0.3){$1$} \htext(-0.7 0.75){$0$}
\htext(-0.5 1.5){$2$}
\htext(-0.5 2.25){$3$}
\htext(-0.5 2.75){$3$}
\htext(-0.5 3.5){$2$}
\htext(-0.7 4.75){$0$}
\esegment
\move(5 -7)
\bsegment
\move(0 0)\dtri \move(-1 0)\dtri
\move(-1 1)\lvec(-2 0)\lvec(-2 1)\lvec(0 1)\lvec(-1 0)\lvec(-1 4)\lvec(0 4)
\lvec(0 0)\lvec(-2 0)
\move(0 2)\lvec(-1 2)
\move(0 2.5)\lvec(-1 2.5)
\move(0 3)\lvec(-1 3)
\htext(-0.3 0.3){$1$} \htext(-0.7 0.75){$0$}
\htext(-1.3 0.3){$0$} \htext(-1.7 0.75){$1$}
\htext(-0.5 1.5){$2$}
\htext(-0.5 2.25){$3$}
\htext(-0.5 2.75){$3$}
\htext(-0.5 3.5){$2$}
\esegment
\move(-1.5 23.8)\avec(-1.5 21.7)\htext(-1.1 22.9){$0$}
\move(-1.5 20.2)\avec(-1.5 18.2)\htext(-1.1 19.4){$2$}
\move(-1.75 15.7)\avec(-3.4 13.7)\htext(-2.9 14.9){$1$}
\move(-1.25 15.7)\avec(0.3 14.2)\htext(-0.1 15.2){$3$}
\move(-3.5 11.2)\avec(-3.5 8.7)\htext(-3.9 10.1){$3$}
\move(0.1 11.2)\avec(-3.1 8.7)\htext(-1.9 10.2){$1$}
\move(0.5 11.2)\avec(0.5 9.2)\htext(0.9 10.3){$3$}
\move(-4.2 5.7)\avec(-5.8 2.7)\htext(-5.4 4.4){$2$}
\move(-3.8 5.7)\avec(-1.7 3.2)\htext(-2.4 4.7){$3$}
\move(0.3 5.7)\avec(-1.3 3.2)\htext(-0.8 4.7){$1$}
\move(0.7 5.7)\avec(2.3 4.2)\htext(1.8 5.2){$2$}
\move(-6.2 -0.3)\avec(-6.9 -4.3)\htext(-6.9 -2.1){$0$}
\move(-5.8 -0.3)\avec(-4.1 -4.3)\htext(-4.6 -2.1){$3$}
\move(-1.8 -0.3)\avec(-0.9 -3.8)\htext(-1.0 -1.9){$2$}
\move(2.3 -0.3)\avec(1.6 -1.8)\htext(1.7 -0.8){$0$}
\move(2.7 -0.3)\avec(4.4 -2.8)\htext(3.9 -1.3){$1$}
\vtext(-7.5 -8){$\cdots$}
\vtext(-4 -8){$\cdots$}
\vtext(-1 -8){$\cdots$}
\vtext(1.5 -8){$\cdots$}
\vtext(4 -8){$\cdots$}
\move(-10 -8.5)\move(7 25)
\end{texdraw}
\end{center}

\newpage

(g) The crystal $\Y(\La_3)$ for $B_3^{(1)}$ 

\vskip 3mm

\begin{center}
\begin{texdraw}
\fontsize{7}{7}\selectfont
\drawdim em
\setunitscale 1.75
%\htext(-5 20){\normalsize $\bullet$ $\Y(\La_3)$ for $B_3^{(1)}$}
\htext(-0.4 22.5){$Y_{\La_3}$}
\move(0 17)
\bsegment
\move(0 0)\lvec(-1 0)\lvec(-1 0.5)\lvec(0 0.5) \ifill f:0.7
\move(0 0)\lvec(0 1)\lvec(-1 1)\lvec(-1 0)\lvec(0 0)
\move(0 0.5)\lvec(-1 0.5)
\htext(-0.5 0.26){$3$}
\htext(-0.5 0.76){$3$}
\esegment
\move(0 12)
\bsegment
\move(0 0)\lvec(-1 0)\lvec(-1 0.5)\lvec(0 0.5) \ifill f:0.7
\move(0 0)\lvec(-1 0)\lvec(-1 2)\lvec(0 2)\lvec(0 0)
\move(0 0.5)\lvec(-1 0.5)
\move(0 1)\lvec(-1 1)
\htext(-0.5 0.26){$3$}
\htext(-0.5 0.76){$3$}
\htext(-0.5 1.5){$2$}
\esegment
\move(-4 6)
\bsegment
\move(0 0)\lvec(-1 0)\lvec(-1 0.5)\lvec(0 0.5) \ifill f:0.7
\move(0 2)\lvec(-1 2)\lvec(0 3)\lvec(0 0)\lvec(-1 0)\lvec(-1 3)\lvec(0 3)
\move(0 0.5)\lvec(-1 0.5)
\move(0 1)\lvec(-1 1)
\htext(-0.5 0.26){$3$}
\htext(-0.5 0.76){$3$}
\htext(-0.5 1.5){$2$}
\htext(-0.7 2.75){$0$}
\esegment
\move(0 6)
\bsegment
\move(0 0)\lvec(-1 0)\lvec(-1 0.5)\lvec(0 0.5) \ifill f:0.7
\move(0 2)\lvec(-1 2)\lvec(0 3)\lvec(0 0)\lvec(-1 0)\lvec(-1 3)\lvec(0 3)
\move(0 0.5)\lvec(-1 0.5)
\move(0 1)\lvec(-1 1)
\htext(-0.5 0.26){$3$}
\htext(-0.5 0.76){$3$}
\htext(-0.5 1.5){$2$}
\htext(-0.3 2.3){$1$}
\esegment
\move(4 6)
\bsegment
\move(0 0)\lvec(-2 0)\lvec(-2 0.5)\lvec(0 0.5) \ifill f:0.7
\move(-1 0)\lvec(-1 2)\lvec(0 2)\lvec(0 0)\lvec(-2 0)\lvec(-2 1)\lvec(0 1)
\move(0 0.5)\lvec(-2 0.5)
\htext(-0.5 0.26){$3$}
\htext(-1.5 0.26){$3$}
\htext(-0.5 0.76){$3$}
\htext(-1.5 0.76){$3$}
\htext(-0.5 1.5){$2$}
\esegment
\move(-4 0)
\bsegment
\move(0 0)\lvec(-1 0)\lvec(-1 0.5)\lvec(0 0.5) \ifill f:0.7
\move(0 2)\lvec(-1 2)\lvec(0 3)\lvec(0 0)\lvec(-1 0)\lvec(-1 3)\lvec(0 3)
\move(0 0.5)\lvec(-1 0.5)
\move(0 1)\lvec(-1 1)
\htext(-0.5 0.26){$3$}
\htext(-0.5 0.76){$3$}
\htext(-0.5 1.5){$2$}
\htext(-0.3 2.3){$1$}
\htext(-0.7 2.75){$0$}
\esegment
\move(0 0)
\bsegment
\move(0 0)\lvec(-2 0)\lvec(-2 0.5)\lvec(0 0.5) \ifill f:0.7
\move(-1 0)\lvec(-1 3)\lvec(0 3)\lvec(0 0)\lvec(-2 0)\lvec(-2 1)\lvec(0 1)
\move(0 0.5)\lvec(-2 0.5)
\move(0 2)\lvec(-1 2)\lvec(0 3)
\htext(-0.5 0.26){$3$}
\htext(-1.5 0.26){$3$}
\htext(-0.5 0.76){$3$}
\htext(-1.5 0.76){$3$}
\htext(-0.5 1.5){$2$}
\htext(-0.7 2.75){$0$}
\esegment
\move(4 0)
\bsegment
\move(0 0)\lvec(-2 0)\lvec(-2 0.5)\lvec(0 0.5) \ifill f:0.7
\move(-1 0)\lvec(-1 3)\lvec(0 3)\lvec(0 0)\lvec(-2 0)\lvec(-2 1)\lvec(0 1)
\move(0 0.5)\lvec(-2 0.5)
\move(0 2)\lvec(-1 2)\lvec(0 3)
\htext(-0.5 0.26){$3$}
\htext(-1.5 0.26){$3$}
\htext(-0.5 0.76){$3$}
\htext(-1.5 0.76){$3$}
\htext(-0.5 1.5){$2$}
\htext(-0.3 2.3){$1$}
\esegment
\move(-5 -6)
\bsegment
\move(0 0)\lvec(-1 0)\lvec(-1 0.5)\lvec(0 0.5) \ifill f:0.7
\move(0 2)\lvec(-1 2)\lvec(0 3)\lvec(0 0)\lvec(-1 0)\lvec(-1 3)\lvec(0 3)
\move(0 3)\lvec(0 4)\lvec(-1 4)\lvec(-1 3)
\move(0 0.5)\lvec(-1 0.5)
\move(0 1)\lvec(-1 1)
\htext(-0.5 0.26){$3$}
\htext(-0.5 0.76){$3$}
\htext(-0.5 1.5){$2$}
\htext(-0.3 2.3){$1$}
\htext(-0.7 2.75){$0$}
\htext(-0.5 3.5){$2$}
\esegment
\move(-2 -6)
\bsegment
\move(0 0)\lvec(-2 0)\lvec(-2 0.5)\lvec(0 0.5) \ifill f:0.7
\move(-1 0)\lvec(-1 3)\lvec(0 3)\lvec(0 0)\lvec(-2 0)\lvec(-2 1)\lvec(0 1)
\move(0 0.5)\lvec(-2 0.5)
\move(0 2)\lvec(-1 2)\lvec(0 3)
\htext(-0.5 0.26){$3$}
\htext(-1.5 0.26){$3$}
\htext(-0.5 0.76){$3$}
\htext(-1.5 0.76){$3$}
\htext(-0.5 1.5){$2$}
\htext(-0.3 2.3){$1$}
\htext(-0.7 2.75){$0$}
\esegment
\move(1 -6)
\bsegment
\move(0 0)\lvec(-2 0)\lvec(-2 0.5)\lvec(0 0.5) \ifill f:0.7
\move(-1 0)\lvec(-1 3)\lvec(0 3)\lvec(0 0)\lvec(-2 0)\lvec(-2 2)\lvec(0 2)
\move(0 0.5)\lvec(-2 0.5)
\move(0 1)\lvec(-2 1)
\move(-1 2)\lvec(0 3)
\htext(-0.5 0.26){$3$}
\htext(-1.5 0.26){$3$}
\htext(-0.5 0.76){$3$}
\htext(-1.5 0.76){$3$}
\htext(-0.5 1.5){$2$}
\htext(-1.5 1.5){$2$}
\htext(-0.7 2.75){$0$}
\esegment
\move(4 -6)
\bsegment
\move(0 0)\lvec(-2 0)\lvec(-2 0.5)\lvec(0 0.5) \ifill f:0.7
\move(-1 0)\lvec(-1 3)\lvec(0 3)\lvec(0 0)\lvec(-2 0)\lvec(-2 2)\lvec(0 2)
\move(0 0.5)\lvec(-2 0.5)
\move(0 1)\lvec(-2 1)
\move(-1 2)\lvec(0 3)
\htext(-0.5 0.26){$3$}
\htext(-1.5 0.26){$3$}
\htext(-0.5 0.76){$3$}
\htext(-1.5 0.76){$3$}
\htext(-0.5 1.5){$2$}
\htext(-1.5 1.5){$2$}
\htext(-0.3 2.3){$1$}
\esegment
\vtext(-5.5 -6.8){$\cdots$}
\vtext(-3 -6.8){$\cdots$}
\vtext(0 -6.8){$\cdots$}
\vtext(3 -6.8){$\cdots$}
\move(-0.5 22)\avec(-0.5 18.2)\htext(-0.1 20.2){$3$}
\move(-0.5 16.7)\avec(-0.5 14.2)\htext(-0.1 15.6){$2$}
\move(-0.9 11.7)\avec(-4.2 9.2)\htext(-3 10.6){$0$}
\move(-0.5 11.7)\avec(-0.5 9.2)\htext(-0.1 10.6){$1$}
\move(-0.1 11.7)\avec(3.2 8.2)\htext(1.8 10.6){$3$}
\move(-4.5 5.7)\avec(-4.5 3.2)\htext(-4.9 4.6){$1$}
\move(-4.1 5.7)\avec(-0.9 3.2)\htext(-3 5.4){$3$}
\move(-0.9 5.7)\avec(-4.1 3.2)\htext(-1.9 5.4){$0$}
\move(-0.1 5.7)\avec(3.1 3.2)\htext(1 5.4){$3$}
\move(2.9 5.7)\avec(-0.1 3.2)\htext(2 5.4){$0$}
\move(3.5 5.7)\avec(3.5 3.2)\htext(3.9 4.6){$1$}
\move(-4.7 -0.3)\avec(-5.3 -1.8)\htext(-5.3 -0.8){$2$}
\move(-4.3 -0.3)\avec(-2.8 -2.8)\htext(-3.1 -1.4){$3$}
\move(-1.2 -0.3)\avec(-2.5 -2.8)\htext(-2.2 -1.4){$1$}
\move(-0.8 -0.3)\avec(0.3 -2.8)\htext(-0.1 -0.9){$2$}
\move(2.8 -0.3)\avec(-2.2 -2.8)\htext(0.9 -0.9){$0$}
\move(3.5 -0.3)\avec(3.5 -2.8)\htext(3.9 -1.4){$2$}
\move(7 -7.3)\move(-8 22.9)
\end{texdraw}
\end{center}

\end{example}

%%%%%%%%%%%%%%%%%%%%%%%%%%%%%%%%%%%%%%%%%%%%%%%%%%%%%%%%%%%%%%%%%%%%%%
\newpage

%\bibliographystyle{amsplain}
%\bibliography{tetris}

\providecommand{\bysame}{\leavevmode\hbox to3em{\hrulefill}\thinspace}

\end{document}

%%%%%%%%%%%%%%%%%%%%%%%%%%%%%%%%%%%%%%%%%%%%%%%%%%%%%%%%%%%%%%%%%%%%%%%%%%

\vskip 3mm 
(g) $C_n^{(1)}$ ($n\geq2$)

\vskip 2mm
\begin{center}
\begin{texdraw}
\small
\drawdim em
\setunitscale 0.3
\move(0 0)
\lcir r:1
\move(1 0)\lvec(9 0)
\move(10 0)\lcir r:1
\move(11 0)\lvec(19 0)
\htext(23.6 0){$\cdots$}
\move(28 0)\lvec(36 0)
\move(37 0)\lcir r:1
\move(38 0.4)\lvec(46 0.4)
\move(38 -0.4)\lvec(46 -0.4)
\move(39.2 1.2)\lvec(38 0)\lvec(39.2 -1.2)
\move(47 0)\lcir r:1
\move(-1 0.4)\lvec(-9 0.4)
\move(-1 -0.4)\lvec(-9 -0.4)
\move(-10 0)\lcir r:1
\move(-2.2 1.2)\lvec(-1 0)\lvec(-2.2 -1.2)
\htext(-10 -3){$0$}
\htext(0 -3){$1$}
\htext(10 -3){$2$}
\htext(37 -3){$n-1$}
\htext(47 -3){$n$}
\move(-12 -4.6)\move(49 3)
\end{texdraw}
\end{center}

\vskip 2mm
\begin{equation*}
\begin{aligned}\mbox{}
& c = h_0 + h_1 +\cdots+h_{n-1}+h_n,\\
& \delta = \alpha_0 + 2\alpha_1 + \cdots +2\alpha_{n-1}+\alpha_n,\\
& \lambda = \Lambda_i \ \  (i\in I).
\end{aligned}
\end{equation*}

%crystal graph 

{\savebox{\tmppic}{\begin{texdraw}
\drawdim mm
\textref h:C v:C
\arrowheadsize l:2.4 w:1.1 \arrowheadtype t:F
\setunitscale 1
\htext(0 0){$b$}
\move(2 -0.5)\ravec(6 0)
\htext(10 0){$b'$}
\htext(4.5 0.8){$_i$}
\move(-1.1 -1.6)\move(11 1.5)
\end{texdraw}}%
\raisebox{-0.1\height}{\usebox{\tmppic}}%
}